\documentclass[11pt, twoside]{article}

\usepackage[text={375pt, 540pt}, centering, headheight=15pt]{geometry}
\usepackage{fancyhdr}
\usepackage{amssymb, amsmath}
\usepackage{graphicx}
\usepackage{hyperref}
\usepackage{breakurl}
\usepackage[labelsep=none]{caption}

\renewcommand{\headrulewidth}{0.5pt}

\hypersetup{pdfstartview=FitH, pdfpagemode=UseNone, urlbordercolor = {0 1 1}}

\newenvironment{definition}{\emph{Definition.}}{}
\newcounter{theorem}
\newenvironment{theorem}{\refstepcounter{theorem}\emph{Theorem \thetheorem.}}{}
\newenvironment{proof}{\emph{Proof.}}{$\square$}
\newcounter{lemma}
\newenvironment{lemma}{\refstepcounter{lemma}\emph{Lemma \thelemma.}}{}
\newcounter{corollary}
\newenvironment{corollary}{\refstepcounter{corollary}\emph{Corollary \thecorollary.}}{}
\newcounter{question}
\newenvironment{question}{\refstepcounter{question}\emph{Question \thequestion.}}{}

\allowdisplaybreaks

\def \tthbar {\texttt{\raisebox{0.75ex}{\scalebox{1}[0.75]{-}}\kern-0.5em h}}

\begin{document}

\vspace*{\bigskipamount}

\noindent{\LARGE \bfseries The Second Leaper Theorem}\bigskip

\noindent{\scshape Nikolai Beluhov}\bigskip

\begin{center} \parbox{336pt}{\footnotesize \emph{Abstract.} A $(p, q)$-leaper is a fairy chess piece that, from a square $a$, can move to any of the squares $a + (\pm p, \pm q)$ or $a + (\pm q, \pm p)$. Let $L$ be a $(p, q)$-leaper with $p + q$ odd and $C$ a cycle of $L$ within a $(p + q) \times (p + q)$ chessboard. We show that there exists a second leaper $M$, distinct from $L$, such that a Hamiltonian cycle $D$ of $M$ exists over the squares of $C$. We give descriptions of $C$ and $M$ in terms of continued fractions. We introduce the notion of a direction graph, roughly a leaper graph from which all information has been abstracted away save for the directions of the moves, and we study $C$ and $D$ in terms of direction graphs. We introduce the notion of a dual generalized chessboard, a generalized chessboard $B$ of more than one square such that the leaper graph of a leaper $L$ over $B$ is connected and isomorphic to the leaper graph of a second leaper $M$, distinct from $L$, over $B$, and we give constructions for dual generalized chessboards.} \end{center}

\section{Preliminaries}

\emph{Fairy chess} is the study of chess problems featuring unusual boards, pieces, or stipulations.

\medskip

\begin{definition} A \emph{regular chessboard} is a rectangular grid of unit squares. A \emph{generalized chessboard} is a set of unit squares in the plane, with sides parallel to the coordinate axes, obtained from each other by means of integer translations. \end{definition}

\medskip

For instance, every polyomino is a generalized chessboard.

The \emph{infinite chessboard} is the one obtained by dissecting all of the plane into unit squares by means of two pencils of parallel lines.

We reserve the term \emph{square} for a unit square regarded as a part of a chessboard. We refer to both regular and generalized chessboards as \emph{boards} for short.

Consider any object $O$ that consists of squares, possibly together with some structure imposed on them. For instance, $O$ may be a square, a board, or a graph whose vertices are squares. We write $O + v$ for the copy of $O$ under a translation $v$.

\medskip

\begin{definition} Let $p$ and $q$ be nonnegative integers, at least one of them positive. A \emph{$(p, q)$-leaper} $L$ is a fairy chess piece that, from a square $a$, can move to any of the squares $a + (\pm p, \pm q)$ or $a + (\pm q, \pm p)$. \end{definition}

\medskip

For instance, in orthodox chess the knight is a $(1, 2)$-leaper and the king is a combination of a $(0, 1)$-leaper and a $(1, 1)$-leaper.

We refer to $p$ and $q$ as the \emph{proportions} of $L$, and to translations of the form $(\pm p, \pm q)$ or $(\pm q, \pm p)$ as $L$-translations.

\medskip

\begin{definition} The \emph{leaper graph} of a leaper $L$ over a board $B$ is the graph whose vertices are the squares of $B$ and whose edges join the pairs of squares of $B$ that are joined by a move of $L$. \end{definition}

\medskip

\begin{definition} A leaper $L$ is \emph{free} over a board $B$ if the leaper graph of $L$ over $B$ is connected. \end{definition}

\medskip

\begin{definition} An \emph{open tour} of a leaper $L$ over a board $B$ is a Hamiltonian path in the leaper graph of $L$ over $B$. A \emph{closed tour} of $L$ over $B$ is a Hamiltonian cycle in the leaper graph of $L$ over $B$. \end{definition}

\medskip

The problem of constructing knight tours over regular chessboards dates back to at least the ninth century. For a detailed historical overview, see the \emph{History} sections of \cite{J2}.

The study of leaper tours beyond the knight appears to have commenced sometime in the late nineteenth and early twentieth century. The concept of a free leaper was introduced in \cite{JW} by George Jelliss and Theophilus Willcocks. General questions about leaper graphs were raised in \cite{JW} and George Jelliss' \cite{J1}, the latter also establishing a number of general properties of leaper graphs. In \cite{K}, Donald Knuth completely solved the question of whether a given leaper is free over a given regular chessboard, and studied the question of whether there exists a closed tour of a given leaper over a given regular chessboard. For a detailed historical overview, see the \emph{Leapers at Large} section of \cite{J2}.

\medskip

\begin{definition} A $(p, q)$-leaper $L$ is \emph{orthogonal} if $p = 0$ or $q = 0$, \emph{diagonal} if $p = q$, and \emph{skew} if $p \neq 0$, $q \neq 0$, and $p \neq q$. \end{definition}

\medskip

There exist four possible directions for the moves of an orthogonal leaper (east, north, west, and south), four possible directions for the moves of a diagonal leaper (northeast, northwest, southwest, and southeast), and eight possible directions for the moves of a skew leaper (east-northeast, north-northeast, \ldots, east-southeast).

\medskip

\begin{definition} A $(p, q)$-leaper $L$ is \emph{basic} if $p + q$ is odd and $p$ and $q$ are relatively prime. \end{definition}

\medskip

Equivalently, a $(p, q)$-leaper $L$ is basic if $p - q$ and $p + q$ are relatively prime.

All basic leapers apart from the $(0, 1)$-leaper are skew leapers.

There are a number of reasons to single out the class of basic leapers.

\begin{figure}[ht!]\hspace*{\fill}\includegraphics[scale=0.66]{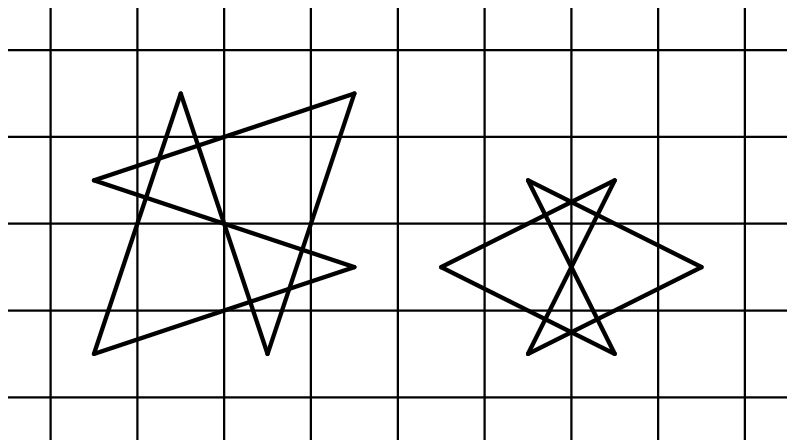}\hspace*{\fill}\caption{}\label{1312}\end{figure}

Firstly, let $L$ be an arbitrary leaper. Then there exists a unique basic leaper $L'$ such that every connected leaper graph of $L$ is a scaled and rotated copy of a leaper graph of $L'$. For instance, when $L$ is a $(1, 3)$-leaper, $L'$ is a $(1, 2)$-leaper, as in Figure \ref{1312}.

More precisely, let $L$ be a $(p, q)$-leaper, $d$ the greatest common divisor of $p$ and $q$, and $p = dp'$ and $q = dq'$. If $p' + q'$ is odd, then $L'$ is a $(p', q')$-leaper, the scaling factor equals $d$, and without loss of generality the angle of rotation equals $0^\circ$. If $p' + q'$ is even, then $L'$ is a $(\frac{1}{2}|p' - q'|, \frac{1}{2}(p' + q'))$ leaper, the scaling factor equals $\sqrt{2}d$, and without loss of generality the angle of rotation equals $45^\circ$, as in Figure \ref{1312}.

Therefore, if a problem is only concerned with the intrinsic properties of leaper graphs, it suffices to study basic leapers.

Secondly, a leaper $L$ is free over the infinite chessboard (or, equivalently, over all sufficiently large regular chessboards) if and only if it is basic.

\section{The Second Leaper Theorem}

Consider a $(p, q)$-leaper $L$, $p < q$, over a $(p + q) \times (p + q)$ chessboard $B$. The central $(q - p) \times (q - p)$ squares of $B$ are of degree zero in the corresponding leaper graph $G$, and all remaining squares are of degree two. Therefore, $G$ consists of $(q - p)^2$ isolated vertices and a number of disjoint cycles.

The following theorem was discovered by the author in January 2006.

\medskip

\begin{theorem} Let $L$ be a $(p, q)$-leaper with $p + q$ odd and $C$ a cycle of $L$ within a $(p + q) \times (p + q)$ chessboard. Then there exists a second leaper $M$, distinct from $L$, such that a Hamiltonian cycle of $M$ exists over the squares of $C$. \label{sl} \end{theorem}

\medskip

A number of examples are in order.

\begin{figure}[ht!]\hspace*{\fill}\includegraphics[scale=0.66]{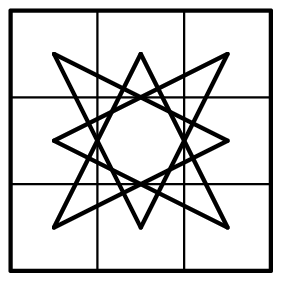}\hfill\includegraphics[scale=0.66]{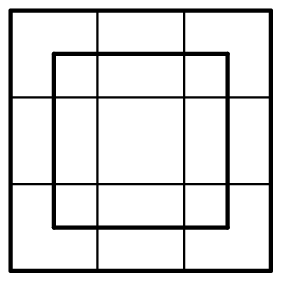}\hspace*{\fill}\caption{}\label{1233}\end{figure}

\begin{figure}[ht!]\hspace*{\fill}\includegraphics[scale=0.66]{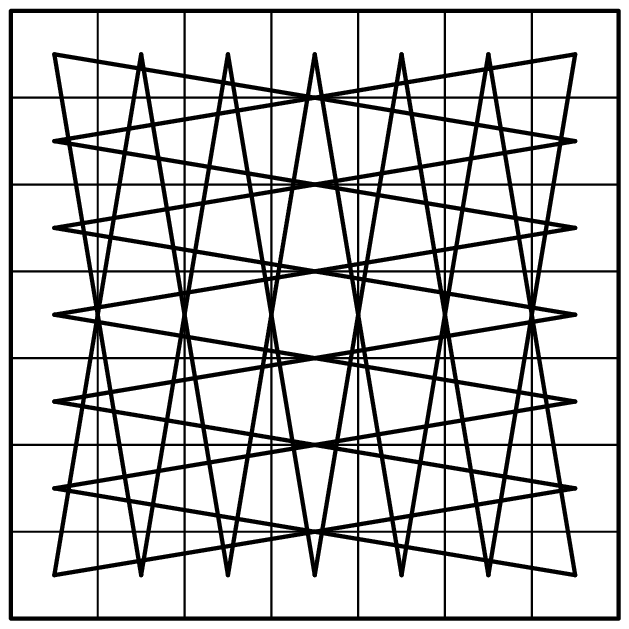}\hfill\includegraphics[scale=0.66]{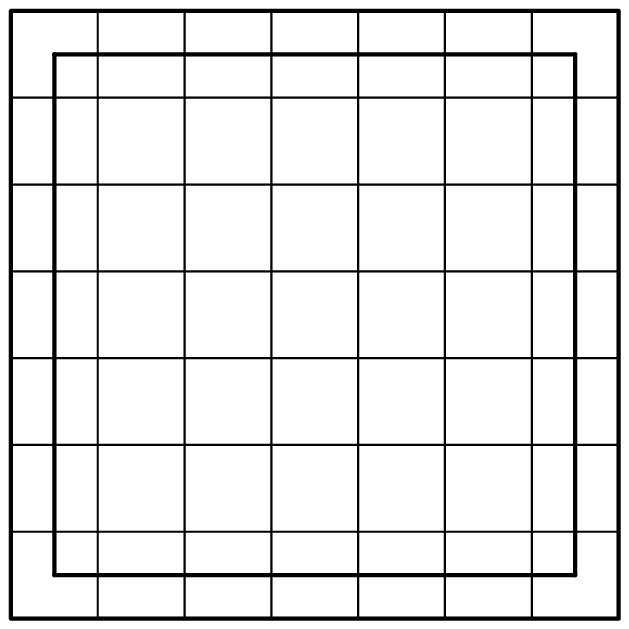}\hspace*{\fill}\caption{}\label{1677}\end{figure}

When $p = 1$ and $q = 2$, $G$ contains a single eight-square cycle $C$ which is also toured by a $(0, 1)$-leaper (Figure \ref{1233}). In general, when $p = 1$ and $q = 2k$, $G$ contains a single $4k$-square cycle which is also toured by a $(0, 1)$-leaper (Figure \ref{1677}).

\begin{figure}[t!]\hspace*{\fill}\includegraphics[scale=0.66]{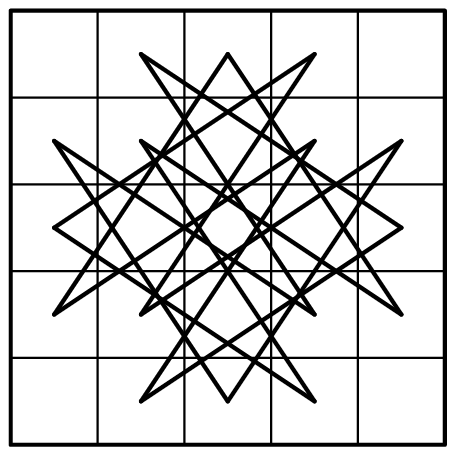}\hfill\includegraphics[scale=0.66]{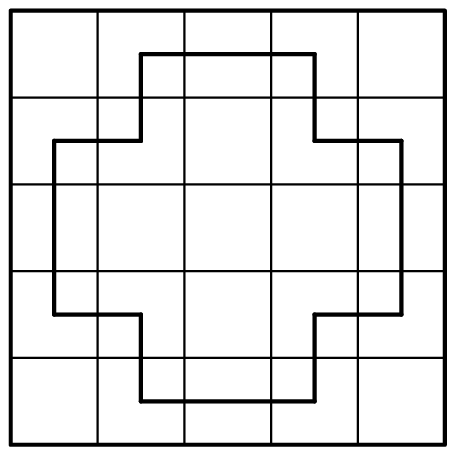}\hspace*{\fill}\newline\bigskip\newline\hspace*{\fill}\includegraphics[scale=0.66]{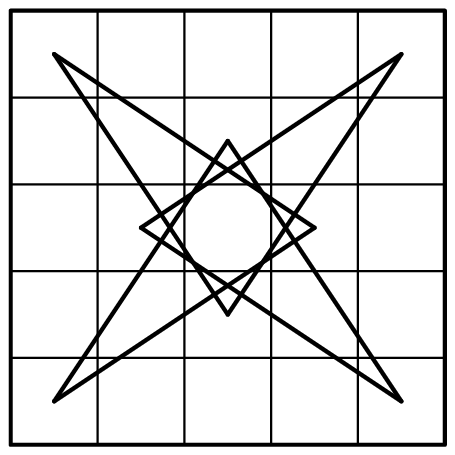}\hfill\includegraphics[scale=0.66]{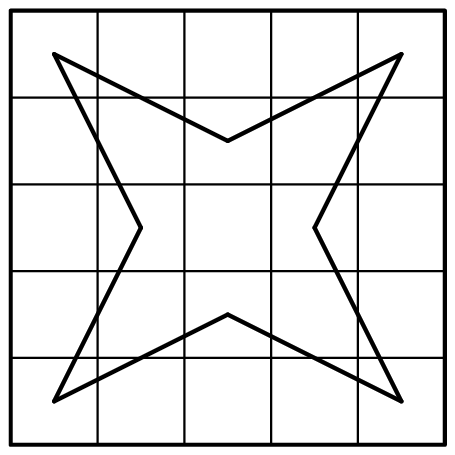}\hspace*{\fill}\caption{}\label{2355}\end{figure}

\begin{figure}[p!]\hspace*{\fill}\includegraphics[scale=0.66]{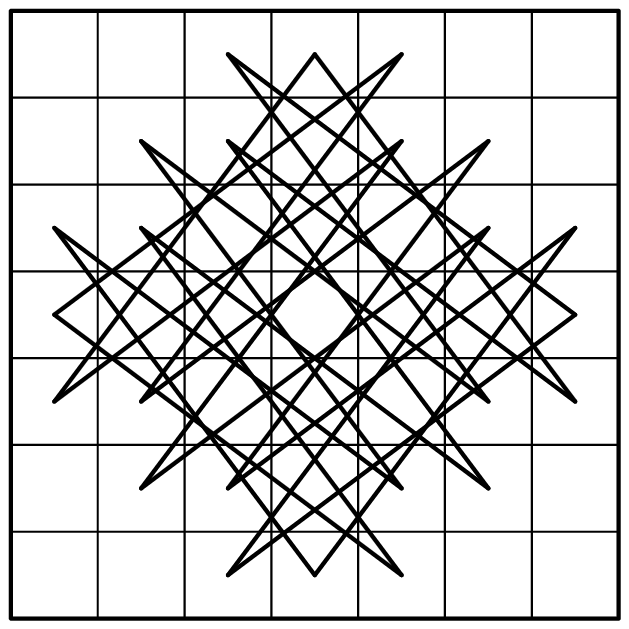}\hfill\includegraphics[scale=0.66]{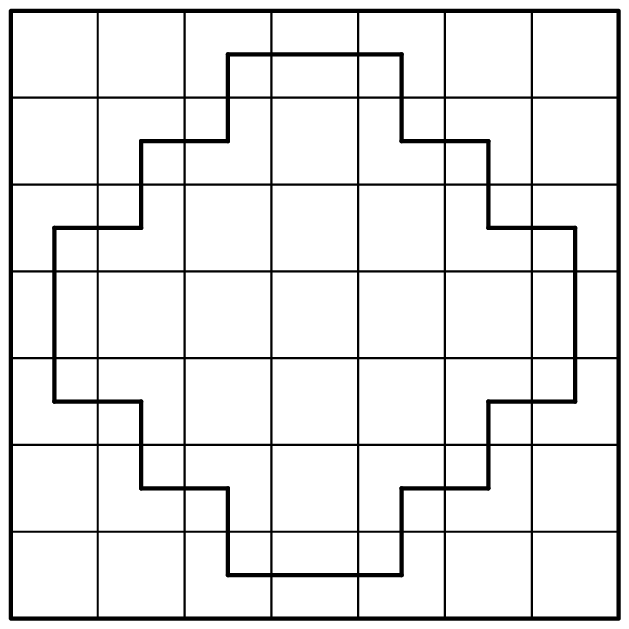}\hspace*{\fill}\newline\bigskip\newline\hspace*{\fill}\includegraphics[scale=0.66]{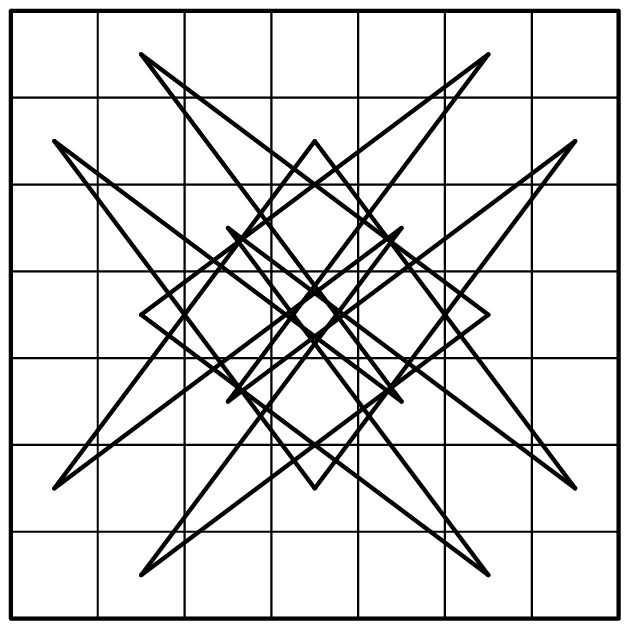}\hfill\includegraphics[scale=0.66]{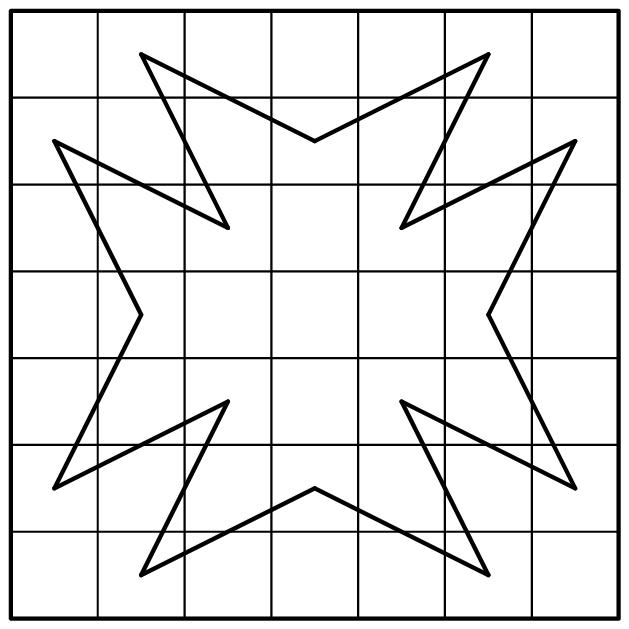}\hspace*{\fill}\newline\bigskip\newline\hspace*{\fill}\includegraphics[scale=0.66]{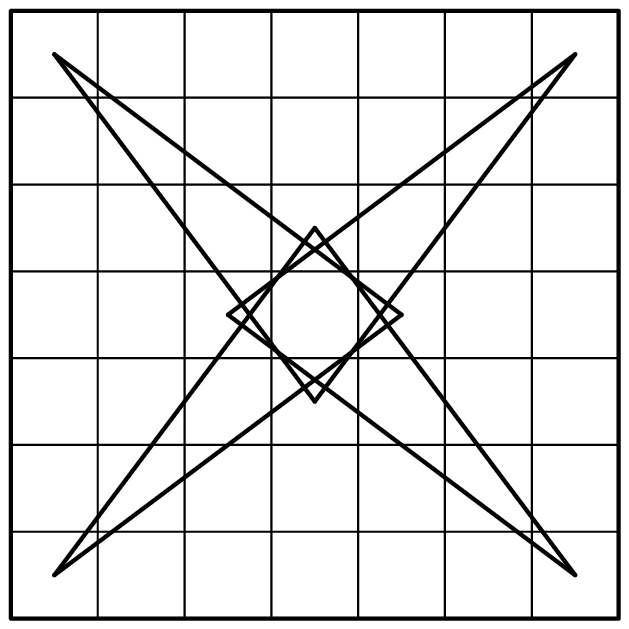}\hfill\includegraphics[scale=0.66]{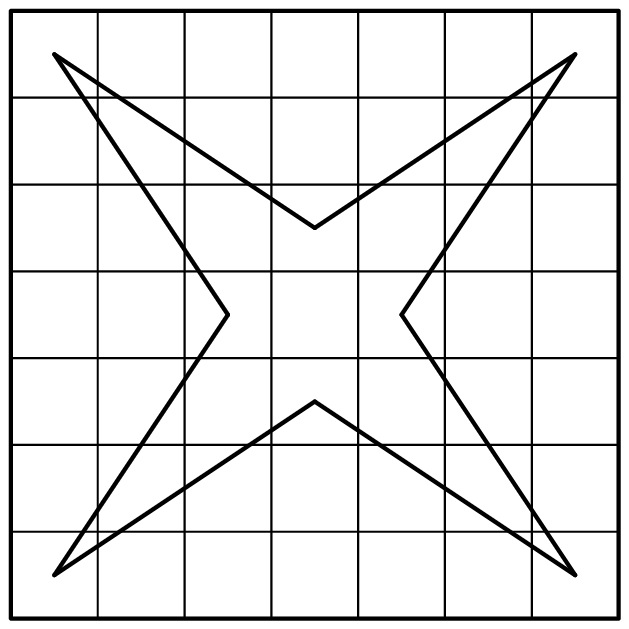}\hspace*{\fill}\caption{}\label{3477}\end{figure}

When $p = 2$ and $q = 3$, $G$ contains one 16-square cycle which is also toured by a $(0, 1)$-leaper and one eight-square cycle which is also toured by a $(1, 2)$-leaper (Figure \ref{2355}). In general, when $q = p + 1$, $G$ contains $p$ cycles, one toured by each $(r, r + 1)$-leaper, $r = 0, 1, \ldots, p - 1$ (Figure \ref{3477}).

\begin{figure}[t!]\hspace*{\fill}\includegraphics[scale=0.66]{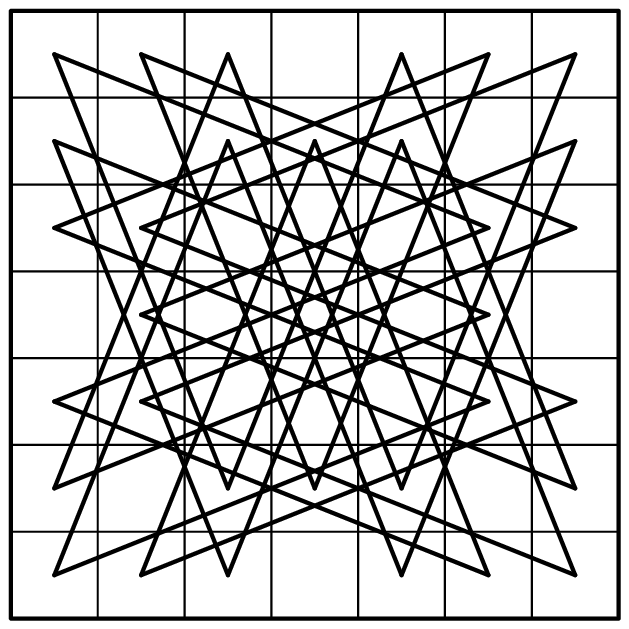}\hfill\includegraphics[scale=0.66]{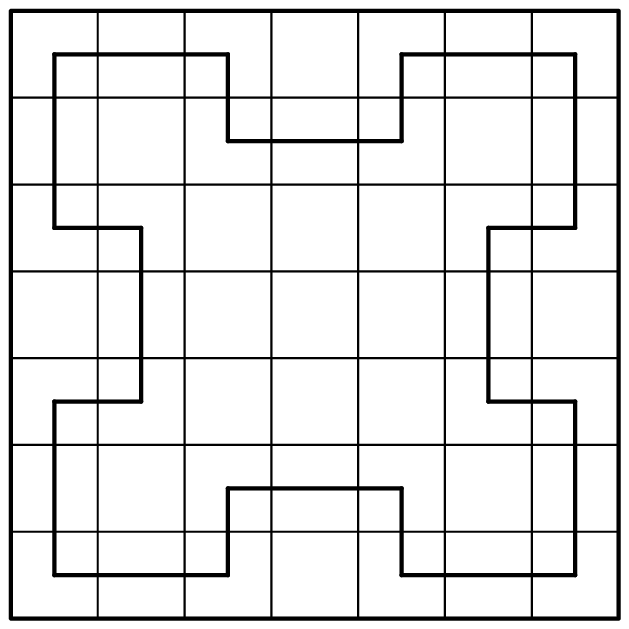}\hspace*{\fill}\newline\bigskip\newline\hspace*{\fill}\includegraphics[scale=0.66]{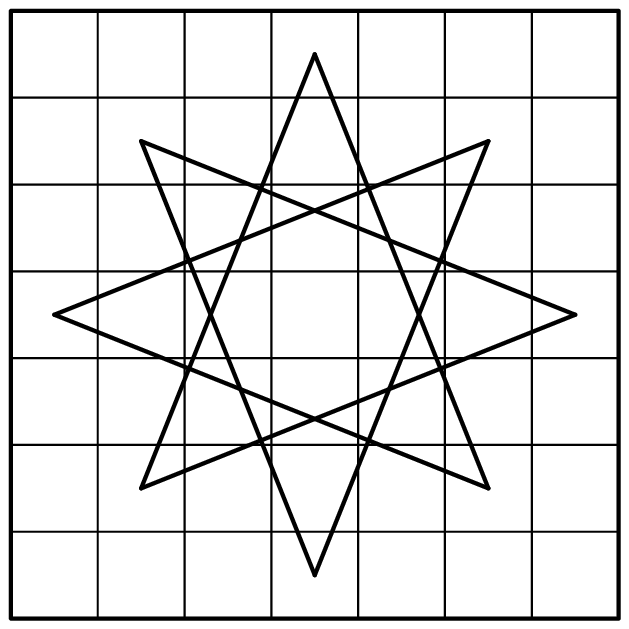}\hfill\includegraphics[scale=0.66]{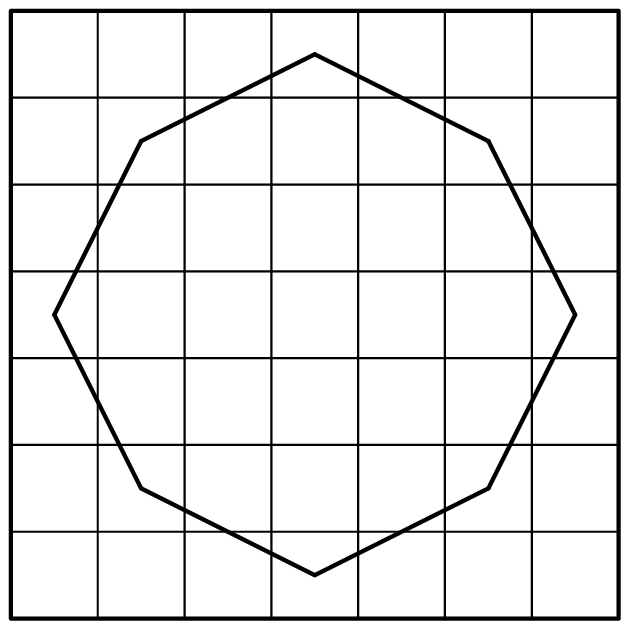}\hspace*{\fill}\caption{}\label{2577}\end{figure}

When $p = 2$ and $q = 5$, $G$ contains one 32-square cycle which is also toured by a $(0, 1)$-leaper and one eight-square cycle which is also toured by a $(1, 2)$-leaper (Figure \ref{2577}).

\medskip

\begin{proof} We begin by introducing several useful notions.

\medskip

\begin{definition} An \emph{$(m, n)$-frame} $F$, $m < n$, is an $(m + n) \times (m + n)$ chessboard with an $(n - m) \times (n - m)$ hole in the center. The \emph{sections} of $F$, $F_\texttt{E}$, $F_\texttt{NE}$, \ldots, $F_\texttt{SE}$ as in Figure \ref{framesect}, are the eight rectangular subboards that $F$ is dissected into when the sides of the central hole are extended. \end{definition}

\begin{figure}[ht!]\hspace*{\fill}\includegraphics{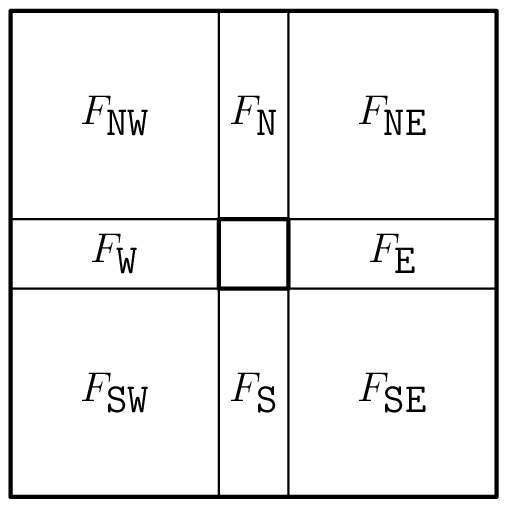}\hspace*{\fill}\caption{}\label{framesect}\end{figure}

\medskip

Given an $(m, n)$-frame $F$, we write $L_F$ for the corresponding $(m, n)$-leaper and $G_F$ for the leaper graph of $L_F$ over $F$.

The \texttt{E}, \texttt{NE}, \ldots, \texttt{SE} stand for the eight directions east, northeast, \ldots, southeast. Given a direction $i$ out of \texttt{E}, \texttt{NE}, \ldots, \texttt{SE} and an integer $k$, we write $-i$ for the direction opposite $i$ and $i + k$ for the direction $k$ steps counterclockwise from $i$. For instance, $-\texttt{E} = \texttt{W}$, $\texttt{NE} + 1 = \texttt{N}$, for all directions $i$, $-i = i + 4$, and, for all directions $i$ and integers $k$, $-(i + k) = -i + k$.

Given a path $w$ through a number of squares, we write $-w$ for the path obtained from $w$ by traversing it in the opposite direction.

We proceed to shed some light onto the structure of a second leaper's Hamiltonian cycle.

\medskip

\begin{definition} A cycle $D$ of a leaper $M$ within a frame $F$ is \emph{proper} if $D$ is the concatenation of eight disjoint nonempty paths, \[D = a^\texttt{E}a^\texttt{NE}\ldots a^\texttt{SE},\] such that, for all directions $i$, $a^i$ lies within $F_i$ and the $L_F$-translation that maps $F_i \cup F_{i + 1}$ onto $F_{-(i + 1)} \cup F_{-i}$ also maps $a^ia^{i + 1}$ onto $-a^{-(i + 1)}{-a^{-i}}$. \end{definition}

\medskip

Let us write this out in more detail. In a proper cycle, we have \[ a^\texttt{E}a^\texttt{NE} + (-n, -m) = -a^\texttt{SW}{-a^\texttt{W}}, \] \[ a^\texttt{NE}a^\texttt{N} + (-m, -n) = -a^\texttt{S}{-a^\texttt{SW}}, \] \[ a^\texttt{N}a^\texttt{NW} + (m, -n) = -a^\texttt{SE}{-a^\texttt{S}}, \] and \[ a^\texttt{NW}a^\texttt{W} + (n, -m) = -a^\texttt{E}{-a^\texttt{SE}}.\]

Eventually we will see that, in the statement of the theorem, ``a Hamiltonian cycle of $M$'' can be replaced with ``a proper Hamiltonian cycle of $M$''.

From here on, the plan of the proof is as follows. First we define three transformations, $f$, $g$, and $h$, that \emph{lift} smaller frames to larger ones. Given a subset $S$ of the squares of a smaller frame $F$, each lifting transformation constructs from it a subset $T$ of the squares of a larger frame $H$.

This is done in such a way that if $S$ is the vertex set of an $L_F$-cycle within $F$ then $T$ is the vertex set of an $L_H$-cycle within $H$ and if $S$ is the vertex set of a proper cycle of a leaper $M$ within $F$ then $T$ is the vertex set of a proper cycle of $M$ within $H$.

The proof is then completed by induction on $p + q$.

We begin with $f$.

\medskip

\begin{definition} Let $F$ be an $(m, n)$-frame, $m < n$, and $H$ an $(m, 2m + n)$-frame. Place $F$ and $H$ so that their centers coincide. (Or, equivalently, so that the outer contour of $F$ coincides with the inner contour of $H$.) We define the transformation $f$, lifting $F$ to $H$, as follows.

An \emph{$f$-translation} is any of the eight translations corresponding to a move of a $(0, m + n)$-leaper or an $(m, m)$-leaper. Given a direction $i$, we write $v^f_i$ for the $f$-translation pointing in direction $i$. For instance, $v^f_\texttt{E} = (m + n, 0)$ and $v^f_\texttt{NE} = (m, m)$.

Let $S$ be a subset of the squares of $F$. Then $f(S)$ is a subset of the squares of $H$, the disjoint union of \[ (S + v^f_i) \cap H_i \] when $i$ ranges over \texttt{E}, \texttt{NE}, \ldots, \texttt{SE}. \end{definition}

\medskip

Equivalently, $f(S)$ is the disjoint union of a number of translation copies of the intersections of $S$ with the sections of $F$ as in Figure \ref{fdef}. Each subboard $I$ of $H$ labeled $F_i$ in Figure \ref{fdef} contains a copy of $S \cap F_i$ under the $f$-translation that maps $F_i$ onto $I$.

\begin{figure}[ht!]\hspace*{\fill}\includegraphics{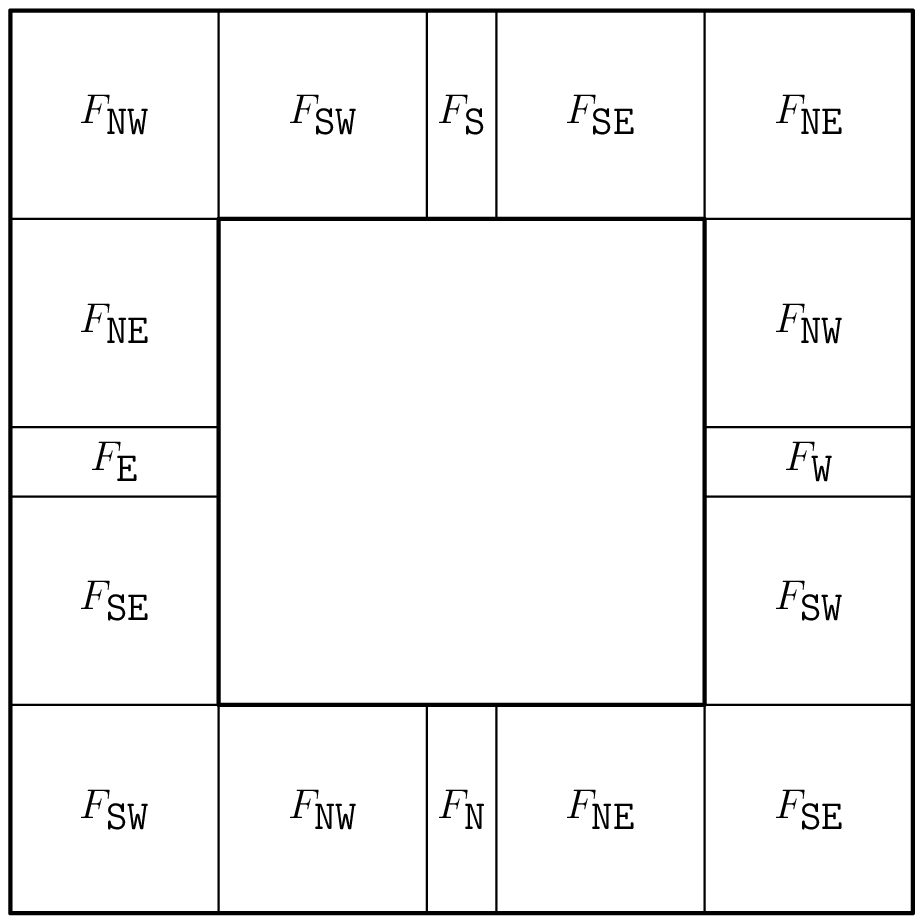}\hspace*{\fill}\caption{}\label{fdef}\end{figure}

We proceed to establish the two key properties of $f$.

\medskip

\begin{lemma} Let $S$ be the vertex set of an $L_F$-cycle within $F$. Then $f(S)$ is the vertex set of an $L_H$-cycle within $H$. \label{f1} \end{lemma}

\medskip

\begin{proof} Let $S$ be the vertex set of the $L_F$-cycle $C$.

Replace each square $a$ in $C$ with a path within $H$ through the images of $a$ under $f$ as follows.

\medskip

\emph{Case 1.} $a$ belongs to $F_\texttt{E}$. Replace $a$ with the single-square path \[ a + v^f_\texttt{W}. \]

\medskip

\emph{Case 2.} $a$ belongs to $F_\texttt{NE}$. Then the two squares adjacent to $a$ in $C$ are $a + (-n, -m)$ and $a + (-m, -n)$. Replace $a$ with the three-square path \[ (a + v^f_\texttt{W})(a + v^f_\texttt{NE})(a + v^f_\texttt{S}), \] traversed in such a direction that the endpoint $a + v^f_\texttt{W}$ is on the side of $a + (-n, -m)$ and the endpoint $a + v^f_\texttt{S}$ is on the side of $a + (-m, -n)$.

\medskip

All other cases are obtained from Cases 1 and 2 by symmetry.

The concatenation of the resulting sequence of paths is an $L_H$-cycle of vertex set $f(S)$. \end{proof}

\medskip

\begin{lemma} Let $S$ be the vertex set of a proper cycle of a leaper $M$ within $F$. Then $f(S)$ is the vertex set of a proper cycle of $M$ within $H$. \label{f2} \end{lemma}

\medskip

\begin{proof} Let $S$ be the vertex set of the proper cycle \[D = a^\texttt{E}a^\texttt{NE}\ldots a^\texttt{SE}\] of $M$ within $F$.

Let \[b^\texttt{E} = -a^\texttt{SW}{-a^\texttt{W}}{-a^\texttt{NW}} + v^f_\texttt{E}\] and define $b^\texttt{N}$, $b^\texttt{W}$, and $b^\texttt{S}$ analogously. Also let \[b^\texttt{NE} = a^\texttt{NE} + v^f_\texttt{NE}\] and define $b^\texttt{NW}$, $b^\texttt{SW}$, and $b^\texttt{SE}$ analogously. We proceed to show that \[ b^\texttt{E}b^\texttt{NE}\ldots b^\texttt{SE} \] is a proper cycle of $M$ of vertex set $f(S)$.

We need to show that, for all directions $i$, the final square of $b^i$ is joined by an $M$-move to the opening square of $b^{i + 1}$. We establish the case $i = \texttt{E}$, and all other cases are obtained from it by symmetry.

Since $D$ is a proper cycle, \[ -a^\texttt{W}{-a^\texttt{NW}} + (n, -m) = a^\texttt{SE}a^\texttt{E}. \]

It follows that the suffix $-a^\texttt{W}{-a^\texttt{NW}} + v^f_\texttt{E}$ of $b^\texttt{E}$ coincides with $a^\texttt{SE}a^\texttt{E} + v^f_\texttt{NE}$. Therefore, the final square of $b^\texttt{E}$ and the opening square of $b^\texttt{NE}$ are in the same relative position as the final square of $a^\texttt{E}$ and the opening square of $a^\texttt{NE}$.

We are left to show that, for all directions $i$, the $L_H$-translation that maps $H_i \cup H_{i + 1}$ onto $H_{-(i + 1)} \cup H_{-i}$ also maps $b^ib^{i + 1}$ onto $-b^{-(i + 1)}{-b^{-i}}$. We establish the case $i = \texttt{E}$, and all other cases are obtained from it by symmetry.

We have \[ b^\texttt{E}b^\texttt{NE} = (-a^\texttt{SW} + v^f_\texttt{E})(-a^\texttt{W}{-a^\texttt{NW}} + v^f_\texttt{E})(a^\texttt{NE} + v^f_\texttt{NE}) \] and \[ -b^\texttt{SW}{-b^\texttt{W}} = (-a^\texttt{SW} + v^f_\texttt{SW})(a^\texttt{SE}{a^\texttt{E}} + v^f_\texttt{W})(a^\texttt{NE} + v^f_\texttt{W}). \]

Since $D$ is a proper cycle, \[-a^\texttt{W}{-a^\texttt{NW}} + (n, -m) = a^\texttt{SE}{a^\texttt{E}}.\]

Therefore, a $(-2m - n, -m)$ translation maps $b^\texttt{E}b^\texttt{NE}$ onto $-b^\texttt{SW}{-b^\texttt{W}}$. \end{proof}

\medskip

Before we continue to $g$ and $h$, we need to introduce one more species of subdivision of a frame.

\medskip

\begin{definition} Let $m < n$, $3m \ge n$, and $F$ be an $(m, n)$-frame. The \emph{shell} of $F$, $F^-$, is the union of eight equal, symmetrically placed square subboards of $F$ defined as follows.

When $2m \ge n$, $F^-$ is the union of eight subboards of $F$ of size $(n - m) \times (n - m)$, one in each corner and four adjacent by side to the central hole, as in Figure \ref{shelldefa}.

When $2m \le n \le 3m$, $F^-$ is the union of eight subboards of $F$ of size $(3m - n) \times (3m - n)$, one in the middle of each outer side and four adjacent by corner to the central hole, as in Figure \ref{shelldefb}.

In both cases, the \emph{core} of $F$, $F^+$, is $F \setminus F^-$. \end{definition}

\begin{figure}[p!]\hspace*{\fill}\includegraphics{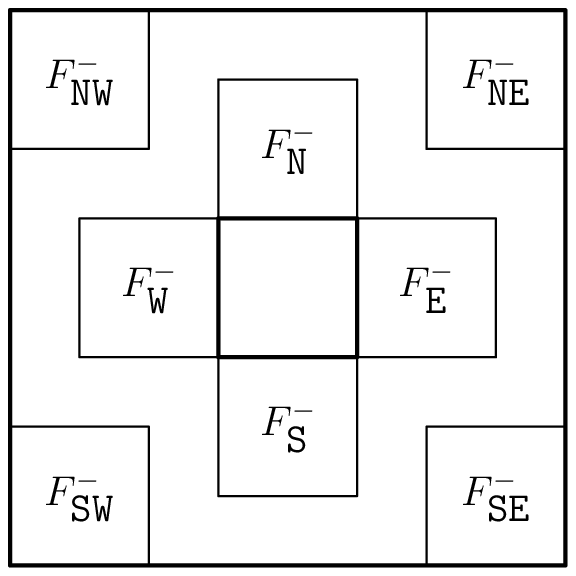}\hspace*{\fill}\caption{}\label{shelldefa}\end{figure}

\begin{figure}[p!]\hspace*{\fill}\includegraphics{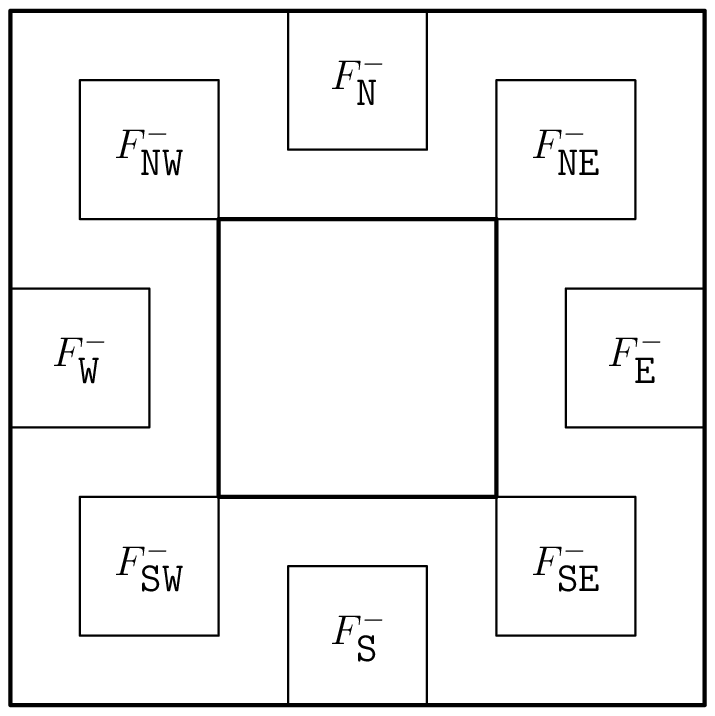}\hspace*{\fill}\caption{}\label{shelldefb}\end{figure}

\medskip

When $2m = n$, both parts of the definition give the same shell which coincides with the complete frame.

Let, for all directions $i$, $F^-_i$ be the intersection of the shell of $F$ and the section $F_i$ of $F$, as in Figures \ref{shelldefa} and \ref{shelldefb}. Then, for all $i$, an $L_F$-translation maps $F^-_i$ onto $F^-_{i + 3}$. Thus the shell of $F$ is the union of the vertex sets of a number of disjoint eight-square $L_F$-cycles and the core of $F$ is the union of the vertex sets of all other $L_F$-cycles within $F$.

We go on to $g$ and $h$.

\medskip

\begin{definition} Let $F$ be an $(m, n)$-frame, $m < n$, and $H$ an $(n, 2n - m)$-frame. Place $F$ and $H$ so that their centers coincide. (Or, equivalently, so that their inner contours coincide.) We define the transformation $g$, lifting $F$ to $H$, as follows.

A \emph{$g$-translation} is any of the eight translations corresponding to a move of a $(0, n - m)$-leaper or an $(n, n)$-leaper. Given a direction $i$, we write $v^g_i$ for the $g$-translation pointing in direction $i$. For instance, $v^g_\texttt{E} = (n - m, 0)$ and $v^f_\texttt{NE} = (n, n)$.

Let $S$ be a subset of the squares of $F$. Then $g(S)$ is a subset of the core of $H$, the disjoint union of \[ (S + v^g_i) \cap H_i \] when $i$ ranges over \texttt{E}, \texttt{NE}, \ldots, \texttt{SE}. \end{definition}

\medskip

Equivalently, $g(S)$ is the disjoint union of a number of translation copies of the intersections of $S$ with the sections of $F$ as in Figure \ref{gdef}.

\begin{figure}[ht!]\hspace*{\fill}\includegraphics{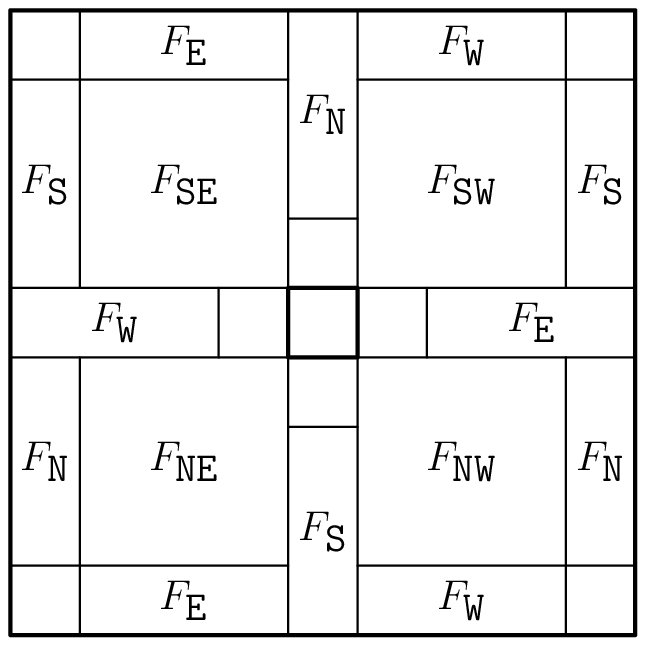}\hspace*{\fill}\caption{}\label{gdef}\end{figure}

\medskip

\begin{lemma} Let $S$ be the vertex set of an $L_F$-cycle within $F$. Then $g(S)$ is the vertex set of an $L_H$-cycle within $H$. \label{g1} \end{lemma}

\medskip

\begin{proof} Analogous to the proof of Lemma \ref{f1}.

Let $S$ be the vertex set of the $L_F$-cycle $C$ and $a$ a square in $C$. The transformation rules for $a$ are as follows.

\medskip

\emph{Case 1.} $a$ belongs to $F_\texttt{E}$. Then the two squares adjacent to $a$ in $C$ are $a + (-n, m)$ and $a + (-n, -m)$. Replace $a$ with the three-square path \[ (a + v^f_\texttt{NW})(a + v^f_\texttt{E})(a + v^f_\texttt{SW}), \] traversed in such a direction that the endpoint $a + v^f_\texttt{NW}$ is on the side of $a + (-n, m)$ and the endpoint $a + v^f_\texttt{SW}$ is on the side of $a + (-n, -m)$.

\medskip

\emph{Case 2.} $a$ belongs to $F_\texttt{NE}$. Replace $a$ with the single-square path \[ a + v^f_\texttt{SW}. \]

\medskip

All other cases are obtained from Cases 1 and 2 by symmetry. \end{proof}

\medskip

\begin{lemma} Let $S$ be the vertex set of a proper cycle of a leaper $M$ within $F$. Then $g(S)$ is the vertex set of a proper cycle of $M$ within $H$. \label{g2} \end{lemma}

\medskip

The proof is analogous to the proof of Lemma \ref{f2}.

\medskip

\begin{definition} Let $F$ be an $(m, n)$-frame, $m < n$, and $H$ an $(n, m + 2n)$-frame. Place $F$ and $H$ so that their centers coincide. (Or, equivalently, so that the outer contour of $F$ coincides with the inner contour of $H$.) We define the transformation $h$, lifting $F$ to $H$, as follows.

An \emph{$h$-translation} is any of the eight translations corresponding to a move of a $(0, m + n)$-leaper or an $(n, n)$-leaper. Given a direction $i$, we write $v^h_i$ for the $h$-translation pointing in direction $i$. For instance, $v^h_\texttt{E} = (m + n, 0)$ and $v^h_\texttt{NE} = (n, n)$.

Let $S$ be a subset of the squares of $F$. Then $h(S)$ is a subset of the core of $H$, the disjoint union of \[ (S + v^h_i) \cap H_i \] when $i$ ranges over \texttt{E}, \texttt{NE}, \ldots, \texttt{SE}. \end{definition}

\medskip

Equivalently, $h(S)$ is the disjoint union of a number of translation copies of the intersections of $S$ with the sections of $F$ as in Figure \ref{hdef}.

\begin{figure}[ht!]\hspace*{\fill}\includegraphics{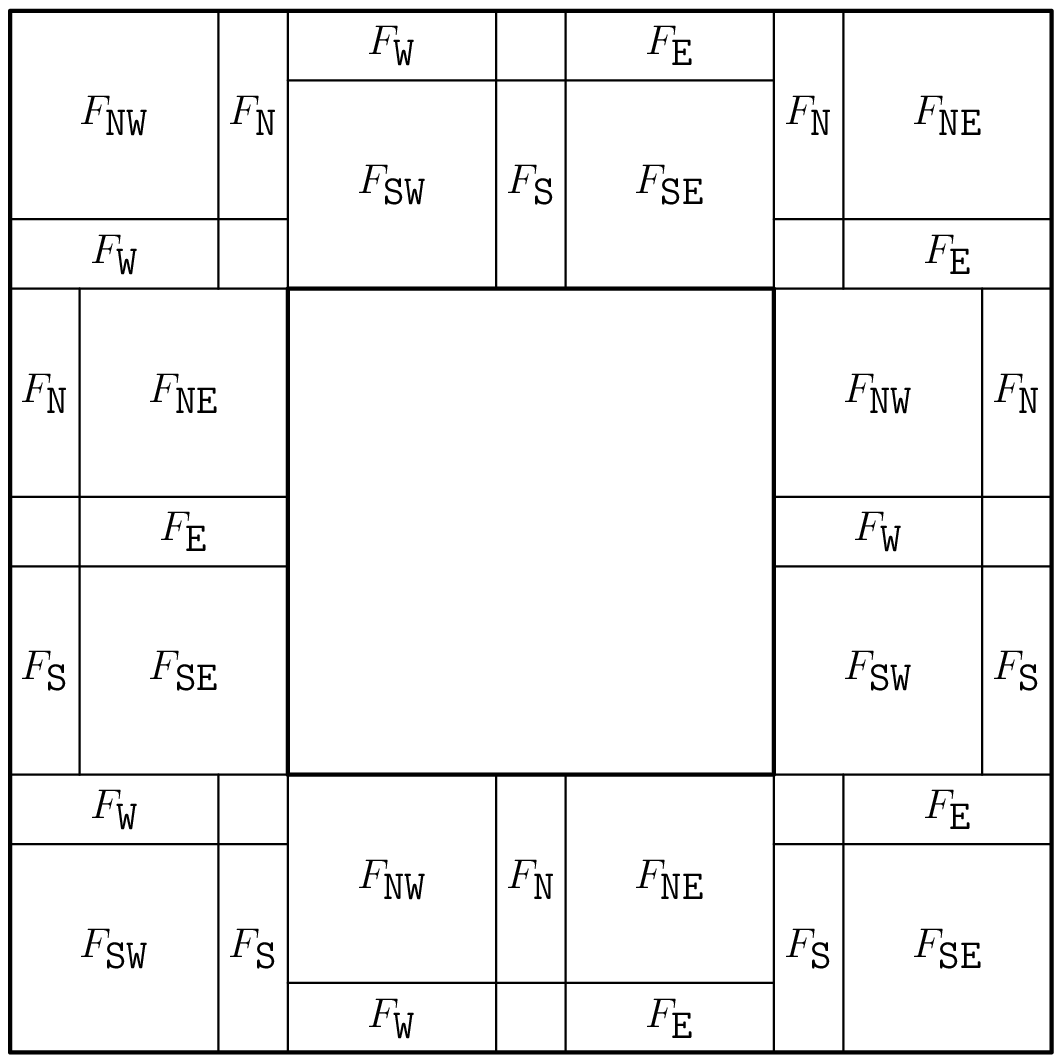}\hspace*{\fill}\caption{}\label{hdef}\end{figure}

\medskip

\begin{lemma} Let $S$ be the vertex set of an $L_F$-cycle within $F$. Then $h(S)$ is the vertex set of an $L_H$-cycle within $H$. \label{h1} \end{lemma}

\medskip

\begin{proof} Analogous to the proof of Lemma \ref{f1}.

Let $S$ be the vertex set of the $L_F$-cycle $C$ and $a$ a square in $C$. The transformation rules for $a$ are as follows.

\medskip

\emph{Case 1.} $a$ belongs to $F_\texttt{E}$. Then the two squares adjacent to $a$ in $C$ are $a + (-n, m)$ and $a + (-n, -m)$. Replace $a$ with the five-square path \[ (a + v^f_\texttt{N})(a + v^f_\texttt{SE})(a + v^f_\texttt{W})(a + v^f_\texttt{NE})(a + v^f_\texttt{S}), \] traversed in such a direction that the endpoint $a + v^f_\texttt{N}$ is on the side of $a + (-n, m)$ and the endpoint $a + v^f_\texttt{S}$ is on the side of $a + (-n, -m)$.

\medskip

\emph{Case 2.} $a$ belongs to $F_\texttt{NE}$. Then the two squares adjacent to $a$ in $C$ are $a + (-n, -m)$ and $a + (-m, -n)$. Replace $a$ with the three-square path \[ (a + v^f_\texttt{S})(a + v^f_\texttt{NE})(a + v^f_\texttt{W}), \] traversed in such a direction that the endpoint $a + v^f_\texttt{S}$ is on the side of $a + (-n, -m)$ and the endpoint $a + v^f_\texttt{W}$ is on the side of $a + (-m, -n)$.

\medskip

All other cases are obtained from Cases 1 and 2 by symmetry. \end{proof}

\medskip

\begin{lemma} Let $S$ be the vertex set of a proper cycle of a leaper $M$ within $F$. Then $h(S)$ is the vertex set of a proper cycle of $M$ within $H$. \label{h2} \end{lemma}

\medskip

The proof is analogous to the proof of Lemma \ref{f2}.

We pause for a moment to point out a deep connection between the three lifting transformations.

Let $m$ and $n$ be arbitrary integers.

Introduce a Cartesian coordinate system $Oxy$ over the infinite chessboard such that the integer points are the centers of the squares if $m + n$ is odd, and the vertices of the squares if it is even, and write $(x, y)$ for the square centered at $(x, y)$.

Define the \emph{standard} $(m, n)$-frame to be the set of all squares $(x, y)$ such that at least one of $|x|$ and $|y|$ exceeds the smaller of $\frac{1}{2}|m - n|$ and $\frac{1}{2}(m + n)$, and both of $|x|$ and $|y|$ are less than the larger. The \texttt{E}-section of a standard $(m, n)$-frame would be the set of all squares $(x, y)$ such that $x$ lies between $\frac{1}{2}(n - m)$ and $\frac{1}{2}(m + n)$ and $y$ lies between $\frac{1}{2}(m - n)$ and $\frac{1}{2}(n - m)$, the \texttt{NE}-section the set of all squares $(x, y)$ such that both $x$ and $y$ lie between $\frac{1}{2}(n - m)$ and $\frac{1}{2}(m + n)$, and so on. In general, a standard frame is not the disjoint union of its sections.

Let, then, $0 \le m < n$. The standard $(n, -m)$-frame consists of the same squares as the standard $(m, n)$-frame. Therefore, we can view an $(m, n)$-frame $F$ as an overlapping inside-out $(n, -m)$-frame $F'$. When we plug the values $(n, -m)$ into the definition of $f$ and lift $F'$ by means of $f$, the result is precisely the same as when we lift $F$ by means of $g$.

Furthermore, the standard $(n, m)$-frame consists of the same squares as the standard $(m, n)$-frame. Therefore, we can view an $(m, n)$-frame $F$ as an overlapping $(n, m)$-frame $F''$. When we plug the values $(n, m)$ into the definition of $f$ and lift $F''$ by means of $f$, the result is precisely the same as when we lift $F$ by means of $h$.

So, all three lifting transformations are in a sense forms of the same fundamental lifting transformation.

We have built all the tools we need and are ready to tackle the theorem.

Let $d$ be the greatest common divisor of $p$ and $q$, $p = dp'$, and $q = dq'$. Then every cycle of $L$ within a $(p + q) \times (p + q)$ chessboard is a scaled, by a factor of $d$, copy of a cycle of the basic $(p', q')$-leaper $L'$ within a $(p' + q') \times (p' + q')$ chessboard. Therefore, it suffices to consider the case when $L$ is a basic leaper.

Let, from here on, $p < q$ and $L$ be a basic leaper.

We will show, by induction on $p + q$, that there exists a second $(r, s)$-leaper $M$, $r + s < p + q$, such that a proper Hamiltonian cycle of $M$ exists over the squares of $C$.

When $p = 1$ and $q = 2$, the theorem holds with $r = 0$ and $s = 1$ as in Figure \ref{1233}.

Let, from here on, $p + q > 3$ and $H$ be a $(p, q)$-frame. We distinguish three cases for the proportions of $H$.

\medskip

\emph{Case $f$.} $3p < q$. Let $m = p$ and $n = q - 2p$. Then $m < n$, $m + n < p + q$, $p = m$, and $q = 2m + n$.

Let $F$ be an $(m, n)$-frame. Then $f$ lifts $F$ to $H$.

By Lemma \ref{f1}, $f$ lifts each $L_F$-cycle within $F$ to an $L_H$-cycle within $H$. Since $f$ maps the set of all squares of $F$ onto the set of all squares of $H$, there exists an $L_F$-cycle $A$ within $F$ such that $f$ lifts $A$ to $C$.

By the induction hypothesis, there exists a proper Hamiltonian cycle of an $(r, s)$-leaper $M$, $r + s < m + n$, over the squares of $A$. By Lemma \ref{f2}, we are done.

\medskip

\emph{Case $g$.} $2p > q$. Let $m = 2p - q$ and $n = p$. Then $m < n$, $m + n < p + q$, $p = n$, and $q = 2n - m$.

Suppose first that $C$ lies within the core of $H$. Let $F$ be an $(m, n)$-frame. Then $g$ lifts $F$ to $H$ and the proof continues as in Case $f$.

Suppose, then, that $C$ lies within the shell of $H$. Then $C$ consists of eight squares, one in each $H^-_i$, and a proper Hamiltonian cycle of an $(m, n)$-leaper $M$ exists over the squares of $C$.

\medskip

\emph{Case $h$.} $2p < q < 3p$. Let $m = q - 2p$ and $n = p$. Then $m < n$, $m + n < p + q$, $p = n$, and $q = m + 2n$.

Suppose first that $C$ lies within the core of $H$. Let $F$ be an $(m, n)$-frame. Then $h$ lifts $F$ to $H$ and the proof continues as in Case $f$.

Suppose, then, that $C$ lies within the shell of $H$. The proof continues as in the second part of Case $g$.

\medskip

This completes the proof of the theorem. \end{proof}

\section{Descents and Even Continued Fractions}

Let us look more closely into the concluding part of the proof of Theorem \ref{sl}. Let $L$ be a basic $(p, q)$-leaper distinct from the $(0, 1)$-leaper and the $(1, 2)$-leaper, and $H$ a $(p, q)$-frame. Then there exist a unique $(m, n)$-frame $F$ and a unique lifting transformation that lifts $F$ to $H$. When we continue this process backwards, eventually it bottoms out at a $(1, 2)$-frame. Therefore, there exists a unique sequence of lifting transformations that lifts a $(1, 2)$-frame to $H$.

\medskip

\begin{definition} The \emph{descent} of a skew basic $(p, q)$-leaper $L$ is the unique string $e_1e_2\ldots e_l$, composed of the characters \texttt{f}, \texttt{g}, and \texttt{h}, such that successively applying lifting transformations of types $e_l$, $e_{l - 1}$, \ldots, $e_1$ to a $(1, 2)$-frame lifts it to a $(p, q)$-frame. \end{definition}

\medskip

This induces a one-to-one mapping between strings composed of the characters \texttt{f}, \texttt{g}, and \texttt{h} and skew basic leapers.

Let us look into several examples. The descent of a $(1, 2)$-leaper is the empty string. The descent of a $(1, 2r)$-leaper is \texttt{ff}\ldots \texttt{f}, where the character \texttt{f} occurs $r - 1$ times. The descent of an $(r, r + 1)$-leaper is \texttt{gg}\ldots \texttt{g}, where the character \texttt{g} occurs $r - 1$ times. And the descent of an $(18, 41)$-leaper is \texttt{hfgh}.

The proof of Theorem \ref{sl} is essentially a proof by induction on descent.

Let $G$ be the leaper graph of a basic $(p, q)$-leaper $L$ over a $(p + q) \times (p + q)$ chessboard. Let us look at the cycles of $G$ from the point of view of the descent of $L$.

Consider a $(p, q)$-frame $F$. Lifting $F$ by means of $f$ extends all existing cycles without altering the second leapers that tour them, and does not create any new cycles. Lifting $F$ by means of either $g$ or $h$ extends all existing cycles without altering the second leapers that tour them, and creates $(q - p)^2$ new eight-square cycles within the shell of the larger frame, all of which are translation copies of each other and each of which is also toured by a $(p, q)$-leaper. Hence the following theorem.

\medskip

\begin{theorem} Let $e = e_1e_2\ldots e_l$ be the descent of the skew basic $(p, q)$-leaper $L$.

Consider all suffixes $e_{i + 1}e_{i + 2}\ldots e_l$ of $e$ such that $e_i$ is either \texttt{g} or \texttt{h}. Let $k$ be one larger than the number of such suffixes, the $i$-th such suffix, $i = 1$, 2, \ldots, $k - 1$, be the descent of a $(p_i, q_i)$-leaper, and $p_k = 0$ and $q_k = 1$.

Then the leaper graph $G$ of $L$ over a $(p + q) \times (p + q)$ chessboard consists of $(q - p)^2$ isolated vertices and a number of disjoint cycles of $k$ distinct types. For all $i$, $G$ contains $(q_i - p_i)^2$ cycles of type $i$, all of which are translation copies of each other and each of which is also toured by a $(p_i, q_i)$-leaper. \label{descent} \end{theorem}

\medskip

It is possible to word Theorem \ref{descent} in terms of continued fractions without referencing descents. In order to do so, let us track what happens to the ratio $r = \frac{q}{p}$ of the proportions of a $(p, q)$-frame $F$ when we lift $F$.

The transformation $f$ lifts $F$ to a $(p, 2p + q)$-frame $H$, and \[\frac{2p + q}{p} = 2 + \frac{q}{p} = 2 + r.\]

The transformation $g$ lifts $F$ to a $(q, 2q - p)$-frame $H$, and \[\frac{2q - p}{q} = 2 - \frac{p}{q} = 2 - \frac{1}{r}.\]

The transformation $h$ lifts $F$ to a $(q, p + 2q)$-frame $H$, and \[\frac{p + 2q}{q} = 2 + \frac{p}{q} = 2 + \frac{1}{r}.\]

Chaining the right-hand sides of the above equations gives us an expression of the form \[\frac{q}{p} = c_1 \pm \cfrac{1}{c_2 \pm \cfrac{1}{\ddots \pm \frac{1}{c_k}}},\] where the terms $c_1$, $c_2$, \ldots, $c_k$ and the $\pm$ signs are determined as follows.

Let $e = e_1e_2\ldots e_l$ be the descent of the skew basic $(p, q)$-leaper $L$ and \[e = e'_1e''_1e'_2e''_2\ldots e''_{k - 1}e'_k\] a partitioning of $e$ into (possibly empty) substrings such that, for all $i$, $e'_i$ is a (possibly empty) run of the character \texttt{f} of length $c'_i$ and $e''_i$ consists of a single character, either \texttt{g} or \texttt{h}. Then $c_i = 2c'_i + 2$ and the sign following $c_i$ is $-$ if $e''_i = \texttt{g}$ and $+$ if $e''_i = \texttt{h}$.

\medskip

\begin{definition} An \emph{even continued fraction} is an expression of the form \[[c_1\pm, c_2\pm, \ldots, c_k] = c_1 \pm \cfrac{1}{c_2 \pm \cfrac{1}{\ddots \pm \frac{1}{c_k}}},\] where $c_i$ is an even integer for all $i$, $c_1$ is nonnegative, and $c_2$, $c_3$, \ldots, $c_k$ are all positive. \end{definition}

\medskip

A nonnegative rational number $\frac{q}{p}$ with $p$ and $q$ relatively prime possesses a finite even continued fraction representation if and only if $p + q$ is odd. Moreover, that representation is unique.

We are ready to give the alternative form of Theorem \ref{descent}.

\medskip

\begin{theorem} Let $L$ be a skew basic $(p, q)$-leaper. Let \[\frac{q}{p} = [c_1\pm, c_2\pm, \ldots, c_k]\] be the representation of $\frac{q}{p}$ as an even continued fraction, \[\frac{q_i}{p_i} = [c_{i + 1}\pm, c_{i + 2}\pm, \ldots, c_k]\] for $i = 1$, 2, \ldots, $k - 1$, where $\frac{q_i}{p_i}$ is irreducible, and $p_k = 0$ and $q_k = 1$.

Then the leaper graph $G$ of $L$ over a $(p + q) \times (p + q)$ chessboard consists of $(q - p)^2$ isolated vertices and a number of disjoint cycles of $k$ distinct types. For all $i$, $G$ contains $(q_i - p_i)^2$ cycles of type $i$, all of which are translation copies of each other and each of which is also toured by a $(p_i, q_i)$-leaper. \label{ecfrac} \end{theorem}

\medskip

Here follow a couple of noteworthy corollaries of Theorems \ref{descent} and \ref{ecfrac}.

\medskip

\begin{corollary} The leaper graph of a $(p, q)$-leaper $L$, $p < q$, over a $(p + q) \times (p + q)$ chessboard contains a single cycle (or, equivalently, the leaper graph of $L$ over a $(p, q)$-frame is nonempty and connected) if and only if $p = 1$ and $q$ is even. \end{corollary}

\medskip

\begin{corollary} Given a skew basic $(p, q)$-leaper $L$, there exists a unique cycle $C$ of $L$ within a $(p + q) \times (p + q)$ chessboard such that there exists a Hamiltonian cycle of a $(0, 1)$-leaper over the squares of $C$. \end{corollary}

\medskip

We go on to give a concise expression for the length of each cycle of type $i$ in terms of even continued fractions. First, however, we need to lay some groundwork.

\medskip

\begin{theorem} Let $L$ be a $(p, q)$-leaper with $p + q$ odd, $C$ a cycle of $L$ within a $(p + q) \times (p + q)$ chessboard $B$, and $F$ the $(p, q)$-frame obtained from $B$ by removing the central $(q - p)^2$ squares.

Then $C$ visits an odd number of squares (and, in particular, at least one square) in each section of $F$. Furthermore, $C$ visits the same number of squares in each of the four side sections $F_\texttt{E}$, $F_\texttt{N}$, $F_\texttt{W}$, and $F_\texttt{S}$ of $F$, and the same number of squares in each of the four corner sections $F_\texttt{NE}$, $F_\texttt{NW}$, $F_\texttt{SW}$, and $F_\texttt{SE}$ of $F$. \label{seclen} \end{theorem}

\medskip

\begin{proof} By induction on descent. \end{proof}

\medskip

\begin{theorem} Let, in the setting of Theorem \ref{ecfrac}, for all $i$ \[\frac{l_i}{d_i} = [c_1\pm, 2-, c_2\pm, 2-, \ldots, c_i],\] where $\frac{l_i}{d_i}$ is irreducible. Then the length of each cycle of type $i$ is $4l_i$. \label{len} \end{theorem}

\medskip

\begin{proof} Let $F$ be the $(p, q)$-frame obtained from the $(p + q) \times (p + q)$ chessboard by removing all squares isolated in $G$, and $C$ a cycle in $G$.

Let $d$ be the number of squares that $C$ visits within $F_\texttt{E}$ and $l$ the number of squares that $C$ visits within $F_\texttt{E} \cup F_\texttt{NE}$. By Theorem \ref{seclen}, the length of $C$ is $4l$.

When $C$ is contained within the shell of $F$ (provided that $F$ does indeed possess a shell), it is of type 1, $c_1 = 2$, $d = d_1 = 1$, and $l = l_1 = 2$.

Let us track what happens to the ratio $r = \frac{l}{d}$ when we lift $F$ to a larger frame $H$. Let $E$ be the image of $C$ within $H$, and define $d'$ and $l'$ analogously to $d$ and $l$, but based on $E$ and $H$.

When we lift $F$ by means of $f$, $d' = d$, $l' = 2d + l$, and \[\frac{l'}{d'} = \frac{2d + l}{d} = 2 + \frac{l}{d} = 2 + r.\]

When we lift $F$ by means of $g$, $d' = 2l - d$, $l' = 3l - 2d$, and \[\frac{l'}{d'} = \frac{3l - 2d}{2l - d} = 2 - \cfrac{1}{2 - \frac{d}{l}} = 2 - \cfrac{1}{2 - \frac{1}{r}}.\]

When we lift $F$ by means of $h$, $d' = 2l - d$, $l' = 5l - 2d$, and \[\frac{l'}{d'} = \frac{5l - 2d}{2l - d} = 2 + \cfrac{1}{2 - \frac{d}{l}} = 2 + \cfrac{1}{2 - \frac{1}{r}}.\]

Chaining the right-hand sides of the above equations yields the theorem. \end{proof}

\medskip

Occasionally, there exists a third leaper besides the second one.

\medskip

\begin{theorem} Let, in the setting of Theorems \ref{ecfrac} and \ref{len}, $c_k = 2$ (or, equivalently, the final character $e_l$ in the descent $e_1e_2\ldots e_l$ of $L$ be either \texttt{g} or \texttt{h}) and $l_k$ be indivisible by three. Then there exists a Hamiltonian cycle of a $(1, 2)$-leaper over the squares of the unique cycle of type $k$. \label{tl} \end{theorem}

\medskip

Let us look into a couple of examples.

\begin{figure}[ht!]\hspace*{\fill}\includegraphics[scale=0.66]{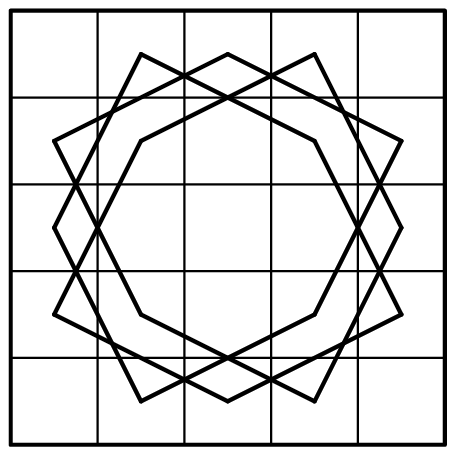}\hspace*{\fill}\caption{}\label{2355e}\end{figure}

When $p = 2$ and $q = 3$, $k = 2$, $c_2 = 2$, and $l_2 = 16$. Figure \ref{2355e} shows that the unique cycle of type 2, previously depicted in Figure \ref{2355}, top, is also toured by a $(1, 2)$-leaper.

\begin{figure}[ht!]\hspace*{\fill}\includegraphics[scale=0.66]{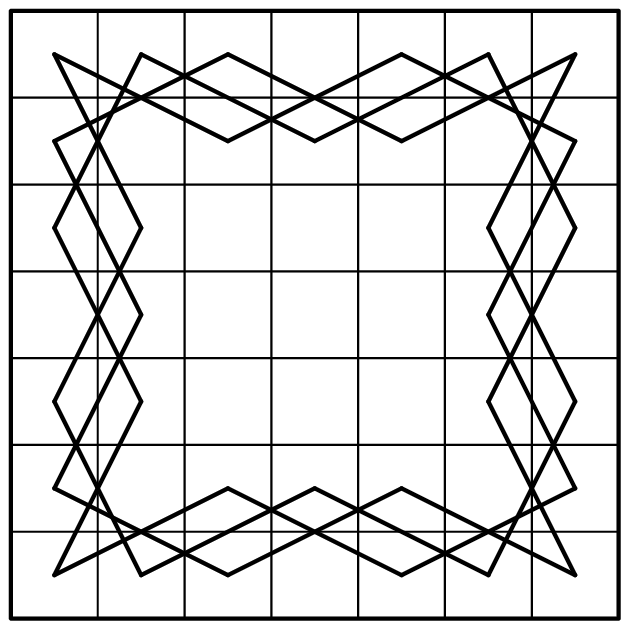}\hspace*{\fill}\caption{}\label{2577e}\end{figure}

When $p = 2$ and $q = 5$, $k = 2$, $c_2 = 2$, and $l_2 = 32$. Figure \ref{2577e} shows that the unique cycle of type 2, previously depicted in Figure \ref{2577}, top, is also toured by a $(1, 2)$-leaper.

\medskip

\begin{proof} When $k = 1$, the theorem holds as in Figure \ref{1233}. Let, from here on, $k \ge 2$.

Let $a_1a_2\ldots a_{l_k}$ be a proper Hamiltonian cycle of a $(0, 1)$-leaper over the squares of the unique cycle $C$ in $G$ of type $k$, as constructed in the proof of Theorem \ref{sl}, $a_{i + l_k} \equiv a_i$ for all integers $i$, and $F$ the $(p, q)$-frame obtained from the $(p + q) \times (p + q)$ chessboard by removing all squares isolated in $G$.

\medskip

\begin{lemma} The cycle $C$ visits at least three squares in each section of $F$ and the squares $a_i$ and $a_{i + 3}$ are linked by a $(1, 2)$-move for all integers $i$. \label{knight} \end{lemma}

\medskip

\begin{proof} When $p = 2$ and either $q = 3$ or $q = 5$, the claim holds as in Figures \ref{2355e} and \ref{2577e}. The proof continues by induction on descent.

Suppose that the lemma holds for $C$. Lift $F$ to a larger frame $H$ by means of any of the three lifting transformations. Let $E$ be the image of $C$ within $H$ and $b_1b_2\ldots b_s$ a proper Hamiltonian cycle of a $(0, 1)$-leaper over the squares of $E$, as constructed in the proof of Theorem \ref{sl}, with $b_{j + s} \equiv b_j$ for all integers $j$.

By the definitions of $f$, $g$, and $h$ and the induction hypothesis, $E$ visits at least three squares in each section of $H$.

It follows that, for all $j$, $b_j$ and $b_{j + 3}$ are either within the same subboard of $H$ in Figure \ref{fdef}, \ref{gdef}, or \ref{hdef}, or within the union of two subboards of $H$ visited in direct succession by $b_1b_2\ldots b_s$.

Therefore, there always exist two squares $a_i$ and $a_{i + 3}$ within $F$ such that an $f$, $g$, or $h$-translation maps $a_i$ and $a_{i + 3}$ onto $b_j$ and $b_{j + 3}$.

By the induction hypothesis, $a_i$ and $a_{i + 3}$ are linked by a $(1, 2)$-move. Therefore, so are $b_j$ and $b_{j + 3}$. \end{proof}

\medskip

Consider the sequence of squares \[a_3, a_6, \ldots, a_{l_k},\] the $i$-th term of which equals $a_{3i}$ for $i = 1$, $2$, \ldots, $l_k$.

By Lemma \ref{knight}, this sequence is a cycle of a $(1, 2)$-leaper. Furthermore, as $l_k$ is indivisible by three, it contains each square of $C$ precisely once. Therefore, it is a Hamiltonian cycle of a $(1, 2)$-leaper over the squares of $C$. \end{proof}

\medskip

The proof of Theorem \ref{tl} does not appear to generalize beyond the $(1, 2)$-leaper. For more on this topic, see Questions \ref{multleap} and \ref{multgraph}.

\section{Connectedness}

We take a brief detour to expand on an application of the technique of induction on descent.

In \cite{K}, Donald Knuth shows that a $(p, q)$-leaper $L$, $p \le q$, is free over a regular chessboard $B$ (of size greater than $1 \times 1$) if and only if $L$ is basic and $B$ contains a $(p + q) \times 2q$ subboard. We proceed to give a proof of the more involved ``if'' part by induction on descent.

\medskip

\begin{theorem} (Donald Knuth, \cite{K}) Let $L$ be a basic $(p, q)$-leaper, $p < q$. Then $L$ is free over a $(p + q) \times 2q$ chessboard. \label{connectedness} \end{theorem}

\medskip

\begin{proof} Let $B$ be a $(p + q) \times 2q$ chessboard. We view $B$ as the union of $q - p + 1$ translation copies $F + (t, 0)$ of a $(p, q)$-frame $F$, $t = 0$, 1, \ldots, $q - p$.

Let $C_1$, $C_2$, \ldots, $C_c$ be all $L_F$-cycles within $F$. Then the leaper graph of $L$ over $B$ is the union of the translation copies $C_k + (t, 0)$ for $k = 1$, 2, \ldots, $c$ and $t = 0$, 1, \ldots, $q - p$.

Construct the graph $K$ as follows. The vertices of $K$ are the ordered pairs $(k, t)$ where $k = 1$, 2, \ldots, $c$ and $t = 0$, 1, \ldots, $q - p$. Two vertices $(k', t')$ and $(k'', t'')$ in $K$ are joined by an edge if and only if there exists a square in $B$ that belongs to both $C_{k'} + (t', 0)$ and $C_{k''} + (t'', 0)$. Then the leaper graph of $L$ over $B$ is connected if and only if $K$ is.

When $p = 1$ and $q = 2$, $K$ consists of two vertices joined by an edge.

Suppose, then, that $K$ is connected. Lift $F$ to an $(m, n)$-frame $H$ by means of any of the three lifting transformations.

For simplicity, we extend the definitions of shell and core as follows. If $3m \le n$, the shell of $H$ is empty and the core of $H$ is $H$.

Let $d$ be the number of shell cycles in $H$ (zero if the lifting transformation is $f$ and $H$ only possesses a shell in the extended sense) and $D_1$, $D_2$, \ldots, $D_{c + d}$ all $L_H$-cycles within $H$, so that lifting $C_i$ gives $D_{i + d}$ for all $i$.

Define $E$ and $N$ analogously to $B$ and $K$, but based on $H$. It suffices to show that $N$ is connected.

Let $E^+$ be the union of the translation copies $H^+ + (t, 0)$, $t = 0$, 1, \ldots, $n - m$, and $N^+$ the spanning subgraph of $N$ over all vertices $(k, t)$ such that $k > d$, corresponding to cycles in the core of $H$.

Let $(k', t')$ and $(k'', t'')$, $t' < t''$, be joined by an edge in $K$. Then there exist two squares $a'$ and $a''$ in $F$ that are joined by a $(t'' - t', 0)$ move such that $a'$ belongs to $C_{k'}$ and $a''$ belongs to $C_{k''}$.

Since $t'' - t' \le q - p$ does not exceed the side of the central hole of $F$, both of $a'$ and $a''$ belong to one of the following subboards of $F$: \[ F_\texttt{E}, F_\texttt{NE} \cup F_\texttt{N}, F_\texttt{N} \cup F_\texttt{NW}, F_\texttt{W}, F_\texttt{SW} \cup F_\texttt{S}, F_\texttt{S} \cup F_\texttt{SE}. \]

By the definitions of $f$, $g$, and $h$, there exist in $H^+$ two images $b'$ and $b''$ of $a'$ and $a''$ such that $b'$ and $b''$ are joined by a $(t'' - t', 0)$ move, $b'$ belongs to $D_{k' + d}$ and $b''$ belongs to $D_{k'' + d}$.

It follows that, whenever an edge joins $(k', t')$ and $(k'', t'')$ in $K$ and $t'' - t' = u'' - u'$, $t' < t''$, $0 \le u' < u'' \le n - m$, an edge joins $(k' + d, u')$ and $(k'' + d, u'')$ in $N^+$. Therefore, $N^+$ is connected.

We are left to take care of the shell cycles of $H$. However, since all of the translation copies $H^-_\texttt{N} + (t, 0)$, $t = 0$, 1, \ldots, $n - m$, are subboards of $E^+$, every vertex $(k, t)$ of $N$ such that $k \le d$, corresponding to a cycle in the shell of $H$, is joined by an edge to a vertex in $N^+$. Therefore, $N$ is connected. \end{proof}

\section{Direction Graphs}

\begin{figure}[ht!]\hspace*{\fill}\includegraphics[scale=0.66]{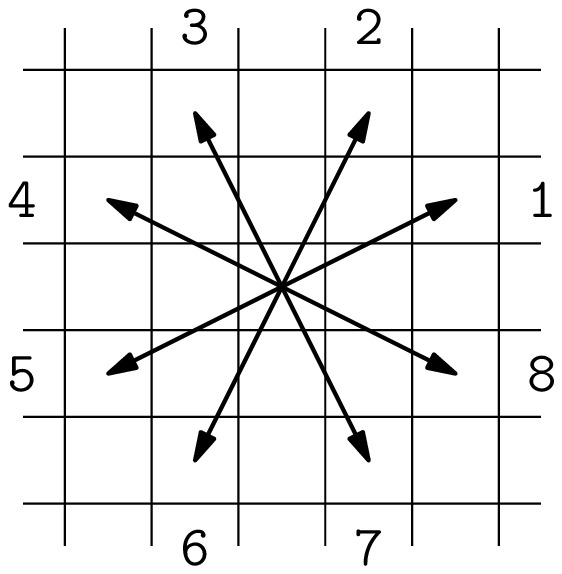}\hspace*{\fill}\caption{}\label{dir}\end{figure}

Let $L$ be a skew leaper. Then there exist eight possible \emph{skew directions} for the moves of $L$, east-northeast, north-northeast, \ldots, east-southeast, which we label \texttt{1} through \texttt{8} starting from east-northeast and proceeding counterclockwise, as in Figure \ref{dir}.

Given a skew direction $i$ and an integer $k$, we write $-i$ for the direction opposite $i$ and $i + k$ for the direction $k$ steps counterclockwise from $i$. For instance, $-\texttt{1} = \texttt{5}$, $\texttt{2} + 1 = \texttt{3}$, for all skew directions $i$, $-i = i + 4$, and, for all skew directions $i$ and integers $k$, $-(i + k) = -i + k$.

We associate a $2 \times 2$ matrix with each skew direction as follows.

\medskip

\begin{definition} The \emph{direction matrix} of a move of a skew $(p, q)$-leaper $L$, $p < q$, of direction $i$ from a square $a$ to a square $b = a + v$ is the unique matrix $A_i$ out of the following eight, \begin{align*} A_\texttt{1} &= \left(\begin{array}{cc} 0 & 1\\ 1 & 0 \end{array}\right), & A_\texttt{2} &= \left(\begin{array}{cc} 1 & 0\\ 0 & 1 \end{array}\right),\\ A_\texttt{3} &= \left(\begin{array}{cc} -1 & 0\\ 0 & 1 \end{array}\right), & A_\texttt{4} &= \left(\begin{array}{cc} 0 & -1\\ 1 & 0 \end{array}\right),\\ A_\texttt{5} &= \left(\begin{array}{cc} 0 & -1\\ -1 & 0 \end{array}\right), & A_\texttt{6} &= \left(\begin{array}{cc} -1 & 0\\ 0 & -1 \end{array}\right),\\ A_\texttt{7} &= \left(\begin{array}{cc} 1 & 0\\ 0 & -1 \end{array}\right), & A_\texttt{8} &= \left(\begin{array}{cc} 0 & 1\\ -1 & 0 \end{array}\right), \end{align*} such that \[v^T = A_i\left(\begin{array}{c} p\\ q \end{array}\right).\] \end{definition}

\medskip

For all skew directions $i$, $A_{-i} = -A_i$.

\medskip

\begin{definition} A \emph{direction graph} is a symmetric directed graph $\Phi$ whose arcs are labeled with skew directions in such a way that the sum of the associated direction matrices over every simple cycle in $\Phi$ is the zero matrix. \end{definition}

\medskip

In particular, in a direction graph $\Phi$, for every arc pointing from $x$ to $y$ and labeled $i$, the arc pointing from $y$ to $x$ is labeled $-i$. Hence the following definition.

\medskip

\begin{definition} A labeled oriented graph $\Psi$ \emph{represents} a direction graph $\Phi$ if $\Psi$ is obtained from $\Phi$ by removing exactly one arc out of each symmetric pair. \end{definition}

\medskip

A direction graph is completely determined by any labeled oriented graph that represents it.

\medskip

\begin{definition} A direction graph $\Phi$ is \emph{extracted from} a leaper graph $G$ of a skew leaper $L$ if there exists a one-to-one mapping $\sigma$ between the vertices of $G$ and the vertices of $\Phi$ such that an $L$-move of direction $i$ leads from a vertex $a$ to a vertex $b$ in $G$ if and only if an arc labeled $i$ points from $\sigma(a)$ to $\sigma(b)$ in $\Phi$. \end{definition}

\medskip

In other words, a direction graph is extracted from a leaper graph by abstracting away all information (such as the precise positions of the squares and the proportions of the leaper) save for the directions of the moves.

\begin{figure}[ht!]\hspace*{\fill}\includegraphics[scale=0.66]{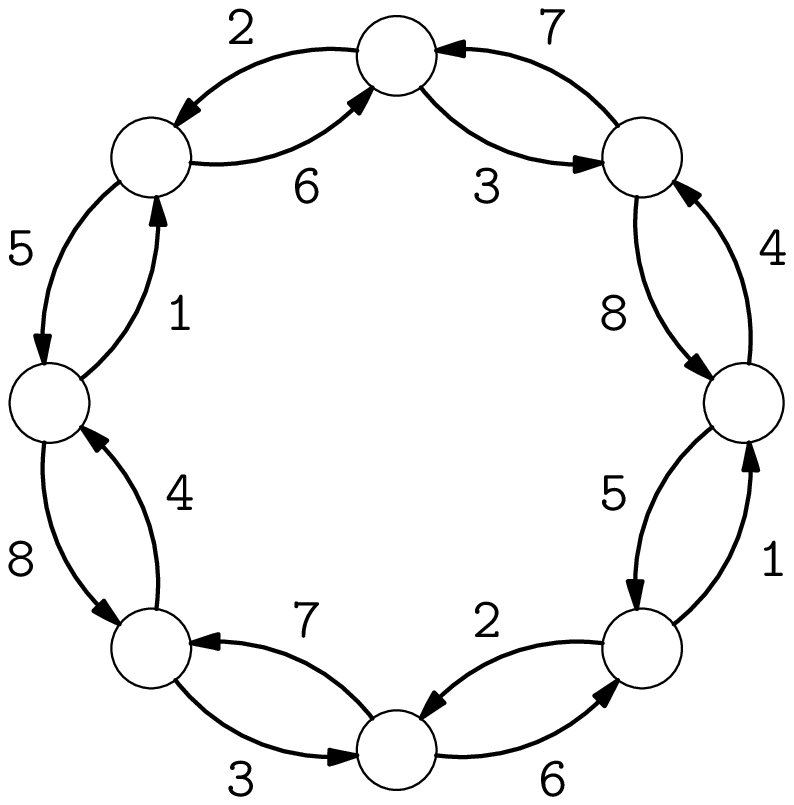}\hspace*{\fill}\caption{}\label{1233dir}\end{figure}

For instance, the direction graph extracted from the $(1, 2)$-cycle in Figure \ref{1233} is depicted in Figure \ref{1233dir}, and is represented by the oriented cycle labeled \texttt{47258361}.

It may happen that a direction graph cannot be extracted from a given leaper graph.

\medskip

\begin{definition} A simple oriented cycle $C$ of a skew leaper is \emph{trivial} if the sum of the associated direction matrices over $C$ is the zero matrix, and \emph{nontrivial} otherwise. A simple cycle $C$ of a skew leaper is trivial if any (or, equivalently, both) of its two orientations are, and nontrivial otherwise. \end{definition}

\medskip

Given a leaper graph $G$, it is possible to extract a direction graph $\Phi$ from $G$ if and only if $G$ does not contain a nontrivial cycle.

For instance, the $(1, 2)$-cycle and the $(1, 3)$-cycle in Figure \ref{1312} are both nontrivial, and it is not possible to extract a direction graph from either.

\medskip

\begin{definition} Let $\Phi$ be a direction graph and $L$ a leaper. A graph $G$ is an \emph{$L$-instantiation} of $\Phi$ if the vertices of $G$ are squares and there exists a mapping $\tau$ from the vertices of $\Phi$ to the vertices of $G$ such that, for every arc in $\Phi$ pointing from $x$ to $y$ and labeled $i$, an $L$-move of direction $i$ leads from $\tau(x)$ to $\tau(y)$. \end{definition}

\medskip

Let $\Phi$ be the direction graph extracted from a leaper graph $G$ of a leaper $L$. Then $G$ is an $L$-instantiation of $\Phi$.

Given a direction graph $\Phi$ and a leaper $L$, there always exists an $L$-in\-stan\-ti\-a\-tion of $\Phi$, unique up to translation provided that $\Phi$ is connected. It is well-defined for arbitrary leapers $L$, as when $L$ is an orthogonal or diagonal $(p, q)$-leaper, $p \le q$, we can define an $L$-move of direction $i$ as an $L$-move of translation $v$ such that $v^T = A_i \left( \begin{smallmatrix} p\\ q \end{smallmatrix} \right)$. However, it may happen that the mapping $\tau$ is not one-to-one or that the leaper graph of $L$ over the vertex set of $G$ is distinct from $G$.

We go on to delineate the class of direction graphs for which instantiation is well-behaved.

Let $x$ and $y$ be vertices in a direction graph $\Phi$. Then the sum of the associated direction matrices over every path from $x$ to $y$ in $\Phi$ is the same.

\medskip

\begin{definition} Let $\Phi$ be a direction graph and $x$ and $y$ vertices in the same connected component of $\Phi$. The \emph{distance} from $x$ to $y$ in $\Phi$ is the sum of the associated direction matrices over any path from $x$ to $y$ in $\Phi$. \end{definition}

\medskip

\begin{definition} A direction graph $\Phi$ is \emph{coherent} if, for every pair of vertices $x$ and $y$ in the same connected component of $\Phi$, the distance from $x$ to $y$ in $\Phi$ equals the zero matrix if and only if $x$ and $y$ coincide, and a direction matrix if and only if $x$ and $y$ are joined by an arc. \end{definition}

\medskip

Every direction graph extracted from a leaper graph is coherent.

Given a coherent direction graph $\Phi$, let $L$ be a skew $(p, q)$-leaper, $p < q$, and $G$ an $L$-in\-stan\-ti\-a\-tion of $\Phi$ such that the images under $\tau$ of vertices in different connected components of $\Phi$ do not coincide and are not joined by an $L$-move. Furthermore, let $x$ and $y$ be vertices in the same connected component of $\Phi$, $A$ the distance from $x$ to $y$ in $\Phi$, and $v$ the translation defined by $v^T = A \left( \begin{smallmatrix} p\\ q \end{smallmatrix} \right)$. Then $\tau(x) + v = \tau(y)$.

When $A$ is distinct from the zero matrix, there exists at most one skew basic $(p, q)$-leaper $L$ such that $v$ is the zero translation, and when $A$ is distinct from the direction matrix $A_i$, there exists at most one skew basic $(p, q)$-leaper $L$ such that $v$ is an $L$-translation of direction $i$.

It follows that the mapping $\tau$ is one-to-one and the leaper graph of $L$ over the vertex set of $G$ is $G$ for all but finitely many skew basic leapers $L$. In other words, $G$ is a leaper graph of $L$ and $\Phi$ is extracted from $G$ for all but finitely many skew basic leapers $L$.

Every leaper graph of an orthogonal or diagonal leaper is an instantiation of a direction graph. A leaper graph of a skew leaper is an instantiation of a direction graph if and only if it does not contain a nontrivial cycle.

\medskip

\begin{definition} A permutation $\pi$ of the eight skew directions is an \emph{equivalence permutation} if there exist two real-coefficient $2 \times 2$ matrices $P$ and $Q$, \emph{inducing} $\pi$, such that \[A_{\pi(i)} = PA_iQ\] for all skew directions $i$. \end{definition}

\medskip

Given an equivalence permutation $\pi$ induced by $P$ and $Q$ and a directed graph $\Phi$ whose arcs are labeled with skew directions, we write $\pi(\Phi)$ or $P\Phi Q$ for the labeled directed graph obtained from $\Phi$ by applying $\pi$ to every arc label. If $\Phi$ is a direction graph, then so is $\pi(\Phi)$. Furthermore, if $\Phi$ is a coherent direction graph, then so is $\pi(\Phi)$.

\medskip

\begin{definition} Two direction graphs $\Phi_1$ and $\Phi_2$ are \emph{equivalent} if there exists an equivalence permutation $\pi$ such that $\pi(\Phi_1) = \Phi_2$. \end{definition}

\medskip

We proceed to look at Theorems \ref{sl}, \ref{descent}, and \ref{ecfrac} from the point of view of direction graphs.

Let $L$ be a skew basic $(p, q)$-leaper, $p < q$, $C$ a cycle of $L$ within a $(p + q) \times (p + q)$ chessboard, and $F$ the associated $(p, q)$-frame.

Impose an orientation on $C$ and consider a square $a$ of $C$ in the section $F_\texttt{E}$ of $F$. The two squares adjacent to $a$ in $C$ are $a + (-q, p)$ and $a + (-q, -p)$. Therefore, the directions of the moves to and from $a$ in $C$ are either \texttt{1} and \texttt{4}, or \texttt{8} and \texttt{5}.

Analogous reasoning applies to all sections of $F$. It follows that the directions of the moves to and from every square $a$ in $C$ are of the form $i$ and $i \pm 3$.

\medskip

\begin{definition} Let $L$ be a skew basic $(p, q)$-leaper and $C$ a cycle of $L$ within a $(p + q) \times (p + q)$ chessboard.

Impose an orientation on $C$ and label every square $a$ in $C$ $+_\texttt{s}$, $+_\texttt{c}$, $-_\texttt{s}$, or $-_\texttt{c}$ as follows. Let $i$ be the direction of the move to $a$ in $C$. Then the $+$ or $-$ signifies whether the direction of the move from $a$ in $C$ is $i + 3$ or $i - 3$, and the \texttt{s} or \texttt{c} signifies whether $a$ belongs to a side or a corner section of the associated frame $F$.

Equivalently, label $a$ according to the following table. \[ \begin{tabular}{c c} Label of $a$ & Directions of moves to and from $a$\\ \hline $+_\texttt{s}$ & \texttt{14}, \texttt{36}, \texttt{58}, \texttt{72}\\  $+_\texttt{c}$ & \texttt{25}, \texttt{47}, \texttt{61}, \texttt{83} \\ $-_\texttt{s}$ & \texttt{27}, \texttt{41}, \texttt{63}, \texttt{85} \\ $-_\texttt{c}$ & \texttt{16}, \texttt{38}, \texttt{52}, \texttt{74} \end{tabular} \] \end{definition}

Then the cyclic string formed by the vertex labels of $C$ is a \emph{signature} of $C$.

\medskip

For instance, a signature of the $(1, 2)$-cycle in Figure \ref{1233} is $+_\texttt{s}+_\texttt{c}+_\texttt{s}+_\texttt{c}+_\texttt{s}+_\texttt{c}+_\texttt{s}+_\texttt{c}$. This is also a signature of every shell cycle, such as the $(2, 3)$-cycle in Figure \ref{2355}, bottom, and the $(2, 5)$-cycle in Figure \ref{2577}, bottom.

Let us track how a signature evolves when we lift a cycle.

Lift $C$ by means of $f$ to $f(C)$. By the proof of Lemma \ref{f1}, a signature of $f(C)$ is obtained from a signature of $C$ by replacing every character with a string composed of the characters $+_\texttt{s}$, $+_\texttt{c}$, $-_\texttt{s}$, and $-_\texttt{c}$, subject to the following system of rewriting rules. \begin{align*} +_\texttt{s} &\to -_\texttt{s} & +_\texttt{c} &\to +_\texttt{s}+_\texttt{c}+_\texttt{s}\\ -_\texttt{s} &\to +_\texttt{s} & -_\texttt{c} &\to -_\texttt{s}-_\texttt{c}-_\texttt{s} \end{align*}

\begin{figure}[ht!]\hspace*{\fill}\includegraphics[scale=0.66]{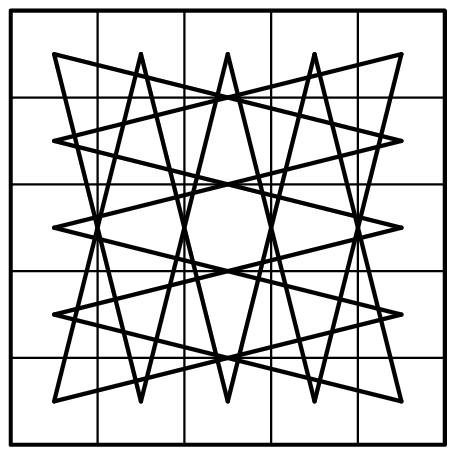}\hfill\includegraphics[scale=0.66]{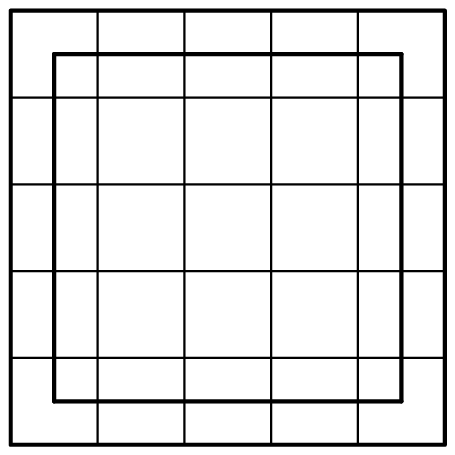}\hspace*{\fill}\caption{}\label{1455}\end{figure}

We refer to this transformation as an \emph{$f$-rewrite}. For instance, the $f$-rewrite of $+_\texttt{s}+_\texttt{c}+_\texttt{s}+_\texttt{c}+_\texttt{s}+_\texttt{c}+_\texttt{s}+_\texttt{c}$ is $-_\texttt{s}+_\texttt{s}+_\texttt{c}+_\texttt{s}-_\texttt{s}+_\texttt{s}+_\texttt{c}+_\texttt{s}-_\texttt{s}+_\texttt{s}+_\texttt{c}+_\texttt{s}-_\texttt{s}+_\texttt{s}+_\texttt{c}+_\texttt{s}$, and this is a signature of the $(1, 4)$-cycle in Figure \ref{1455}.

Next lift $C$ by means of $g$ to $g(C)$. By the proof of Lemma \ref{g1}, a signature of $g(C)$ is obtained from a signature of $C$ by replacing every character with a string composed of the characters $+_\texttt{s}$, $+_\texttt{c}$, $-_\texttt{s}$, and $-_\texttt{c}$, subject to the following system of rewriting rules. \begin{align*} +_\texttt{s} &\to +_\texttt{c}+_\texttt{s}+_\texttt{c} & +_\texttt{c} &\to -_\texttt{c}\\ -_\texttt{s} &\to -_\texttt{c}-_\texttt{s}-_\texttt{c} & -_\texttt{c} &\to +_\texttt{c} \end{align*}

We refer to this transformation as a \emph{$g$-rewrite}. For instance, the $g$-rewrite of $+_\texttt{s}+_\texttt{c}+_\texttt{s}+_\texttt{c}+_\texttt{s}+_\texttt{c}+_\texttt{s}+_\texttt{c}$ is $+_\texttt{c}+_\texttt{s}+_\texttt{c}-_\texttt{c}+_\texttt{c}+_\texttt{s}+_\texttt{c}-_\texttt{c}+_\texttt{c}+_\texttt{s}+_\texttt{c}-_\texttt{c}+_\texttt{c}+_\texttt{s}+_\texttt{c}-_\texttt{c}$, and this is a signature of the $(2, 3)$-cycle in Figure \ref{2355}, top.

Lastly, lift $C$ by means of $h$ to $h(C)$. By the proof of Lemma \ref{h1}, a signature of $h(C)$ is obtained from a signature of $C$ by replacing every character with a string composed of the characters $+_\texttt{s}$, $+_\texttt{c}$, $-_\texttt{s}$, and $-_\texttt{c}$, subject to the following system of rewriting rules. \begin{align*} +_\texttt{s} &\to +_\texttt{s}+_\texttt{c}+_\texttt{s}+_\texttt{c}+_\texttt{s} & +_\texttt{c} &\to -_\texttt{s}-_\texttt{c}-_\texttt{s}\\ -_\texttt{s} &\to -_\texttt{s}-_\texttt{c}-_\texttt{s}-_\texttt{c}-_\texttt{s} & -_\texttt{c} &\to +_\texttt{s}+_\texttt{c}+_\texttt{s} \end{align*}

We refer to this transformation as an \emph{$h$-rewrite}. For instance, the $h$-rewrite of $+_\texttt{s}+_\texttt{c}+_\texttt{s}+_\texttt{c}+_\texttt{s}+_\texttt{c}+_\texttt{s}+_\texttt{c}$ is $+_\texttt{s}+_\texttt{c}+_\texttt{s}+_\texttt{c}+_\texttt{s}-_\texttt{s}-_\texttt{c}-_\texttt{s}+_\texttt{s}+_\texttt{c}+_\texttt{s}+_\texttt{c}+_\texttt{s}-_\texttt{s}-_\texttt{c}-_\texttt{s}+_\texttt{s}+_\texttt{c}+_\texttt{s}+_\texttt{c}+_\texttt{s}-_\texttt{s}-_\texttt{c}-_\texttt{s}+_\texttt{s}+_\texttt{c}+_\texttt{s}+_\texttt{c}+_\texttt{s}-_\texttt{s}-_\texttt{c}-_\texttt{s}$, and this is a signature of the $(2, 5)$-cycle in Figure \ref{2577}, top.

Given a string $e = e_1e_2\ldots e_l$ composed of the characters \texttt{f}, \texttt{g}, and \texttt{h} and a string $s$ composed of the characters $+_\texttt{s}$, $+_\texttt{c}$, $-_\texttt{s}$, and $-_\texttt{c}$, we write $R_e(s)$ for the string obtained by successively applying rewrites of types $e_l$, $e_{l - 1}$, \ldots, $e_1$ to $s$.

\medskip

\begin{definition} Let $L$ be a skew basic $(p, q)$-leaper, $C$ a cycle of $L$ within a $(p + q) \times (p + q)$ chessboard, and $e = e_1e_2\ldots e_l$ a prefix  of the descent of $L$ such that $C$ is the product of successively applying lifting transformations of types $e_l$, $e_{l - 1}$, \ldots, $e_1$ to a shell cycle. Then $e$ is the \emph{descent} of $C$. \end{definition}

\medskip

Let $e$ be the descent of $C$. Then a signature of $C$ is \[R_e(+_\texttt{s}+_\texttt{c}+_\texttt{s}+_\texttt{c}+_\texttt{s}+_\texttt{c}+_\texttt{s}+_\texttt{c}).\]

We proceed to show that $C$ is trivial and that a signature of $C$ completely determines the direction graph extracted from $C$.

\medskip

\begin{theorem} Let $L$ be a $(p, q)$-leaper with $p + q$ odd and $C$ a cycle of $L$ within a $(p + q) \times (p + q)$ chessboard. Then $C$ possesses orientation-preserving fourfold rotational symmetry together with orientation-reversing axial symmetry along a vertical, a horizontal, and two diagonal axes. \label{symm} \end{theorem}

\medskip

\begin{proof} By induction on descent. \end{proof}

\medskip

By Theorem \ref{symm}, the moves of $C$ occur in pairs symmetric with respect to the center of symmetry of $C$. The sum of the associated direction matrices over each such pair is the zero matrix. Therefore, $C$ is trivial.

Let $C = a_1a_2\ldots a_{4n}$ and $s_i$ be the label of square $a_i$ in a signature $s_1s_2\ldots s_{4n}$ of $C$ for $i = 1$, 2, \ldots, $4n$.

Suppose that $s_1 = +_\texttt{s}$. All other cases are analogous.

Since $s_1 = +_\texttt{s}$, the direction of the move to $a_1$ in $C$ is \texttt{1}, \texttt{3}, \texttt{5}, or \texttt{7}. Since, for $i = 1$, 2, \ldots, $4n - 1$, the direction of the move to $a_i$ in $C$ and the label $s_i$ of $a_i$ completely determine the direction of the move from $a_i$ to $a_{i + 1}$ in $C$, each of those four possibilities yields a unique possibility for the directions of all moves of $C$.

The labeling rules in the definition of a signature are invariant under a $90^\circ$ rotation. By Theorem \ref{symm}, it follows that $s_i = s_{i + n} = s_{i + 2n} = s_{i + 3n}$ for $i = 1$, 2, \ldots, $n$. Again by Theorem \ref{symm}, the moves to $a_1$, $a_{n + 1}$, $a_{2n + 1}$, and $a_{3n + 1}$ in $C$ are copies of each other under multiple-of-quarter-turn rotations. Therefore, their directions are \texttt{1}, \texttt{3}, \texttt{5}, and \texttt{7}.

It follows that all four possibilities for the directions of all moves of $C$ are cyclic shifts of each other and yield the same direction graph. Hence the following theorem.

\medskip

\begin{theorem} Let $L$ be a skew basic $(p, q)$-leaper, $C$ a cycle of $L$ within a $(p + q) \times (p + q)$ chessboard, and $e$ the descent of $C$. Then $C$ is trivial and the direction graph extracted from $C$ depends only on $e$. Equivalently, let $C$ be a cycle of type $k$. Then the direction graph extracted from $C$ depends only on the $k$-th convergent $[c_1\pm, c_2\pm, \ldots, c_k]$ of the even continued fraction representation of $\frac{q}{p}$. \label{fund} \end{theorem}

\medskip

Theorem \ref{fund} provides the basis for the following definition.

\medskip

\begin{definition} The \emph{fundamental direction cycle} $\Phi(e)$ of descent $e$ is the direction graph extracted from a cycle of descent $e$ of a skew basic $(p, q)$-leaper within a $(p + q) \times (p + q)$ chessboard. Equivalently, the fundamental direction cycle $\Phi[c_1\pm, c_2\pm, \ldots, c_k]$ of type $[c_1\pm,\allowbreak c_2\pm,\allowbreak \ldots,\allowbreak c_k]$ is the direction graph extracted from a cycle of type $k$ of a skew basic $(p, q)$-leaper within a $(p + q) \times (p + q)$ chessboard, where $[c_1\pm, c_2\pm, \ldots, c_k]$ is the $k$-th convergent of the even continued fraction representation of $\frac{q}{p}$. \end{definition}

\medskip

The two notations for a fundamental direction cycle are related as follows. The fundamental direction cycles of descent $e$ and type $[c_1\pm, c_2\pm, \ldots, c_k]$ coincide if $e$ admits a partitioning into (possibly empty) substrings of the form $e = e'_1e''_1e'_2e''_2\ldots e''_{k - 1}e'_k$, where $e'_i$ is a run of the character \texttt{f} of length $\frac{c_i}{2} - 1$ for $i = 1$, $2$, \ldots, $k$ and $e''_i$ consists of a single character, \texttt{g} if the sign following $c_i$ is $-$ and \texttt{h} if it is $+$, for $i = 1$, $2$, \ldots, $k - 1$.

We go on to second leaper cycles.

\medskip

\begin{definition} Let $L$ be a $(p, q)$-leaper with $p + q$ odd and $C$ a cycle of $L$ within a $(p + q) \times (p + q)$ chessboard. Then the \emph{canonical second leaper} associated with $C$ and the \emph{canonical second leaper cycle} over the squares of $C$ are the ones constructed in the proof of Theorem \ref{sl}. \end{definition}

\medskip

Let $D$ be a canonical second leaper cycle of a canonical second leaper $M$ over the squares of $C$.

Then $D$ is the product of successively applying a series of lifting transformations to a canonical second leaper cycle $D'$ of $M$ within the shell of a frame $F'$. If $F'$ is a $(1, 2)$-frame, then $M$ is a $(0, 1)$-leaper and no direction graph can be extracted from $D$. Otherwise, $F'$ is the product of lifting a smaller frame by means of either $g$ or $h$.

\medskip

\begin{definition} Let $L$ be a skew basic $(p, q)$-leaper, $D$ a canonical second leaper cycle of a skew basic canonical second leaper over the squares of a cycle of $L$ within a $(p + q) \times (p + q)$ chessboard, and $e = e_1e_2\ldots e_l$ a prefix of the descent $e_1e_2\ldots e_m$, $l < m$, of $L$ such that $D$ is the product of successively applying lifting transformations of types $e_l$, $e_{l - 1}$, \ldots, $e_1$ to a shell canonical second leaper cycle. Then $e$ is the \emph{descent} of $D$ and $e_{l + 1}$ is the \emph{origin} of $D$. \end{definition}

\medskip

The origin of a canonical second leaper cycle of a skew basic canonical second leaper is always \texttt{g} or \texttt{h}.

When $D$ is of origin \texttt{g}, the direction graph extracted from $D'$ is represented by the oriented cycle labeled \texttt{25476183}. For instance, such is the case with the $(1, 2)$-cycle in Figure \ref{2355}.

When $D$ is of origin \texttt{h}, the direction graph extracted from $D'$ is represented by the oriented cycle labeled \texttt{34567812}. For instance, such is the case with the $(1, 2)$-cycle in Figure \ref{2577}.

Impose an orientation on $D$. By induction on descent following the proofs of Lemmas \ref{f2}, \ref{g2}, and \ref{h2}, two skew directions occur as the directions of the moves to and from a square in $D$ if and only if they occur as the directions of the moves to and from a square in one of the two orientations of $D'$.

\medskip

\begin{definition} Let $L$ be a skew basic $(p, q)$-leaper and $D$ a canonical second leaper cycle of a skew basic canonical second leaper over the squares of a cycle of $L$ within a $(p + q) \times (p + q)$ chessboard.

Impose an orientation on $D$ and label every square $a$ of $D$ $+_\texttt{s}$, $+_\texttt{c}$, $-_\texttt{s}$, or $-_\texttt{c}$ as follows. Refer to the table \[ \begin{tabular}{c c} Label of $a$ & Directions of moves to and from $a$\\ \hline $+_\texttt{s}$ & \texttt{25}, \texttt{47}, \texttt{61}, \texttt{83}\\  $+_\texttt{c}$ & \texttt{18}, \texttt{32}, \texttt{54}, \texttt{76} \\ $-_\texttt{s}$ & \texttt{16}, \texttt{38}, \texttt{52}, \texttt{74} \\ $-_\texttt{c}$ & \texttt{23}, \texttt{45}, \texttt{67}, \texttt{81} \end{tabular} \] if $D$ is of origin \texttt{g}, and to the table \[ \begin{tabular}{c c} Label of $a$ & Directions of moves to and from $a$\\ \hline $+_\texttt{s}$ & \texttt{12}, \texttt{34}, \texttt{56}, \texttt{78}\\  $+_\texttt{c}$ & \texttt{23}, \texttt{45}, \texttt{67}, \texttt{81} \\ $-_\texttt{s}$ & \texttt{21}, \texttt{43}, \texttt{65}, \texttt{87} \\ $-_\texttt{c}$ & \texttt{18}, \texttt{32}, \texttt{54}, \texttt{76} \end{tabular} \] if $D$ is of origin \texttt{h}.

Then the cyclic string formed by the vertex labels of $D$ is a \emph{signature} of $D$. \end{definition}

\medskip

For instance, $+_\texttt{s}+_\texttt{c}+_\texttt{s}+_\texttt{c}+_\texttt{s}+_\texttt{c}+_\texttt{s}+_\texttt{c}$ is a signature of every shell canonical second leaper cycle, such as the $(1, 2)$-cycle in Figure \ref{2355} and the $(1, 2)$-cycle in Figure \ref{2577}.

\medskip

\begin{theorem} Let $L$ be a $(p, q)$-leaper with $p + q$ odd and $D$ a canonical second leaper cycle over the squares of a cycle of $L$ within a $(p + q) \times (p + q)$ chessboard. Then $D$ possesses orientation-preserving fourfold rotational symmetry together with orientation-reversing axial symmetry along a vertical, a horizontal, and two diagonal axes. Furthermore, all of the aforementioned symmetries preserve the partitioning of $D$ into eight disjoint nonempty paths, one within each section of the associated frame, as in the definition of a proper cycle. \label{secondsymm} \end{theorem}

\medskip

\begin{proof} By induction on descent. \end{proof}

\medskip

Analogously to the case of an $L$-cycle, it follows from Theorem \ref{secondsymm} that $D$ is trivial and that the signature of $D$ completely determines the direction graph extracted from $D$.

Let us track how a signature evolves when we lift a cycle.

Given a signature $s$ of $D$, partition $s$ into eight substrings $s^\texttt{E}$, $s^\texttt{NE}$, \ldots, $s^\texttt{SE}$ so that, for each direction $i$ out of \texttt{E}, \texttt{NE}, \ldots, \texttt{SE}, $s^i$ is the string formed by the vertex labels of the portion of $D$ within the section $F_i$ of the associated frame $F$.

By Theorem \ref{secondsymm}, $s^\texttt{E} = s^\texttt{N} = s^\texttt{W} = s^\texttt{S}$, $s^\texttt{NE} = s^\texttt{NW} = s^\texttt{SW} = s^\texttt{SE}$, and each of the eight strings $s^\texttt{E}$, $s^\texttt{NE}$, \ldots, $s^\texttt{SE}$ is a palindrome.

Therefore, a signature $s$ of $D$ is completely determined by its \emph{section pair}, the ordered pair $(s^\texttt{NE}, s^\texttt{E})$.

Let \begin{align*} \overline{+_\texttt{s}} &= -_\texttt{s} & \overline{+_\texttt{c}} &= -_\texttt{c}\\ \overline{-_\texttt{s}} &= +_\texttt{s} & \overline{-_\texttt{c}} &= +_\texttt{c} \end{align*} and \begin{align*} \overline{\overline{+_\texttt{s}}} &= -_\texttt{c} & \overline{\overline{+_\texttt{c}}} &= -_\texttt{s}\\ \overline{\overline{-_\texttt{s}}} &= +_\texttt{c} & \overline{\overline{-_\texttt{c}}} &= +_\texttt{s}. \end{align*} Given a string $w = w_1w_2\ldots w_n$ composed of the characters $+_\texttt{s}$, $+_\texttt{c}$, $-_\texttt{s}$, and $-_\texttt{c}$, we write $\overline{w}$ for the string $\overline{w_1}\,\overline{w_2}\ldots \overline{w_n}$, $[w$ for the string $\overline{\overline{w_1}}w_2\ldots w_n$, and $w]$ for the string $w_1w_2\ldots w_{n - 1}\overline{\overline{w_n}}$.

Lift $D$ by means of $f$ to $f(D)$. By the proof of Lemma \ref{f2}, \[(s^\texttt{NE},\; [\overline{s^\texttt{NE}}\,\overline{s^\texttt{E}}\,\overline{s^\texttt{NE}}])\] is the section pair of a signature of $f(D)$.

This transformation and its analogues for $g$ and $h$ suffice to establish Theorem \ref{secondfund}, but they are not well-suited to a proof of Theorem \ref{flip}. For this reason, we introduce a change of variables.

If $|s^\texttt{E}| \le |s^\texttt{NE}|$, then partition $s^\texttt{NE}$ into three substrings $s^\texttt{NE} = s^\text{Left}s^\text{Corner}s^\text{Right}$ such that $|s^\text{Left}| = |s^\text{Right}|$ and $|s^\text{Corner}| = |s^\texttt{E}|$, and set $s^\text{Side} = s^\text{Right}s^\texttt{E}s^\text{Left}$.

If $|s^\texttt{E}| \ge |s^\texttt{NE}|$, then partition $s^\texttt{E}$ into three substrings $s^\texttt{E} = s^\text{Left}s^\text{Side}s^\text{Right}$ such that $|s^\text{Left}| = |s^\text{Right}|$ and $|s^\text{Side}| = |s^\texttt{NE}|$, and set $s^\text{Corner} = s^\text{Right}s^\texttt{NE}s^\text{Left}$.

In both cases, $s$ equals, up to a cyclic shift, the concatenation $s^\text{Corner}s^\text{Side}\allowbreak s^\text{Corner}s^\text{Side}\allowbreak s^\text{Corner}s^\text{Side}\allowbreak s^\text{Corner}s^\text{Side}$. We refer to the ordered pair $(s^\text{Corner}, s^\text{Side})$ as the \emph{corner-side pair} of $s$.

We proceed to show that \[(s^\text{Side}s^\text{Corner}s^\text{Side},\; \overline{s^\text{Side}})\] is the corner-side pair of a signature of $f(D)$. We consider the case $|s^\texttt{E}| \le |s^\texttt{NE}|$ in detail, and the opposite case is analogous.

Since \[s^\texttt{NE} = s^\text{Left}s^\text{Corner}s^\text{Right}\] and \begin{align*} [\overline{s^\texttt{NE}}\,\overline{s^\texttt{E}}\,\overline{s^\texttt{NE}}] &= [\overline{s^\text{Left}}\,\overline{s^\text{Corner}}\,\overline{s^\text{Right}}\,\overline{s^\texttt{E}}\,\overline{s^\text{Left}}\,\overline{s^\text{Corner}}\,\overline{s^\text{Right}}] =\\ &= [\overline{s^\text{Left}}\,\overline{s^\text{Corner}}\,\overline{s^\text{Side}}\,\overline{s^\text{Corner}}\,\overline{s^\text{Right}}], \end{align*} a corner-side pair of a signature of $f(D)$ is \[(\overline{s^\text{Corner}}\,\overline{s^\text{Right}}]s^\text{Left}s^\text{Corner}s^\text{Right}[\overline{s^\text{Left}}\,\overline{s^\text{Corner}},\; \overline{s^\text{Side}}).\]

However, since $D$ is a proper cycle, \[\overline{s^\text{Corner}}\,\overline{s^\text{Right}}] = s^\text{Right}s^\texttt{E}\] and \[[\overline{s^\text{Left}}\,\overline{s^\text{Corner}} = s^\texttt{E}s^\text{Left}.\]

Therefore, we can rewrite the above as \[(s^\text{Right}s^\texttt{E}s^\text{Left}s^\text{Corner}s^\text{Right}s^\texttt{E}s^\text{Left},\; \overline{s^\text{Side}}),\] as needed.

\begin{figure}[ht!]\hspace*{\fill}\includegraphics[scale=0.66]{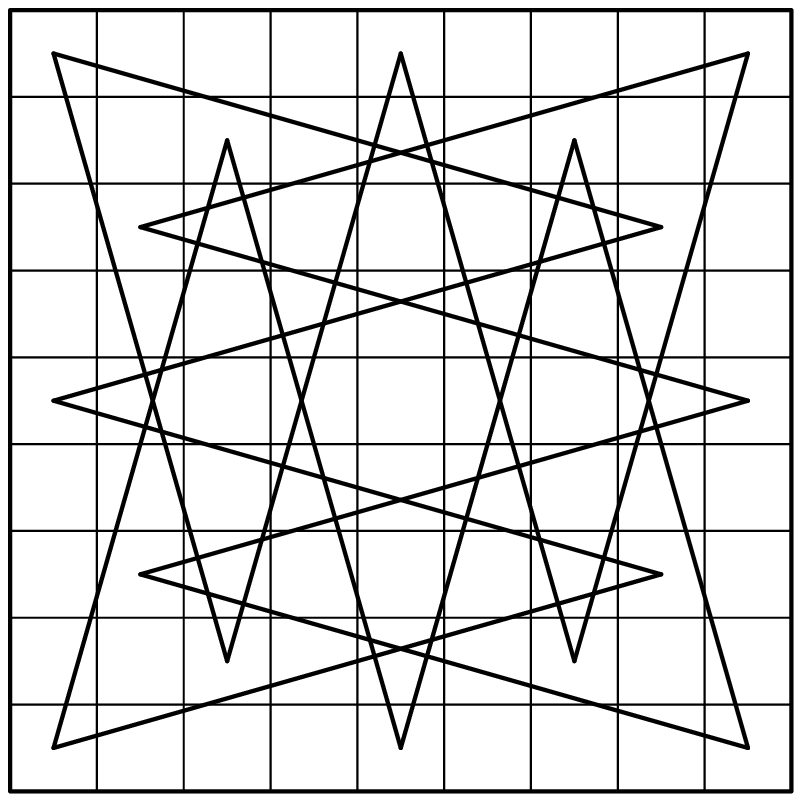}\hfill\includegraphics[scale=0.66]{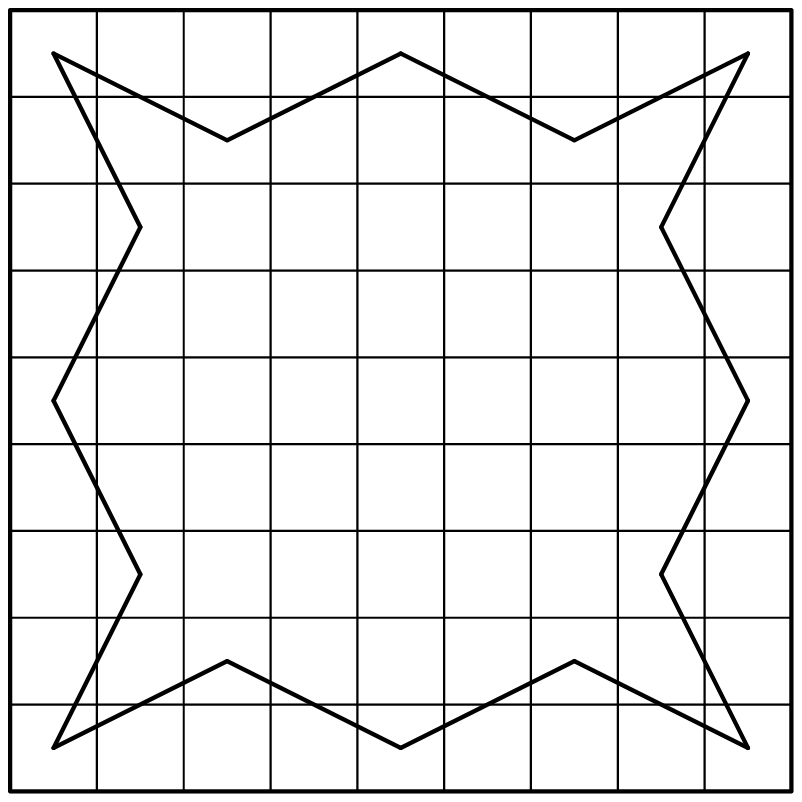}\hspace*{\fill}\caption{}\label{2799}\end{figure}

We refer to the transformation that maps $(s^\text{Corner}, s^\text{Side})$ to $(s^\text{Side}s^\text{Corner}s^\text{Side},\allowbreak \overline{s^\text{Side}})$ as an \emph{$f$-rearrangement}. For instance, the $f$-rearrangement of $(+_\texttt{s}, +_\texttt{c})$ is $(+_\texttt{c}+_\texttt{s}+_\texttt{c}, -_\texttt{c})$, and $+_\texttt{c}+_\texttt{s}+_\texttt{c}-_\texttt{c}+_\texttt{c}+_\texttt{s}+_\texttt{c}-_\texttt{c}+_\texttt{c}+_\texttt{s}+_\texttt{c}-_\texttt{c}+_\texttt{c}+_\texttt{s}+_\texttt{c}-_\texttt{c}$ is a signature of the $(1, 2)$-cycle in Figure \ref{2799}.

Next lift $D$ by means of $g$ to $g(D)$. Analogously, by the proof of Lemma \ref{g2}, \[(\overline{s^\text{Corner}},\; s^\text{Corner}s^\text{Side}s^\text{Corner})\] is the corner-side pair of a signature of $g(D)$.

We refer to this transformation as a \emph{$g$-rearrangement}. For instance, the $g$-rearrangement of $(+_\texttt{s}, +_\texttt{c})$ is $(-_\texttt{s}, +_\texttt{s}+_\texttt{c}+_\texttt{s})$, and $-_\texttt{s}+_\texttt{s}+_\texttt{c}+_\texttt{s}-_\texttt{s}+_\texttt{s}+_\texttt{c}+_\texttt{s}-_\texttt{s}+_\texttt{s}+_\texttt{c}+_\texttt{s}-_\texttt{s}+_\texttt{s}+_\texttt{c}+_\texttt{s}$ is a signature of the $(1, 2)$-cycle in Figure \ref{3477}.

Lastly, lift $D$ by means of $h$ to $h(D)$. Analogously, by the proof of Lemma \ref{h2}, \[(s^\text{Corner}s^\text{Side}s^\text{Corner}s^\text{Side}s^\text{Corner},\; \overline{s^\text{Corner}}\,\overline{s^\text{Side}}\,\overline{s^\text{Corner}})\] is the corner-side pair of a signature of $h(D)$.

\begin{figure}[htp!]\hspace*{\fill}\includegraphics[scale=0.66]{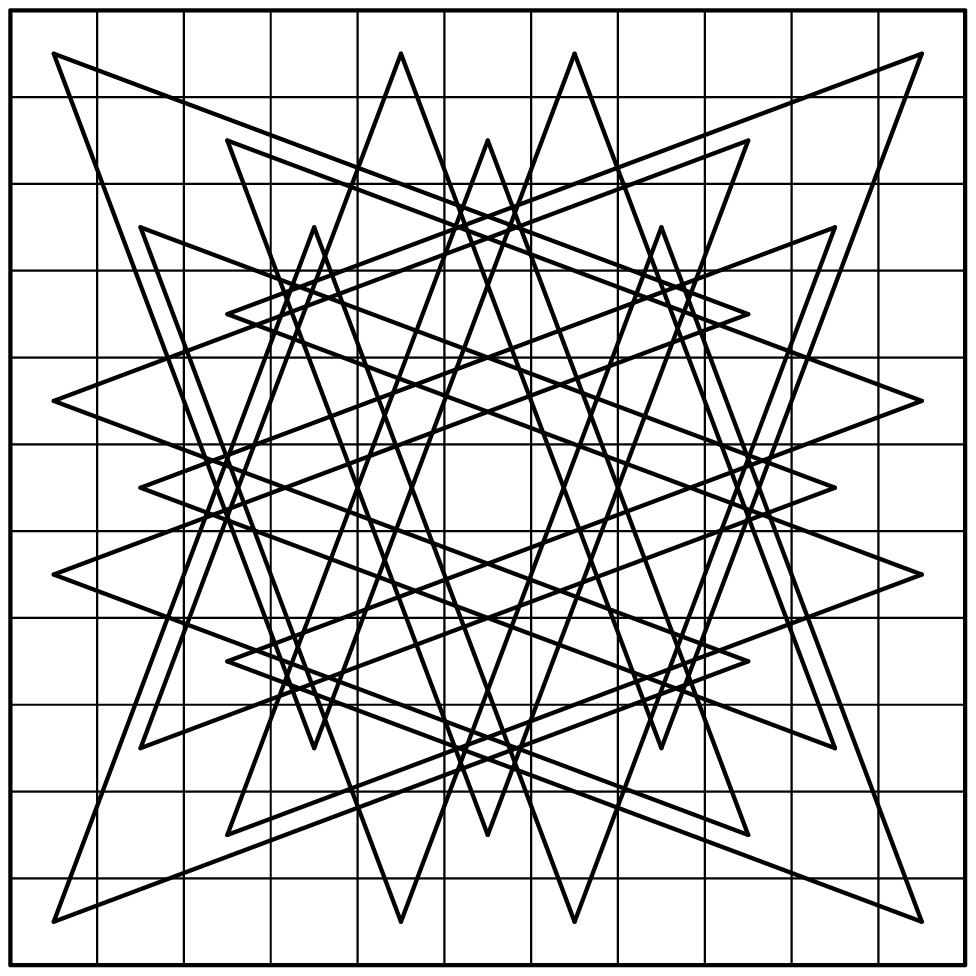}\hspace*{\fill}\newline\bigskip\newline\hspace*{\fill}\includegraphics[scale=0.66]{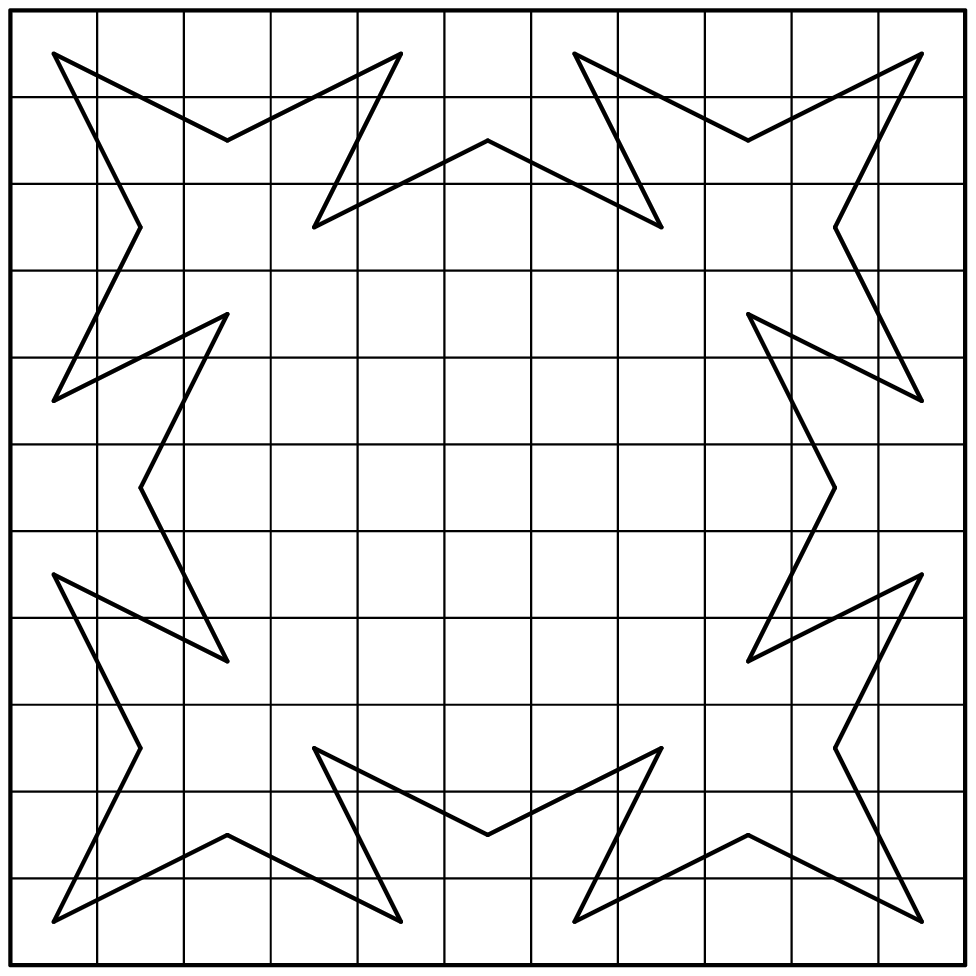}\hspace*{\fill}\caption{}\label{381111}\end{figure}

We refer to this transformation as an \emph{$h$-rearrangement}. For instance, the $h$-rearrangement of $(+_\texttt{s}, +_\texttt{c})$ is $(+_\texttt{s}+_\texttt{c}+_\texttt{s}+_\texttt{c}+_\texttt{s}, -_\texttt{s}-_\texttt{c}-_\texttt{s})$, and $+_\texttt{s}+_\texttt{c}+_\texttt{s}+_\texttt{c}+_\texttt{s}-_\texttt{s}-_\texttt{c}-_\texttt{s}+_\texttt{s}+_\texttt{c}+_\texttt{s}+_\texttt{c}+_\texttt{s}-_\texttt{s}-_\texttt{c}-_\texttt{s}+_\texttt{s}+_\texttt{c}+_\texttt{s}+_\texttt{c}+_\texttt{s}-_\texttt{s}-_\texttt{c}-_\texttt{s}+_\texttt{s}+_\texttt{c}+_\texttt{s}+_\texttt{c}+_\texttt{s}-_\texttt{s}-_\texttt{c}-_\texttt{s}$ is a signature of the $(1, 2)$-cycle in Figure \ref{381111}.

Given a string $e = e_1e_2\ldots e_l$ composed of the characters \texttt{f}, \texttt{g}, and \texttt{h} and an ordered pair of strings $(s', s'')$ composed of the characters $+_\texttt{s}$, $+_\texttt{c}$, $-_\texttt{s}$, and $-_\texttt{c}$, we write $W_e(s', s'')$ for the ordered pair of strings obtained by successively applying rearrangements of types $e_l$, $e_{l - 1}$, \ldots, $e_1$ to $(s', s'')$.

Let $e$ be the descent of $D$. Then \[W_e(+_\texttt{s}, +_\texttt{c})\] is the corner-side pair of a signature of $D$. Hence the following theorem.

\medskip

\begin{theorem} Let $L$ be a skew basic $(p, q)$-leaper, $D$ a canonical second leaper cycle of a skew basic canonical second leaper over the squares of a cycle $C$ of $L$ within a $(p + q) \times (p + q)$ chessboard, $e$ the descent of $D$, and $o$ the origin of $D$. Then $D$ is trivial and the direction graph extracted from $D$ depends only on $e$ and $o$. Equivalently, let $C$ be a cycle of type $k$. Then the direction graph extracted from $D$ depends only on the $k$-th convergent $[c_1\pm, c_2\pm, \ldots, c_k]$ of the even continued fraction representation of $\frac{q}{p}$ and the sign following $c_k$ in that representation. \label{secondfund} \end{theorem}

\medskip

Theorem \ref{secondfund} provides the basis for the following definition.

\medskip

\begin{definition} The \emph{second fundamental direction cycle} $\Phi^\text{II}_o(e)$ of descent $e$ and origin $o$ is the direction graph extracted from a canonical second leaper cycle of descent $e$ and origin $o$. Equivalently, the second fundamental direction cycle $\Phi^\text{II}_\epsilon[c_1\pm, c_2\pm, \ldots, c_k]$ of type $[c_1\pm,\allowbreak c_2\pm,\allowbreak \ldots,\allowbreak c_k]$ and sign $\epsilon$ is the direction graph extracted from a canonical second leaper cycle over the squares of a cycle of type $k$ of a skew basic $(p, q)$-leaper within a $(p + q) \times (p + q)$ chessboard, where $[c_1\pm, c_2\pm, \ldots, c_k]$ is the $k$-th convergent of the even continued fraction representation of $\frac{q}{p}$ and the sign following $c_k$ in that representation is $\epsilon$. \end{definition}

\medskip

The two notations for a second fundamental direction cycle are related as follows. The direction graphs $\Phi^\text{II}_o(e)$ and $\Phi^\text{II}_\epsilon[c_1\pm, c_2\pm, \ldots, c_k]$ coincide if $e$ and $[c_1\pm, c_2\pm, \ldots, c_k]$ are related as in the case of an $L$-cycle, and $\epsilon = -$ if $o = \texttt{g}$ and $\epsilon = +$ if $o = \texttt{h}$.

When $D$ is a cycle of a $(0, 1)$-leaper, Theorem \ref{secondfund} continues to apply in the sense that $D$ is an instantiation of both $\Phi^\text{II}_\texttt{g}(e)$ and $\Phi^\text{II}_\texttt{h}(e)$.

The equivalence permutation \[\pi^\text{II} = \left(\begin{array}{cccccccc} \texttt{1} & \texttt{2} & \texttt{3} & \texttt{4} & \texttt{5} & \texttt{6} & \texttt{7} & \texttt{8}\\ \texttt{4} & \texttt{7} & \texttt{6} & \texttt{1} & \texttt{8} & \texttt{3} & \texttt{2} & \texttt{5} \end{array}\right),\] induced by the unit $2 \times 2$ matrix and $\left( \begin{smallmatrix} 1 & 0\\ 0 & -1 \end{smallmatrix} \right)$, maps the signature labeling rules for a canonical second leaper cycle of origin \texttt{g} to the signature labeling rules for a canonical second leaper cycle of origin \texttt{h}. Since $\Phi^\text{II}_\texttt{g}(e)$ and $\Phi^\text{II}_\texttt{h}(e)$ are determined by a common signature, it follows that they are equivalent under $\pi^\text{II}$.

We conclude this section by showing that the family of all fundamental direction cycles and the family of all second fundamental direction cycles are essentially the same family.

\medskip

\begin{definition} Let $\overline{\texttt{f}} = \texttt{g}$, $\overline{\texttt{g}} = \texttt{f}$, and $\overline{\texttt{h}} = \texttt{h}$. The \emph{flip} of a string $e = e_1e_2\ldots e_l$ composed of the characters \texttt{f}, \texttt{g}, and \texttt{h} is the string $\overline{e_l}\,\overline{e_{l - 1}}\ldots \overline{e_1}$. \end{definition}

\medskip

A string $e'$ is the flip of $e''$ if and only if $e''$ is the flip of $e'$.

\medskip

\begin{theorem} Let the strings $e'$ and $e''$ composed of the characters \texttt{f}, \texttt{g}, and \texttt{h} be flips of each other. Then the fundamental direction cycle of descent $e'$ and the two second fundamental direction cycles of descent $e''$ are equivalent. \label{flip} \end{theorem}

\medskip

In other words, there exists a one-to-two equivalence mapping between fundamental direction cycles and second fundamental direction cycles.

\medskip

\begin{proof} The equivalence permutation \[\pi_\texttt{g} = \left(\begin{array}{cccccccc} \texttt{1} & \texttt{2} & \texttt{3} & \texttt{4} & \texttt{5} & \texttt{6} & \texttt{7} & \texttt{8}\\ \texttt{2} & \texttt{7} & \texttt{8} & \texttt{5} & \texttt{6} & \texttt{3} & \texttt{4} & \texttt{1} \end{array}\right),\] induced by the matrices \[P_\texttt{g} = \left(\begin{array}{cc} \frac{1}{\sqrt{2}} & \frac{1}{\sqrt{2}}\\[5pt] \frac{1}{\sqrt{2}} & -\;\frac{1}{\sqrt{2}} \end{array}\right)\] and \[Q_\texttt{g} = \left(\begin{array}{cc} \frac{1}{\sqrt{2}} & -\;\frac{1}{\sqrt{2}}\\[5pt] \frac{1}{\sqrt{2}} & \frac{1}{\sqrt{2}} \end{array}\right),\] maps the signature labeling rules for an $L$-cycle to the signature labeling rules for a canonical second leaper cycle of origin \texttt{g}.

The equivalence permutation \[\pi_\texttt{h} = \left(\begin{array}{cccccccc} \texttt{1} & \texttt{2} & \texttt{3} & \texttt{4} & \texttt{5} & \texttt{6} & \texttt{7} & \texttt{8}\\ \texttt{7} & \texttt{2} & \texttt{5} & \texttt{8} & \texttt{3} & \texttt{6} & \texttt{1} & \texttt{4} \end{array}\right),\] induced by the matrices \[P_\texttt{h} = Q_\texttt{h} = \left(\begin{array}{cc} \frac{1}{\sqrt{2}} & \frac{1}{\sqrt{2}}\\[5pt] \frac{1}{\sqrt{2}} & -\;\frac{1}{\sqrt{2}} \end{array}\right),\] maps the signature labeling rules for an $L$-cycle to the signature labeling rules for a canonical second leaper cycle of origin \texttt{h}.

(Furthermore, $\pi^\text{II} \circ \pi_\texttt{g} = \pi_\texttt{h}$.)

We are left to show that $\Phi(e')$, $\Phi^\text{II}_\texttt{g}(e'')$, and $\Phi^\text{II}_\texttt{h}(e'')$ are determined by a common signature.

A signature of $\Phi(e')$ is \[R_{e'}(+_\texttt{s}+_\texttt{c}+_\texttt{s}+_\texttt{c}+_\texttt{s}+_\texttt{c}+_\texttt{s}+_\texttt{c}) = R_{e'}(+_\texttt{s})R_{e'}(+_\texttt{c})\ldots R_{e'}(+_\texttt{c}).\]

Therefore, it suffices to show that \[(R_{e'}(+_\texttt{s}), R_{e'}(+_\texttt{c})) = W_{e''}(+_\texttt{s}, +_\texttt{c}).\]

We proceed by induction on descent.

When $e'$ and $e''$ are both the empty string, equality holds.

Suppose, then, that equality holds for $e'$ and $e''$. We need to show that equality holds for $e'\texttt{f}$ and $\texttt{g}e''$, $e'\texttt{g}$ and $\texttt{f}e''$, and $e'\texttt{h}$ and $\texttt{h}e''$. We consider the case of $e'\texttt{f}$ and $\texttt{g}e''$ in detail, and all other cases are analogous.

For all strings $e$ composed of the characters \texttt{f}, \texttt{g}, and \texttt{h} and all strings $s$ composed of the characters $+_\texttt{s}$, $+_\texttt{c}$, $-_\texttt{s}$, and $-_\texttt{c}$, \[R_e(\overline{s}) = \overline{R_e(s)}.\]

Therefore, \begin{align*} (R_{e'\texttt{f}}(+_\texttt{s}),\; R_{e'\texttt{f}}(+_\texttt{c})) &= (R_{e'}(R_\texttt{f}(+_\texttt{s})),\; R_{e'}(R_\texttt{f}(+_\texttt{c})))\\ &= (R_{e'}(-_\texttt{s}),\; R_{e'}(+_\texttt{s}+_\texttt{c}+_\texttt{s}))\\ &= (\overline{R_{e'}(+_\texttt{s})},\; R_{e'}(+_\texttt{s})R_{e'}(+_\texttt{c})R_{e'}(+_\texttt{s}))\\ &= W_\texttt{g}(R_{e'}(+_\texttt{s}),\; R_{e'}(+_\texttt{c}))\\ &= W_\texttt{g}(W_{e''}(+_\texttt{s}, +_\texttt{c}))\\ &= W_{\texttt{g}e''}(+_\texttt{s}, +_\texttt{c}). \end{align*}

This completes the proof of the theorem. \end{proof}

\section{Dual Boards and Dual Direction Graphs}

Theorem \ref{sl} can be strengthened as follows.

\medskip

\begin{theorem} Let $L$ be a $(p, q)$-leaper with $p + q$ odd, $C$ a cycle of $L$ within a $(p + q) \times (p + q)$ chessboard, $M$ the canonical second leaper associated with $C$, and $D$ the canonical second leaper cycle of $M$ over the squares of $C$. Then $D$ is the leaper graph of $M$ over the squares of $C$. \label{dual} \end{theorem}

\medskip

\begin{proof} It suffices to consider the case when $L$ is a skew basic leaper and $p < q$. Let $F$ be the associated $(p, q)$-frame and $r$ and $s$, $r < s$, the proportions of $M$. The proof proceeds by induction on descent.

When $D$ is a shell cycle, the theorem holds.

Suppose, then, that the theorem holds for $D$.

\medskip

\begin{lemma} The proportions of $L$ and $M$ satisfy $s \le p$. Furthermore, $s \le q - p$ unless the origin of $D$ is \texttt{g} and the descent of $D$ is a run of the character \texttt{g}. \label{sp} \end{lemma}

\medskip

\begin{proof} By induction on descent. \end{proof}

\medskip

\emph{Case $f$.} Lift $F$ and $D$ to $H$ and $E$ by means of $f$. For each direction $i$ out of \texttt{E}, \texttt{NE}, \ldots, \texttt{SE}, let $S_i$ be the subboard $(F + v^f_i) \cap H$ of $H$.

Let $a$ and $b$ be two squares in $E$ joined by a move of $M$. By Lemma \ref{sp}, $s \le p$. Therefore, there exists a direction $i$ out of \texttt{E}, \texttt{NE}, \ldots, \texttt{SE} such that both $a$ and $b$ belong to $S_i$. By the fact that $D$ is a proper cycle and the definition of $E$, $a$ and $b$ are the images under $v^f_i$ of two squares $a'$ and $b'$ of $D$.

Since $a'$ and $b'$ are joined by a move of $M$, by the induction hypothesis they are adjacent in $D$. Therefore, by the fact that $D$ is a proper cycle and the definition of $E$, $a$ and $b$ are adjacent in $E$.

\medskip

\emph{Case $g$.} Lift $F$ and $D$ to $H$ and $E$ by means of $g$. For each direction $i$ out of \texttt{E}, \texttt{N}, \texttt{W}, and \texttt{S}, let $S_i$ be the subboard $(F + v^g_i) \cap (H_{i - 1} \cup H_i \cup H_{i + 1})$ of $H$, and, for each direction $i$ out of \texttt{NE}, \texttt{NW}, \texttt{SW}, and \texttt{SE}, let $S_i$ be the subboard $(F + v^g_i) \cap H$ of $H$.

If $s \le q - p$, the proof continues as in Case $f$.

Otherwise, by Lemma \ref{sp} the origin of $D$ is \texttt{g} and the descent of $D$ is a run of the character \texttt{g}.

Let $O$ be the center of symmetry of $D$. Introduce a Cartesian coordinate system $Oxy$ over the infinite chessboard and write $(x, y)$ for the square centered at $(x, y)$.

Let $l$ the length of the descent of $D$. Then the set $S'$ of all squares $((s - r)x', (s - r)y')$ such that $x'$ and $y'$ are integers and $|x'| + |y'| = l + 2$, and the set $S''$ of all squares $(\pm[(s - r)x'' + r], \pm[(s - r)y'' + r])$ such that $x''$ and $y''$ are nonnegative integers and $x'' + y'' = l + 1$, form a partitioning of the vertex set of $E$.

Suppose first that both squares $a$ and $b = a + v$ belong to $S'$. Then both coordinates of the translation $v$ are even and, since $M$ is a basic leaper, $v$ cannot be an $M$-translation. Analogous reasoning applies to the case when both $a$ and $b$ belong to $S''$.

Suppose, then, that $a$ belongs to $S'$. Then all but two of the squares of $S''$ lie outside of the $(2s + 1) \times (2s + 1)$ chessboard centered at $a$. Therefore, every square in $S'$ is joined by a move of $M$ to at most two squares in $S''$.

It follows that two squares of $E$ are joined by a move of $M$ if and only if they are adjacent in $E$.

\medskip

\emph{Case $h$.} Lift $F$ and $D$ to $H$ and $E$ by means of $h$. For each direction $i$ out of \texttt{E}, \texttt{NE}, \ldots, \texttt{SE}, let $S_i$ be the subboard $(F + v^h_i) \cap H$ of $H$. The proof continues as in Case $f$. \end{proof}

\medskip

In other words, in the setting of Theorem \ref{dual}, the leaper graphs of $L$ and $M$ over the board formed by the squares of $C$ are isomorphic.

\medskip

\begin{definition} A board $B$ is \emph{dual} with respect to two distinct leapers $L$ and $M$ if $B$ contains more than one square and the leaper graphs of $L$ and $M$ over $B$ are connected and isomorphic. \end{definition}

\medskip

The notion of a dual board raises a number of questions.

\medskip

\begin{question} Given two distinct leapers $L$ and $M$, does there exist a board dual with respect to $L$ and $M$? \label{lmboard} \end{question}

\medskip

A necessary condition is that $L$ and $M$ are obtained from two basic leapers by means of the same scaling and rotation. Therefore, it suffices to study the case of both of $L$ and $M$ being basic leapers.

By Theorems \ref{descent} and \ref{dual} and Corollary \ref{fghboard}, a sufficient condition for basic $L$ and $M$ is that either one of them is a $(0, 1)$-leaper or both of them are skew leapers and the descent of one of them is a suffix of the descent of the other.

\medskip

\begin{question} Given two distinct leapers $L$ and $M$, does there exist a polyomino board dual with respect to $L$ and $M$? \label{polyomino} \end{question}

\medskip

When one of $L$ and $M$ is a $(0, 1)$-leaper, this question is equivalent to the previous one.

\medskip

\begin{question} Given two distinct leapers $L$ and $M$ such that a board dual with respect to $L$ and $M$ does exist, what is the least number of squares that it may contain? Is the number of squares that it may contain unbounded from above? \label{boardbound} \end{question}

\medskip

A dual board given by Theorem \ref{dual} is never minimal, as removing any square from it yields a board dual with respect to the same pair of leapers.

Corollary \ref{pinwheelunbound} answers one special case of the second part of the question.

\medskip

\begin{question} For what positive integers $n \ge 2$ does there exist a board dual with respect to $n$ pairwise distinct leapers? \label{multleap} \end{question}

\medskip

A strengthening of Theorem \ref{tl} analogous to Theorem \ref{dual} shows that, in the setting of Theorem \ref{tl}, the leaper graph of a $(1, 2)$-leaper over the squares of the unique cycle of type $k$ is a cycle if and only if $e_l = \texttt{g}$ and $l_k$ is indivisible by three. This yields an infinite family of boards dual with respect to three pairwise distinct leapers.

\medskip

Let $L$ be a skew basic $(p, q)$-leaper, $C$ a cycle of $L$ within a $(p, q)$-frame $F$, $M$ the canonical second leaper associated with $C$, and $D$ the canonical second leaper cycle of $M$ over the squares of $C$.

We go on to describe all isomorphisms between $C$ and $D$.

Let $D$ visit $m$ squares in each side section of $F$ and $n$ squares in each corner section of $F$, $D = a^\texttt{E}a^\texttt{NE}\ldots a^\texttt{SE}$ be a partitioning of $D$ into eight disjoint paths such that, for all directions $i$ out of \texttt{E}, \texttt{NE}, \ldots, \texttt{SE}, $a^i$ lies within $F_i$, $a^\texttt{E} = a_1a_2\ldots a_m$, $a^\texttt{NE} = a_{m + 1}a_{m + 2}\ldots a_{m + n}$, \ldots, $a^\texttt{SE} = a_{4m + 3n + 1}a_{4m + 3n + 2}\ldots\allowbreak a_{4(m + n)}$, and, for all integers $j$, $a_{j + 4(m + n)} \equiv a_j$.

Given an integer $x$ and a positive integer $y$, we write $x \% y$ for the remainder of $x$ upon division by $y$, with $0 \le x \% y \le y - 1$.

For $\epsilon = 0$, 1 and all integers $j$, let \[\upsilon(\epsilon, j) = m + n + 2[(\epsilon m - j) \% (m + n)] + 1.\]

Since $D$ is a proper cycle, for every integer $j$ the square $a_j$ is linked by a move of $L$ to both of the squares $a_{j + \upsilon(0, j)}$ and $a_{j + \upsilon(1, j)}$.

By Theorem \ref{seclen}, $m + n$ is even. Therefore, the parity of $j$ differs from the parity of $j + \upsilon(\epsilon, j)$ for $\epsilon = 0$, 1 and all integers $j$.

Furthermore, $a_{j'} = a_{j'' + \upsilon(\epsilon, j'')}$ if and only if $a_{j''} = a_{j' + \upsilon(\epsilon, j')}$, for $\epsilon = 0$, 1 and all integers $j'$ and $j''$.

It follows that drawing an arc from $a_j$ to $a_{j + \upsilon(0, j)}$ for all even $j$ and from $a_j$ to $a_{j + \upsilon(1, j)}$ for all odd $j$ yields an orientation of $C$.

Let $\psi(1) = 1$, \[\psi(j + 1) = \psi(j) + \upsilon(\psi(j) \% 2,\; \psi(j))\] for $j = 1$, 2, \ldots, $4(m + n) - 1$, and $\psi(j + 4(m + n)) = \psi(j)$ for all integers $j$.

Then the mapping \[a_j \to a_{\psi(j)}\] is an isomorphism from $D$ onto $C$, and all other such isomorphisms are obtained from it by means of rotation and reflection.

When $j$ takes on the values 1, 2, \ldots, $4(m + n)$, the difference $\psi(j + 1) - \psi(j) = \upsilon(\psi(j) \% 2,\; \psi(j))$ takes on the same values as $\upsilon(j \% 2, j)$. Those are all integers between $m + n$ and $3(m + n)$ congruent to $n - m - 1$ modulo four, each occurring precisely eight times. Hence the following theorem.

\medskip

\begin{theorem} Let $L$ be a $(p, q)$-leaper with $p + q$ odd, $C$ a cycle of $L$ within a $(p + q) \times (p + q)$ chessboard, $M$ the canonical second leaper associated with $C$, $D$ the canonical second leaper cycle of $M$ over the squares of $C$, and $4\mu$ the length of both $C$ and $D$.

Impose orientations on $C$ and $D$. Given a square $a$ of $D$, let the arc from $a$ in $C$ point to $b$ and the \emph{$D$-displacement} of $a$ be the number of steps along $D$ that lead from $a$ to $b$.

Then the multiset of all $D$-displacements consists of all integers between $\mu$ and $3\mu$ congruent to $\alpha$ modulo four, each occurring precisely eight times, where $\alpha = 1$ or $\alpha = 3$ depending on the orientations of $C$ and $D$. \label{seconddispl} \end{theorem}

\medskip

We proceed to look at the isomorphisms between $C$ and $D$ from the point of view of direction graphs.

Let $r$ and $s$, $r < s$, be the proportions of $M$. Suppose that $M$ is a skew basic leaper, and let $o$ be the origin of $D$ and $e$ the descent of $C$ and $D$.

By the proof of Theorem \ref{secondfund}, the skew direction of the move of $M$ from $a_j$ to $a_{j + 1}$ depends only on $o$, $e$, and $j$.

Therefore, there exists an enumeration \[a^\text{II}_1, a^\text{II}_2, \ldots, a^\text{II}_{4(m + n)}\] (with $a^\text{II}_{j + 4(m + n)} \equiv a^\text{II}_j$ for all integers $j$) of the vertices of the second fundamental direction cycle $\Phi^\text{II}_o(e)$, depending only on $o$ and $e$ and not on the remainder of the descent of $L$, such that $D$ is an instantiation of $\Phi^\text{II}_o(e)$ under the mapping $a^\text{II}_j \to a_j$.

The move of $L$ from $a_j$ to $a_{j + \upsilon(0, j)}$ is of skew direction \texttt{5} when $a_j$ belongs to $F_\texttt{E} \cup F_\texttt{NE}$, \texttt{7} when $a_j$ belongs to $F_\texttt{N} \cup F_\texttt{NW}$, and similarly for all other possible positions of $a_j$. Analogously, the move of $L$ from $a_j$ to $a_{j + \upsilon(1, j)}$ is of skew direction \texttt{6} when $a_j$ belongs to $F_\texttt{NE} \cup F_\texttt{N}$, \texttt{8} when $a_j$ belongs to $F_\texttt{NW} \cup F_\texttt{W}$, and similarly for all other possible positions of $a_j$.

Since $m$ and $n$ are completely determined by $e$, the section of $F$ that $a_j$ belongs to is completely determined by $e$ and $j$. Therefore, the skew direction of the move of $L$ from $a_j$ to $a_{j + \upsilon(\epsilon, j)}$ is completely determined by $e$, $j$, and $\epsilon$.

It follows that there exists an enumeration \[a^\text{I}_1, a^\text{I}_2, \ldots, a^\text{I}_{4(m + n)}\] (with $a^\text{I}_{j + 4(m + n)} \equiv a^\text{I}_j$ for all integers $j$) of the vertices of the fundamental direction cycle $\Phi(e)$, depending only on $e$ and not on the remainder of the descent of $L$, such that $C$ is an instantiation of $\Phi(e)$ under the mapping $a^\text{I}_j \to a_{\psi(j)}$.

Let \begin{align*}A_\texttt{f} &= \left(\begin{array}{cc} 1 & 0\\ 2 & 1 \end{array}\right),\\ A_\texttt{g} &= \left(\begin{array}{cc} 0 & 1\\ -1 & 2 \end{array}\right),\\ A_\texttt{h} &= \left(\begin{array}{cc} 0 & 1\\ 1 & 2 \end{array}\right),\end{align*} and, for all strings $e_1e_2\ldots e_l$ composed of the characters \texttt{f}, \texttt{g}, and \texttt{h}, \[A_{e_1e_2\ldots e_l} = A_{e_1}A_{e_2}\cdots A_{e_l}.\]

Then the proportions of $L$ and $M$ satisfy \[\left(\begin{array}{c} p\\ q \end{array}\right) = A_eA_o\left(\begin{array}{c} r\\ s \end{array}\right).\]

Let $A^\text{I}(j)$ be the direction matrix of the arc from $a^\text{I}_j$ to $a^\text{I}_{j + 1}$ in $\Phi(e)$, $\text{Dist}^\text{I}(j', j'')$ the distance from $a^\text{I}_{j'}$ to $a^\text{I}_{j''}$ in $\Phi(e)$, $A^\text{II}(j)$ the direction matrix of the arc from $a^\text{II}_j$ to $a^\text{II}_{j + 1}$ in $\Phi^\text{II}_o(e)$, and $\text{Dist}^\text{II}(j', j'')$ the distance from $a^\text{II}_{j'}$ to $a^\text{II}_{j''}$ in $\Phi^\text{II}_o(e)$, for all integers $j$, $j'$, and $j''$.

Then \[\text{Dist}^\text{II}(\psi(j), \psi(j + 1))\left(\begin{array}{c} r\\ s \end{array}\right) = A^\text{I}(j)\left(\begin{array}{c} p\\ q \end{array}\right),\] or, equivalently, \[\text{Dist}^\text{II}(\psi(j), \psi(j + 1))\left(\begin{array}{c} r\\ s \end{array}\right) = A^\text{I}(j)A_eA_o\left(\begin{array}{c} r\\ s \end{array}\right),\] for all integers $j$.

Since all objects in this identity apart from $r$ and $s$ depend only on $o$ and $e$, it continues to hold when $r$ and $s$ are replaced with the proportions $r'$ and $s'$, $r' < s'$, of any canonical second leaper. However, every basic leaper occurs as a canonical second leaper.

Therefore, \[\text{Dist}^\text{II}(\psi(j), \psi(j + 1)) = A^\text{I}(j)A_eA_o\] for all integers $j$.

We refer to the above identity as the \emph{duality identity} for $\Phi(e)$ and $\Phi^\text{II}_o(e)$. Its full significance will become evident once we arrive at Theorem \ref{dualfund}.

\medskip

Let $e'$ be the flip of $e$, $C'$ of descent $e'$ a cycle of a skew basic $(p', q')$-leaper $L'$, $p' < q'$, over a $(p' + q') \times (p' + q')$ chessboard, and $D'$ of origin $o$ and descent $e'$ the canonical second leaper cycle of the canonical second leaper $M'$ of proportions $r'$ and $s'$, $r' < s'$, associated with $C'$.

We go on to apply Theorem \ref{flip} to the duality identity for $\Phi(e)$ and $\Phi^\text{II}_o(e)$ in order to obtain the isomorphisms from $C'$ onto $D'$ from the isomorphisms from $D$ onto $C$.

By Theorem \ref{flip} and its proof, there exists an enumeration \[b^\text{I}_1, b^\text{I}_2, \ldots, b^\text{I}_{4(m + n)}\] (with $b^\text{I}_{j + 4(m + n)} \equiv b^\text{I}_j$ for all integers $j$) of the vertices of $\Phi(e')$ such that, for all integers $j$, \[A^\text{II}(j) = P_oA^\text{I}_\text{Flip}(j)Q_o,\] where $A^\text{I}_\text{Flip}(j)$ is the direction matrix associated with the arc from $b^\text{I}_j$ to $b^\text{I}_{j + 1}$ in $\Phi(e')$.

Analogously, by Theorem \ref{flip} and its proof, there exists an enumeration \[b^\text{II}_1, b^\text{II}_2, \ldots, b^\text{II}_{4(m + n)}\] (with $b^\text{II}_{j + 4(m + n)} \equiv b^\text{II}_j$ for all integers $j$) of the vertices of $\Phi^\text{II}_o(e')$ such that, for all integers $j$, \[A^\text{II}_\text{Flip}(j) = P_oA^\text{I}(j)Q_o,\] where $A^\text{II}_\text{Flip}(j)$ is the direction matrix associated with the arc from $b^\text{II}_j$ to $b^\text{II}_{j + 1}$ in $\Phi^\text{II}_o(e')$.

Let, for all integers $j'$ and $j''$, $\text{Dist}^\text{I}_\text{Flip}(j', j'')$ be the distance from $b^\text{I}_{j'}$ to $b^\text{I}_{j''}$ in $\Phi(e')$. Then we can rewrite the duality identity for $\Phi(e)$ and $\Phi^\text{II}_o(e)$ as \[P_o\text{Dist}^\text{I}_\text{Flip}(\psi(j), \psi(j + 1))Q_o = P^{-1}_oA^\text{II}_\text{Flip}(j)Q^{-1}_oA_eA_o\] for all integers $j$.

Since $P_o = P^{-1}_o$ in both cases $o = \texttt{g}$ and $o = \texttt{h}$, this is equivalent to \[\text{Dist}^\text{I}_\text{Flip}(\psi(j), \psi(j + 1)) = A^\text{II}_\text{Flip}(j)Q^{-1}_oA_eA_oQ^{-1}_o.\]

Let $A^\text{Unit}$ be the unit $2 \times 2$ matrix, $e = e_1e_2\ldots e_l$, and $k$ the number of occurrences of the character \texttt{h} in $e$.

In both cases $o = \texttt{g}$ and $o = \texttt{h}$, \[A_\texttt{f}A_oQ^{-1}_oA_\texttt{g} = A_\texttt{g}A_oQ^{-1}_oA_\texttt{f} = A_oQ^{-1}_o,\] \[A_\texttt{h}A_oQ^{-1}_oA_\texttt{h} = -A_oQ^{-1}_o,\] and \[Q^{-1}_oA_oQ^{-1}_oA_o = -A^\text{Unit}.\]

Thus \begin{align*} Q^{-1}_oA_eA_oQ^{-1}_oA_{e'}A_o &= Q^{-1}_oA_{e_1}A_{e_2}\cdots A_{e_l}A_oQ^{-1}_oA_{\overline{e_l}}A_{\overline{e_{l - 1}}}\cdots A_{\overline{e_1}}A_o\\ &= \pm Q^{-1}_oA_{e_1}A_{e_2}\cdots A_{e_{l - 1}}A_oQ^{-1}_oA_{\overline{e_{l - 1}}}A_{\overline{e_{l - 2}}}\cdots A_{\overline{e_1}}A_o\\ &= \cdots = (-1)^kQ^{-1}_oA_oQ^{-1}_oA_o\\ &= (-1)^{k + 1} A^\text{Unit}.\end{align*}

Therefore, \[\text{Dist}^\text{I}_\text{Flip}(\psi(j), \psi(j + 1))A_{e'}A_o = (-1)^{k + 1} A^\text{II}_\text{Flip}(j).\]

Since the proportions of $L'$ and $M'$ satisfy \[\left(\begin{array}{c} p'\\ q' \end{array}\right) = A_{e'}A_{o}\left(\begin{array}{c} r'\\ s' \end{array}\right),\] it follows that \[\text{Dist}^\text{I}_\text{Flip}(\psi(j), \psi(j + 1))\left(\begin{array}{c} p'\\ q' \end{array}\right) = (-1)^{k + 1} A^\text{II}_\text{Flip}(j)\left(\begin{array}{c} r'\\ s' \end{array}\right).\]

Let $C' = b_1b_2\ldots b_{4(m + n)}$ (with $b_{j + 4(m + n)} \equiv b_j$ for all integers $j$) so that $C'$ is an instantiation of $\Phi(e')$ under the mapping $b^\text{I}_j \to b_j$ if $k$ is odd, and $b^\text{I}_j \to b_{j + 2(m + n)}$ if it is even. Then $D'$ is an instantiation of $\Phi^\text{II}_o(e')$ under the mapping $b^\text{II}_j \to b_{\psi(j)}$.

Therefore, the mapping \[b_j \to b_{\psi(j)}\] is an isomorphism from $C'$ onto $D'$. All other such isomorphisms are obtained from it by means of rotation and reflection. Hence the following theorem.

\medskip

\begin{theorem} Let the strings $e'$ and $e''$ composed of the characters \texttt{f}, \texttt{g}, and \texttt{h} be flips of each other. Let $e'$ be the common descent of the cycle $C'$ of a skew basic $(p', q')$-leaper within a $(p' + q') \times (p' + q')$ chessboard and the canonical second leaper cycle $D'$ over the squares of $C'$, $e''$ the common descent of the cycle $C''$ of a skew basic $(p'', q'')$-leaper within a $(p'' + q'') \times (p'' + q'')$ chessboard and the canonical second leaper cycle $D''$ over the squares of $C''$, and $4\mu$ the common length of $C'$, $D'$, $C''$, and $D''$.

Let $\psi$ be a permutation of 1, 2, \ldots, $4\mu$. Then there exists an enumeration $C' = a'_1a'_2\ldots a'_{4\mu}$ of the vertices of $C'$ such that the mapping $a'_j \to a'_{\psi(j)}$ is an isomorphism from $C'$ onto $D'$ if and only if there exists an enumeration $D'' = a''_1a''_2\ldots a''_{4\mu}$ of the vertices of $D''$ such that the mapping $a''_j \to a''_{\psi(j)}$ is an isomorphism from $D''$ onto $C''$. \label{isoflip} \end{theorem}

\medskip

Theorem \ref{isoflip} allows us to extend Theorem \ref{seconddispl} as follows.

\medskip

\begin{theorem} In the setting of Theorem \ref{seconddispl}, given a square $c$ of $C$, let the arc from $c$ in $D$ point to $d$ and the \emph{$C$-displacement} of $c$ be the number of steps along $C$ that lead from $c$ to $d$.

Then the multiset of all $C$-displacements coincides with the multiset of all $D$-displacements. \label{displ} \end{theorem}

\medskip

\begin{proof} By Theorems \ref{seconddispl} and \ref{isoflip}, the multiset of all $C$-displacements consists of all integers between $\mu$ and $3\mu$ congruent to $\beta$ modulo four, each occurring precisely eight times, where $\beta = 1$ or $\beta = 3$ depending on the orientations of $C$ and $D$.

Let the arc from a square $a$ in $C$ point to a square $b$. Then the number of steps along $D$ that lead from $a$ to $b$ is congruent to $\alpha$ modulo four.

Replace each of those steps with a path along $C$ advancing from the same starting vertex to the same ending vertex as that step. Since each step along $D$ is replaced with a number of steps along $C$ congruent to $\beta$ modulo four, the total number of steps along $C$ in the resulting path, leading from $a$ to $b$, is congruent to $\alpha\beta$ modulo four.

Thus $\alpha\beta$ is congruent to one modulo four and $\alpha = \beta$, as needed. \end{proof}

\medskip

The duality identity for $\Phi(e)$ and $\Phi^\text{II}_o(e)$ shows that the underlying source of the duality of a dual board given by Theorem \ref{dual} is the duality of the associated direction graphs.

We proceed to outline the key conditions that a direction graph needs to satisfy so that its instantiations give rise to a family of dual boards in a similar manner.

\medskip

\begin{definition} A direction graph $\Phi$ is \emph{dual} if there exist a direction graph $\Phi'$, a \emph{complement} of $\Phi$, a one-to-one mapping $\eta$ between the vertex sets of $\Phi$ and $\Phi'$, and a real-coefficient $2 \times 2$ \emph{duality matrix} $A_\text{Dual}$, distinct from all direction matrices, such that the following conditions hold.

Both of $\Phi$ and $\Phi'$ are coherent and connected.

The unlabeled symmetric directed graphs obtained from $\Phi$ and $\Phi'$ by discarding all arc labels are isomorphic.

For every arc from $\eta(a)$ to $\eta(b)$ labeled $i$ in $\Phi'$, the distance from $a$ to $b$ in $\Phi$ equals $A_iA_\text{Dual}$. \end{definition}

\medskip

The basic properties of dual direction graphs are summarized in the following theorem.

\medskip

\begin{theorem} Let $\Phi$ be a dual direction graph of complement $\Phi'$ and duality matrix $A_\text{Dual}$. Then:

(a) $\pi(\Phi)$ is a dual direction graph for all equivalence permutations $\pi$.

(b) $A_\text{Dual}$ is a unit-determinant integer-coefficient matrix.

(c) $\Phi'$ is a dual direction graph of complement $\Phi$ and duality matrix $A^{-1}_\text{Dual}$.

(d) The vertex sets of all but finitely many basic leaper instantiations of $\Phi$ are dual boards. \label{dualbasic} \end{theorem}

\medskip

\begin{proof} Let $\eta$ be a one-to-one mapping between the vertex sets of $\Phi$ and $\Phi'$ as in the definition of a dual direction graph.

Let the matrices $P$ and $Q$ induce an equivalence permutation. Then $Q$ is invertible and $P\Phi Q$ is a dual direction graph of complement $P\Phi' Q$ and duality matrix $Q^{-1}A_\text{Dual}Q$. This settles (a).

Let $L$ be a basic $(p, q)$-leaper, $p < q$, and $G$ an $L$-instantiation of $\Phi$ under the mapping $\tau$. Since $\Phi$ is coherent, $G$ is the leaper graph of $L$ over the vertex set of $G$ for all but finitely many basic leapers $L$. Suppose that $L$ is indeed such a basic leaper.

Let $i$ be the label of the arc from $\eta(a)$ to $\eta(b)$ in $\Phi'$ and $A$ the distance from $a$ to $b$ in $\Phi$.

Since $A_\text{Dual} = A^{-1}_iA$, it is an integer-coefficient matrix.

Let \[\left(\begin{array}{c} r'\\ s' \end{array}\right) = A_\text{Dual}\left(\begin{array}{c} p\\ q \end{array}\right),\] $r = \min\{|r'|, |s'|\}$, and $s = \max\{|r'|, |s'|\}.$ Then $r$ and $s$ are nonnegative integers and $r \le s$. Let $M$ be an $(r, s)$-leaper.

The translation $v$ from $\tau(a)$ to $\tau(b)$ satisfies \[v^T = A_iA_\text{Dual} \left(\begin{array}{c} p\\ q \end{array}\right) = A_i \left(\begin{array}{c} r'\\ s' \end{array}\right).\] Therefore, $v$ is an $M$-translation and $\tau(a)$ and $\tau(b)$ are joined by a move of $M$.

Since $\Phi'$ is connected, the leaper graph of $M$ over the vertex set of $G$ is connected. Therefore, $M$ can imitate a move of $L$.

Since $L$ is a basic leaper, the leaper graph of $L$ over the infinite chessboard is connected. Hence the leaper graph of $M$ over the infinite chessboard is also connected and $M$ is also a basic leaper.

It follows that the coordinates of the column vector $A_\text{Dual}\left( \begin{smallmatrix} p\\ q \end{smallmatrix} \right)$ are relatively prime for all relatively prime $p$ and $q$ such that $p + q$ is odd and $p < q$. Therefore, $A_\text{Dual}$ is invertible.

Let $\text{Dist}(x, y)$ be the distance from $x$ to $y$ in $\Phi$ and $\text{Dist}'(x', y')$ the distance from $x'$ to $y'$ in $\Phi'$. It follows from the definition of a dual direction graph that, for all vertices $x$ and $y$ of $\Phi$, \[\text{Dist}(x, y) = \text{Dist}'(\eta(x), \eta(y))A_\text{Dual}.\]

Therefore, \[\text{Dist}'(x', y') = \text{Dist}(\eta^{-1}(x'), \eta^{-1}(y'))A^{-1}_\text{Dual}\] for all vertices $x'$ and $y'$ of $\Phi'$, and in particular for all $x'$ and $y'$ such that $\eta^{-1}(x')$ and $\eta^{-1}(y')$ are adjacent in $\Phi$.

Thus $\Phi'$ is a dual direction graph of complement $\Phi$ and duality matrix $A^{-1}_\text{Dual}$, settling (c).

Consequently, $A^{-1}_\text{Dual}$ is also a duality matrix and thus integer-coefficient. It follows that $A_\text{Dual}$ is unit-determinant, settling (b).

Since $A_\text{Dual}$ is distinct from all direction matrices, $L$ and $M$ are distinct leapers for all but finitely many basic leapers $L$. Suppose that $L$ is indeed such a basic leaper.

There exists a direction matrix $A_j$, unique in the case when $M$ is a skew leaper, such that \[\left(\begin{array}{c} r\\ s \end{array}\right) = A_j\left(\begin{array}{c} r'\\ s' \end{array}\right).\]

Consider the direction graph $\Phi'A^{-1}_j$. The graph $G'$ over the vertex set of $G$ such that an arc points from $\tau(x)$ to $\tau(y)$ in $G'$ if and only if an arc points from $\eta(x)$ to $\eta(y)$ in $\Phi'$ is an $M$-instantiation of $\Phi'A^{-1}_j$ under the mapping $\eta(x) \to \tau(x)$.

Since $\Phi'$ is coherent, so is $\Phi'A^{-1}_j$ for every direction matrix $A_j$. Therefore, $G'$ is the leaper graph of $M$ over the squares of $G'$ for all but finitely many basic leapers $M$ and, since $A_\text{Dual}$ is invertible, for all but finitely many basic leapers $L$. Suppose that $L$ is indeed such a basic leaper.

Then the vertex set of the $L$-instantiation of $\Phi$ is a board dual with respect to $L$ and $M$, settling (d) and completing the proof of the theorem. \end{proof}

\medskip

\begin{theorem} Let $e$ be a string composed of the characters \texttt{f}, \texttt{g}, and \texttt{h}, and $o$ one of the characters \texttt{g} and \texttt{h}. Then both the fundamental direction cycle $\Phi(e)$ and the second fundamental direction cycle $\Phi^\text{II}_o(e)$ are dual direction graphs. Each of them is a complement of the other, the duality matrix of $\Phi^\text{II}_o(e)$ being $A_eA_o$ and that of $\Phi(e)$ its inverse. \label{dualfund} \end{theorem}

\medskip

\begin{proof} By the duality identity for $\Phi(e)$ and $\Phi^\text{II}_o(e)$ and Theorem \ref{dualbasic}. \end{proof}

\medskip

In general, dual boards and dual direction graphs are not related in as straightforward a manner as Theorems \ref{dual} and \ref{dualfund} might suggest. For instance, there exist boards dual with respect to two distinct skew leapers such that the extracted direction graphs are not complements of each other. One example of this is the board formed by the squares of the cycles depicted in Figures \ref{2355}, top and \ref{2355e}, viewed as dual with respect to a $(1, 2)$-leaper and a $(2, 3)$-leaper.

\medskip

\begin{question} Do there exist boards dual with respect to two distinct skew leapers such that at least one of the extracted direction graphs is not a dual direction graph, or such that at least one of the associated leaper graphs contains a nontrivial cycle and does not allow the extraction of a direction graph at all? \label{counterexample} \end{question}

\medskip

Similarly to the notion of a dual board, the notion of a dual direction graph raises a number of questions.

\medskip

\begin{question} Given a unit-determinant integer-coefficient $2 \times 2$ matrix $A$ distinct from all direction matrices, does there exist a dual direction graph of duality matrix $A$? \label{matrixgraph} \end{question}

\medskip

By Theorem \ref{fghhgraph}, a sufficient condition for $A$ is that it is a nonempty product of the matrices $A_\texttt{f}$, $A_\texttt{g}$, and $A_\texttt{h}$.

\medskip

\begin{question} Given a matrix $A$ such that a dual direction graph of duality matrix $A$ does exist, what is the least number of vertices that it may contain? Is the number of vertices that it may contain unbounded from above? \label{graphbound} \end{question}

\medskip

Analogously to the case of a dual board, a dual direction graph given by Theorem \ref{dualfund} is never minimal, as removing any vertex from it yields a dual direction graph of the same duality matrix.

Corollary \ref{pinwheelunbound} answers one special case of the second part of the question.

\medskip

\begin{question} For what positive integers $n \ge 2$ do there exist $n$ dual direction graphs which are pairwise complements? \label{multgraph} \end{question}

\medskip

In particular (as, unlike in the case of a dual board, Theorem \ref{tl} does not yield an example of $n = 3$), do there exist three dual direction graphs which are pairwise complements?

\section[Constructions of Dual Boards and Dual Direction Graphs I]{Constructions of Dual Boards\\ and Dual Direction Graphs I}

We set out to add one more lifting transformation, $\hbar$, to $f$, $g$, and $h$. To this end, first we extend the definitions of $f$, $g$, and $h$ so that they act upon a more general family of cycles.

Introduce a Cartesian coordinate system $Oxy$ over the infinite chessboard such that the centers of the squares are the integer points. For each direction $i$ out of \texttt{E}, \texttt{NE}, \ldots, \texttt{SE}, let $\alpha^i$ be the ray emanating from $O$ and pointing in direction $i$.

We say that a square lies inside, or strictly inside, a region if its center does.

Given a path $w$ through a number of squares, we write $w_\text{Start}$ and $w_\text{End}$ for the opening and final squares of $w$, and, if $w$ is of odd length, $w_\text{Mid}$ for the middle square of $w$.

In the context of a cycle $C$ and an orientation imposed on $C$, given two squares $a$ and $b$ of $C$, we write $[a; b]$ for the portion of $C$ that runs from $a$ to $b$ and $(a; b)$ for the portion of $C$ that runs from the square directly following $a$ to the square directly preceding $b$.

Given a path or a cycle $w$ through a number of squares, we write $\kappa(w)$ for the broken line formed by joining the centers of the squares of $w$ in the order in which they are visited by $w$.

\medskip

\begin{definition} A cycle $D$ of a leaper $M$ is \emph{$(p, q)$-perfect}, $p < q$, if it admits a \emph{perfect partitioning} into eight disjoint nonempty paths, \[D = a^\texttt{E}a^\texttt{NE}\ldots a^\texttt{SE},\] satisfying the following conditions.

\emph{Translation.} \begin{align*} a^\texttt{E}a^\texttt{NE} + (-q, -p) &= -a^\texttt{SW}{-a^\texttt{W}}\\ a^\texttt{NE}a^\texttt{N} + (-p, -q) &= -a^\texttt{S}{-a^\texttt{SW}}\\ a^\texttt{N}a^\texttt{NW} + (p, -q) &= -a^\texttt{SE}{-a^\texttt{S}}\\ a^\texttt{NW}a^\texttt{W} + (q, -p) &= -a^\texttt{E}{-a^\texttt{SE}} \end{align*}

\emph{Symmetry.} For all directions $i$ out of \texttt{E}, \texttt{NE}, \ldots, \texttt{SE}, $a^i$ is of odd length, $a^i_\text{Mid}$ lies on the ray $\alpha^i$, and $a^i$ is symmetric with respect to $\alpha^i$. Furthermore, a $90^\circ$ rotation about $O$ maps $a^i$ onto $a^{i + 2}$.

\emph{Separation.} For all directions $i$ out of \texttt{E}, \texttt{NE}, \ldots, \texttt{SE}, the counterclockwise portion $(a^i_\text{Mid}; a^{i + 1}_\text{Mid})$ of $D$ lies strictly inside the $45^\circ$ angle between the rays $\alpha^i$ and $\alpha^{i + 1}$.

\emph{Simplicity.} The broken line $\kappa(D)$ is the contour of a simple polygon. In particular, $D$ is a simple cycle.

\emph{Coherence.} $D$ is the leaper graph of $M$ over the squares of $D$. \end{definition}

\emph{Protocoherence.} No square of $D$ is joined by a move of $M$ to a square of either $D + (2q, 0)$ or $D + (p + q, p + q)$.

\medskip

\begin{figure}[ht!]\hspace*{\fill}\includegraphics{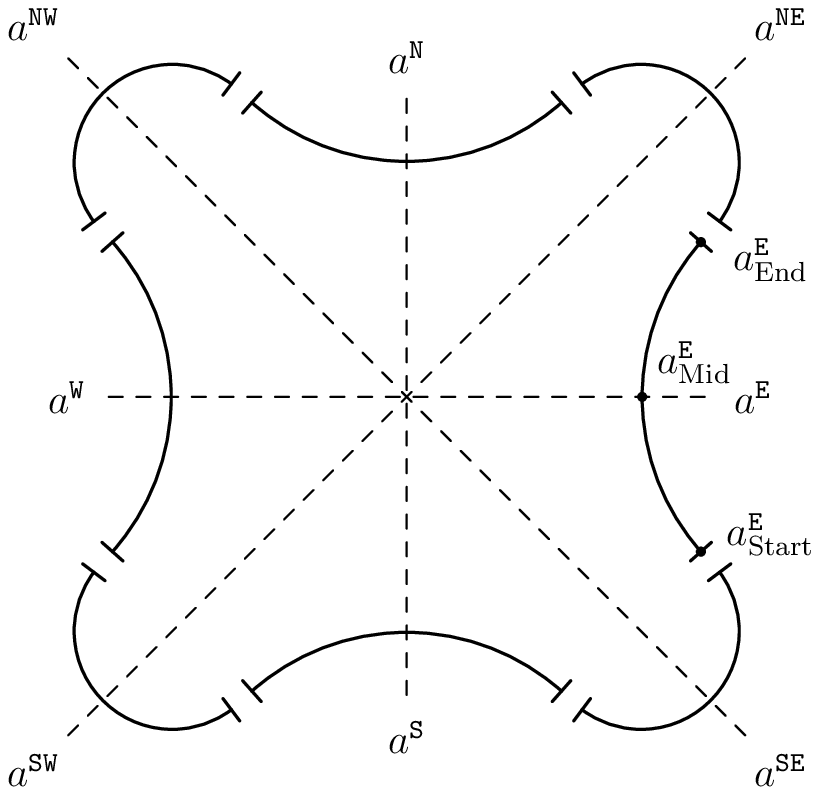}\hspace*{\fill}\caption{}\label{cycle}\end{figure}

Figure \ref{cycle} illustrates the definition of a perfect cycle.

We begin with a couple of lemmas.

\medskip

\begin{lemma} Let $D$ be a $(p, q)$-perfect cycle. Then $D$ lies strictly inside the square bounded by the lines $y = \pm x \pm (p + q)$, and strictly inside the square bounded by the lines $x = \pm q$ and $y = \pm q$. \label{bound} \end{lemma}

\medskip

\begin{proof} Let $D = a^\texttt{E}a^\texttt{NE}\ldots a^\texttt{SE}$ be a perfect partitioning of $D$ and $L$ a $(p, q)$-leaper.

By the translation property of $D$, every square in $a^\texttt{E}a^\texttt{NE}$ is joined by an $L$-move to a square of $a^\texttt{W}a^\texttt{SW}$. By the separation property of $D$, $a^\texttt{W}a^\texttt{SW}$ lies strictly below and to the left of the line $y = -x$. Therefore, $a^\texttt{E}a^\texttt{NE}$ lies strictly below and to the left of the line $y = -x + p + q$. Analogously, so does $a^\texttt{NE}a^\texttt{N}$.

Since the counterclockwise portion $[a^\texttt{E}_\text{Mid}; a^\texttt{N}_\text{Mid}]$ of $D$ is a subpath of $a^\texttt{E}a^\texttt{NE}a^\texttt{N}$, it also lies strictly below and to the left of the line $y = -x + p + q$. By this and the separation property of $D$, $[a^\texttt{E}_\text{Mid}; a^\texttt{N}_\text{Mid}]$ lies inside the triangle $\Delta$ formed by the lines $x = 0$, $y = 0$, and $y = -x + p + q$, with only $a^\texttt{E}_\text{Mid}$ and $a^\texttt{N}_\text{Mid}$ on the contour of $\Delta$, in the interiors of the two legs.

Analogous reasoning applies to $[a^i_\text{Mid}; a^{i + 2}_\text{Mid}]$ for $i = \texttt{N}$, \texttt{W}, and \texttt{S}. Therefore, $D$ lies strictly inside the square bounded by the lines $y = \pm x \pm (p + q)$.

Analogously, $D$ lies strictly inside the square bounded by the lines $x = \pm q$ and $y = \pm q$. \end{proof}

\medskip

\begin{lemma} Let $L$ be a $(p, q)$-leaper, $p < q$, and $D$ a $(p, q)$-perfect cycle. Then two squares of $D$ are joined by an $L$-move if and only if they are joined by an $L$-move by the translation property of $D$. \label{longcoh} \end{lemma}

\medskip

\begin{proof} Let $D = a^\texttt{E}a^\texttt{NE}\ldots a^\texttt{SE}$ be a perfect partitioning of $D$.

Consider any $L$-translation $v$, say, $(p, q)$. It suffices to show that $D$ and $D + v$ do not have any squares in common other than the squares of $a^\texttt{NE}a^\texttt{N}$.

By the separation property of $D$ and Lemma \ref{bound}, the line $y = q$ strictly separates $a^\texttt{NE}a^\texttt{N} + v$ and $D$.

Consider $a^\texttt{NW}a^\texttt{W} + v$.

Let $R$ be the point $(\frac{p - q}{2}, \frac{p + q}{2})$. Consider the cyclic list of broken lines \[ \kappa(a^\texttt{NW}a^\texttt{W}),\; \kappa(a^\texttt{NE}a^\texttt{N}),\; \kappa(a^\texttt{NW}a^\texttt{W} + v),\; \kappa(a^\texttt{NE}a^\texttt{N} + (-q, p)). \]

By the translation and symmetry properties of $D$, $90^\circ$ rotation about $R$ maps each item in this list onto the following one, cyclically. By the simplicity property of $D$, no two items in direct succession have a common point. It follows that no two items have a common point at all. Consequently, $a^\texttt{NW}a^\texttt{W} + v$ has no squares in common with $a^\texttt{NE}a^\texttt{N}$ and $a^\texttt{NW}a^\texttt{W}$.

By the separation property of $D$ and Lemma \ref{bound}, the line $y = 0$ strictly separates $a^\texttt{NW}a^\texttt{W} + v$ and $a^\texttt{SW}a^\texttt{S}$.

Lastly, by the separation property of $D$, the line $y = x$ strictly separates $a^\texttt{NW}a^\texttt{W} + v$ and $a^\texttt{SE}a^\texttt{E}$.

It follows that $a^\texttt{NW}a^\texttt{W} + v$ does not have any squares in common with $D$.

Analogously, neither does $a^\texttt{SE}a^\texttt{E} + v$. \end{proof}

\medskip

By Lemma \ref{longcoh}, a $(p, q)$-perfect cycle admits a unique perfect partitioning. Indeed, $a^\texttt{E}$ is uniquely determined as the intersection of $D$, $D + (q, p)$, and $D + (q, -p)$, and analogously $a^\texttt{NE}$, $a^\texttt{N}$, \ldots, $a^\texttt{SE}$ are uniquely determined as well.

Here, then, follow the extended definitions of $f$, $g$, and $h$.

\medskip

\begin{definition} Let $D$ be a $(p, q)$-perfect cycle of perfect partitioning $D = a^\texttt{E}a^\texttt{NE}\allowbreak\ldots a^\texttt{SE}$. Let \[b^\texttt{E} = -a^\texttt{SW}{-a^\texttt{W}}{-a^\texttt{NW}} + (p + q, 0)\] and \[b^\texttt{NE} = a^\texttt{NE} + (p, p),\] and define $b^\texttt{N}$, $b^\texttt{NW}$, \ldots, $b^\texttt{SE}$ symmetrically. Then the \emph{$f$-lift} of $D$ is the cycle $f(D) = b^\texttt{E}b^\texttt{NE}\allowbreak\ldots b^\texttt{SE}$ together with its partitioning into $b^\texttt{E}$, $b^\texttt{NE}$, \ldots, $b^\texttt{SE}$. \end{definition}

\medskip

\begin{definition} Let $D$ be a $(p, q)$-perfect cycle of perfect partitioning $D = a^\texttt{E}a^\texttt{NE}\allowbreak\ldots a^\texttt{SE}$. Let \[b^\texttt{E} = a^\texttt{E} + (q - p, 0)\] and \[b^\texttt{NE} = -a^\texttt{S}{-a^\texttt{SW}}{-a^\texttt{W}} + (q, q),\] and define $b^\texttt{N}$, $b^\texttt{NW}$, \ldots, $b^\texttt{SE}$ symmetrically. Then the \emph{$g$-lift} of $D$ is the cycle $g(D) = b^\texttt{E}b^\texttt{NE}\allowbreak\ldots b^\texttt{SE}$ together with its partitioning into $b^\texttt{E}$, $b^\texttt{NE}$, \ldots, $b^\texttt{SE}$. \end{definition}

\medskip

\begin{definition} Let $D$ be a $(p, q)$-perfect cycle of perfect partitioning $D = a^\texttt{E}a^\texttt{NE}\allowbreak\ldots a^\texttt{SE}$. Let \[b^\texttt{E} = -a^\texttt{S}{-a^\texttt{SW}}{-a^\texttt{W}}{-a^\texttt{NW}}{-a^\texttt{N}} + (p + q, 0)\] and \[b^\texttt{NE} = a^\texttt{E}a^\texttt{NE}a^\texttt{N} + (q, q),\] and define $b^\texttt{N}$, $b^\texttt{NW}$, \ldots, $b^\texttt{SE}$ symmetrically. Then the \emph{$h$-lift} of $D$ is the cycle $h(D) = b^\texttt{E}b^\texttt{NE}\allowbreak\ldots b^\texttt{SE}$ together with its partitioning into $b^\texttt{E}$, $b^\texttt{NE}$, \ldots, $b^\texttt{SE}$. \end{definition}

\medskip

We are ready to introduce $\hbar$.

\medskip

\begin{definition} Let $D$ be a $(p, q)$-perfect cycle of perfect partitioning $D = a^\texttt{E}a^\texttt{NE}\allowbreak\ldots a^\texttt{SE}$. Let \[b^\texttt{E} = a^\texttt{SE}a^\texttt{E}a^\texttt{NE} + (p + q, 0)\] and \[b^\texttt{NE} = -a^\texttt{SE}{-a^\texttt{S}}{-a^\texttt{SW}}{-a^\texttt{W}}{-a^\texttt{NW}} + (q, q),\] and define $b^\texttt{N}$, $b^\texttt{NW}$, \ldots, $b^\texttt{SE}$ symmetrically. Then the \emph{$\hbar$-lift} of $D$ is the cycle $\hbar(D) = b^\texttt{E}b^\texttt{NE}\allowbreak\ldots b^\texttt{SE}$ together with its partitioning into $b^\texttt{E}$, $b^\texttt{NE}$, \ldots, $b^\texttt{SE}$. \end{definition}

\medskip

The lifting transformations $h$ and $\hbar$ are companions in a way analogous to how $f$ and $g$ are. In a sense, the two lifting transformations in each companion pair are $45^\circ$ rotations of each other.

The following lemma establishes the key properties of the extended $f$, $g$, and $h$ and $\hbar$ as they relate to perfect cycles.

\medskip

\begin{lemma} Let $D$ be a $(p, q)$-perfect cycle of an $(r, s)$-leaper $M$. Then:

(a) Each of the four lifts of $D$ is a cycle of $M$.

(b) $f(D)$ possesses the translation property with parameters $(p, 2p + q)$, $g(D)$ possesses the translation property with parameters $(q, 2q - p)$, and both of $h(D)$ and $\hbar(D)$ possess the translation property with parameters $(q, p + 2q)$.

(c) Each of the four lifts of $D$ possesses the symmetry property.

(d) Each of the four lifts of $D$ possesses the separation property.

(e) Each of the four lifts of $D$ possesses the simplicity property.

(f) Provided that $s \le p$ and $r + s \le q - p$, each of the four lifts of $D$ possesses the coherence property.

(g) Provided that $s \le p$ and $r + s \le q - p$, each of the four lifts of $D$ possesses the protocoherence property with the same parameters as in (b). \label{fghh} \end{lemma}

\medskip

\begin{proof} The proofs of parts (a) and (b) of the lemma in all cases $f$, $g$, $h$, and $\hbar$ are analogous to the proofs of Lemmas \ref{f2}, \ref{g2}, and \ref{h2}.

Part (c) follows by the symmetry property of $D$.

The main motif of the proofs of parts (d), (e), (f), and (g) is that, locally, a lift of $D$ looks precisely like $D$. In other words, for small regions $R$, the subgraph of a lift of $D$ within $R$ is a translation copy of a subgraph of $D$ under one of the associated lifting translations.

Let $D = a^\texttt{E}a^\texttt{NE}\allowbreak\ldots a^\texttt{SE}$ be the perfect partitioning of $D$.

\medskip

\emph{Case $f$.} Let $b^\texttt{E}$, $b^\texttt{NE}$, \ldots, $b^\texttt{SE}$ be as in the definition of $f$.

Consider the counterclockwise portion $(b^\texttt{E}_\text{Mid}; b^\texttt{NE}_\text{Mid})$ of $f(D)$.

By the separation property of $D$ and Lemma \ref{bound}, $-[a^\texttt{NW}_\text{Start}; a^\texttt{W}_\text{Mid}) + (p + q, 0)$ lies in the interior of the $45^\circ$ angle between $\alpha^\texttt{E}$ and $\alpha^\texttt{NE}$. By the separation property of $D$, so does $[a^\texttt{NE}_\text{Start}; a^\texttt{NE}_\text{Mid}) + (p, p)$. Since it is the concatenation of those two paths, so does $(b^\texttt{E}_\text{Mid}; b^\texttt{NE}_\text{Mid})$ as well. This and the symmetry property of $f(D)$ settle (d).

For (e), by the symmetry and separation properties of $f(D)$ it suffices to establish that the broken line $\kappa([b^\texttt{E}_\text{Mid}; b^\texttt{NE}_\text{Mid}])$ does not intersect itself.

By the translation property of $D$, $[b^\texttt{E}_\text{Mid}; b^\texttt{NE}_\text{Mid}]$ is a subpath of $a^\texttt{SE}a^\texttt{E}a^\texttt{NE} + (p, p)$. By the simplicity property of $D$, the broken line tracing the latter does not intersect itself, and this settles (e).

For (f), say that two squares of $f(D)$ form an \emph{illicit pair} if they are joined by a move of $M$, but this move is not an edge of $f(D)$.

Let $i$ be any direction out of \texttt{E}, \texttt{NE}, \ldots, \texttt{SE}.

By the coherence property of $D$, no illicit pair can occur within $b^i$.

Suppose that an illicit pair has one square $a$ in $b^\texttt{E}b^\texttt{NE}$ and another square $b$ in $b^\texttt{N}b^\texttt{NW}$.

By the coherence property of $D$, no illicit pair can occur within $\Pi = a^\texttt{SE}a^\texttt{E}a^\texttt{NE}a^\texttt{N}a^\texttt{NW} + (p, p)$. Thus either $a$ belongs to $\Pi_1 = b^\texttt{E}b^\texttt{NE} \setminus \Pi$ or $b$ belongs to $\Pi_2 = b^\texttt{N}b^\texttt{NW} \setminus \Pi$.

By the translation property of $D$, $\Pi_1 = -a^\texttt{SW} + (p + q, 0)$. Thus, by the separation property of $D$ and Lemma \ref{bound}, the strip bounded by the lines $y = 0$ and $y = p$ strictly separates $\Pi_1$ and $b^\texttt{N}b^\texttt{NW}$. However, since $s \le p$, no illicit pair can have its two squares on different sides of this strip.

By the translation property of $D$, $\Pi_2 = (-a^\texttt{SW} + (0, p + q))b^\texttt{NW}$. Thus, by the separation property of $D$ and Lemma \ref{bound}, the strip bounded by the lines $x = 0$ and $x = p$ strictly separates $b^\texttt{E}b^\texttt{NE}$ and $\Pi_2$. However, since $s \le p$, no illicit pair can have its two squares on different sides of this strip.

We have arrived at a contradiction. By symmetry, it follows that no illicit pair can occur within $b^i \cup b^{i + 1}$, $b^i \cup b^{i + 2}$, and $b^i \cup b^{i + 3}$.

Consider, then, $b^\texttt{E}b^\texttt{NE}$ and $b^\texttt{W}b^\texttt{SW}$. By the separation property of $D$ and Lemma \ref{bound}, the strip bounded by the lines $x = -p$ and $x = p$ strictly separates the former and the latter. As above, by symmetry it follows that no illicit pair can occur within $b^i \cup b^{i + 4}$.

This completes the proof of (f).

For (g), by Lemma \ref{bound} the strip bounded by the lines $x = p + q$ and $x = 3p + q$ strictly separates $f(D)$ and $f(D) + (4p + 2q, 0)$. However, since $s \le p < 2p$, no two squares on different sides of this strip can be joined by a move of $M$.

Suppose, then, that a square $a$ in $f(D)$ is joined by a move of $M$ to a square $b$ in $f(D) + (3p + q, 3p + q)$.

By Lemma \ref{bound}, $f(D) + (3p + q, 3p + q)$ lies strictly to the right of the line $x = 2p$ and strictly above the line $y = 2p$. As above, it follows from this that $a$ lies in the portion of $f(D)$ strictly to the right of the line $x = 0$ and strictly above the line $y = 0$. By the separation property of $f(D)$ and the translation property of $D$, this portion of $f(D)$ is a subpath of $D + (p, p)$.

Analogously, $b$ lies in $D + (2p + q, 2p + q)$. However, by the protocoherence property of $D$, no square of $D + (p, p)$ is joined by a move of $M$ to a square of $D + (2p + q, 2p + q)$. We have arrived at a contradiction and the proof of (g) is complete.

\medskip

\emph{Case $g$.} This case is analogous to Case $f$.

\medskip

\emph{Case $h$.} Let $b^\texttt{E}$, $b^\texttt{NE}$, \ldots, $b^\texttt{SE}$ be as in the definition of $h$.

Consider the counterclockwise portion $(b^\texttt{E}_\text{Mid}; b^\texttt{NE}_\text{Mid})$ of $h(D)$.

By the separation property of $D$ and Lemma \ref{bound}, $-[a^\texttt{N}_\text{Start}; a^\texttt{W}_\text{Mid}) + (p + q, 0)$ lies in the interior of the $45^\circ$ angle between $\alpha^\texttt{E}$ and $\alpha^\texttt{NE}$. By the separation property of $D$ and Lemma \ref{bound}, so does $[a^\texttt{E}_\text{Start}; a^\texttt{NE}_\text{Mid}) + (q, q)$. Since it is the concatenation of those two paths, so does $(b^\texttt{E}_\text{Mid}; b^\texttt{NE}_\text{Mid})$ as well. This and the symmetry property of $h(D)$ settle (d).

As in Case $f$, for (e) it suffices to establish that the broken line $\kappa([b^\texttt{E}_\text{Mid}; b^\texttt{NE}_\text{Mid}])$ does not intersect itself.

Let $o = a^\texttt{NW}_\text{End} + (p + q, 0)$. By the translation property of $D$, $[o; b^\texttt{NE}_\text{Mid}]$ is a subpath of $a^\texttt{S}a^\texttt{SE}a^\texttt{E}a^\texttt{NE} + (q, q)$. By the simplicity property of $D$, no self-intersection can occur within $\kappa([b^\texttt{E}_\text{Mid}; b^\texttt{E}_\text{End}])$ and $\kappa([o; b^\texttt{NE}_\text{Mid}])$. Therefore, any self-intersection within $\kappa([b^\texttt{E}_\text{Mid}; b^\texttt{NE}_\text{Mid}])$ needs to occur between $\kappa([b^\texttt{E}_\text{Mid}; o])$ and $\kappa([b^\texttt{E}_\text{End}; b^\texttt{NE}_\text{Mid}])$.

However, by the separation property of $D$, the strip bounded by the lines $y = -x + p + q$ and $y = -x + 2q$ separates the former and the latter. This settles (e).

For (f), let \[c^\texttt{E} = -a^\texttt{W} + (p + q, 0)\] and \[c^\texttt{NE} = a^\texttt{S}a^\texttt{SE}a^\texttt{E}a^\texttt{NE}a^\texttt{N}a^\texttt{NW}a^\texttt{W} + (q, q),\] and define $c^\texttt{N}$, $c^\texttt{NW}$, \ldots, $c^\texttt{SE}$ symmetrically.

By the translation property of $D$, $h(D) = c^\texttt{E}c^\texttt{NE}\allowbreak\ldots c^\texttt{SE}$, with $c^i$ a subpath of $b^i$ for $i = \texttt{E}$, \texttt{N}, \texttt{W}, and \texttt{S}, and $b^i$ a subpath of $c^i$ for $i = \texttt{NE}$, \texttt{NW}, \texttt{SW}, and \texttt{SE}.

Let $i$ be any direction out of \texttt{E}, \texttt{NE}, \ldots, \texttt{SE}.

By the coherence property of $D$, no illicit pair can occur within $b^i$ and $c^i$.

Suppose that an illicit pair has one square $a$ in $c^\texttt{E}$ and another square $b$ in $c^\texttt{NE}c^\texttt{N}$.

Since no illicit pair can occur within $b^\texttt{E}$, $b$ does not belong to $b^\texttt{E}$.

By the separation property of $D$, the strip bounded by the lines $y = -x + p + q$ and $y = -x + 2q$ strictly separates $c^\texttt{E}$ and $b^\texttt{NE}$. However, since $r + s \le q - p$, no illicit pair can have its two squares on different sides of this strip. Thus $b$ does not belong to $b^\texttt{NE}$.

Since $c^\texttt{E}$ is a subpath of $D + (p + q, 0)$ and $b^\texttt{N}$ is a subpath of $D + (0, p + q)$, by the symmetry and protocoherence properties of $D$ no square of $c^\texttt{E}$ is joined by a move of $M$ to a square of $b^\texttt{N}$. Thus $b$ does not belong to $b^\texttt{N}$ either.

Since $c^\texttt{NE}c^\texttt{N}$ is a subpath of $b^\texttt{E}b^\texttt{NE}b^\texttt{N}$, we have arrived at a contradiction. By symmetry, it follows that no illicit pair can occur within $c^i \cup c^{i + 1}$ for any $i$ and $c^i \cup c^{i + 2}$ for $i = \texttt{E}$, \texttt{N}, \texttt{W}, and \texttt{S}.

Suppose that an illicit pair has one square $a$ in $c^\texttt{NE}$ and another square $b$ in $c^\texttt{NW}$.

By Lemma \ref{bound}, the strip bounded by the lines $y = x$ and $y = x + q - p$ strictly separates $b^\texttt{E}$ and $c^\texttt{NW}$. However, since $r + s \le q - p$, no illicit pair can have its two squares on different sides of this strip.

By the separation property of $D$ and Lemma \ref{bound}, the strip bounded by the lines $y = -x + p + q$ and $y = -x + 2q$ strictly separates $b^\texttt{NE}$ and $c^\texttt{NW}$. As above, no illicit pair can have its two squares on different sides of this strip.

Symmetrically, no illicit pair has one square in $c^\texttt{NE}$ and another square in either $b^\texttt{NW}$ or $b^\texttt{W}$.

We conclude that $a$ must belong to $c^\texttt{NE} \setminus b^\texttt{E}b^\texttt{NE}$ and $b$ must belong to $c^\texttt{NW} \setminus b^\texttt{NW}b^\texttt{W}$. However, both of those are subpaths of $b^\texttt{N}$, which contains no illicit pair.

We have arrived at a contradiction. By symmetry, it follows that no illicit pair can occur within $c^i \cup c^{i + 2}$ for $i = \texttt{NE}$, \texttt{NW}, \texttt{SW}, and \texttt{SE}.

Lastly, by Lemma \ref{bound}, the strip bounded by the lines $x = 0$ and $x = p$ strictly separates $c^\texttt{E}$ from $c^\texttt{NW}$ and $c^\texttt{W}$, and the strip bounded by the lines $y = -x$ and $y = -x + q - p$ strictly separates $c^\texttt{NE}$ from $c^\texttt{W}$ and $c^\texttt{SW}$. As above, by symmetry it follows that no illicit pair can occur within $c^i \cup c^{i + 3}$ and $c^i \cup c^{i + 4}$.

This completes the proof of (f).

For (g), by the separation property of $D$ and Lemma \ref{bound} the strip bounded by the lines $x = 2q$ and $x = 2p + 2q$ strictly separates $h(D)$ and $h(D) + (2p + 4q, 0)$. However, since $s \le p < 2p$, no two squares on different sides of this strip can be joined by a move of $M$.

Suppose, then, that a square $a$ in $h(D)$ is joined by a move of $M$ to a square $b$ in $h(D) + (p + 3q, p + 3q)$.

By the separation property of $D$ and Lemma \ref{bound}, $h(D) + (p + 3q, p + 3q)$ lies strictly to the right of the line $x = p + q$ and strictly above the line $y = p + q$. Since $s \le p$, it follows from this that $a$ lies in the portion of $h(D)$ strictly to the right of the line $x = q$ and strictly above the line $y = q$. By the separation property of $h(D)$ and Lemma \ref{bound}, this portion of $h(D)$ is a subpath of $c^\texttt{NE}$ and thus of $D + (q, q)$.

Analogously, $b$ lies in $D + (p + 2q, p + 2q)$. However, by the protocoherence property of $D$, no square of $D + (q, q)$ is joined by a move of $M$ to a square of $D + (p + 2q, p + 2q)$. We have arrived at a contradiction and the proof of (g) is complete.

\medskip

\emph{Case $\hbar$.} This case is analogous to Case $h$. \end{proof}

\medskip

Lemma \ref{fghh} supplies the induction step in the proof of Lemma \ref{perfect} below.

\medskip

\begin{definition} Let $M$ be an $(r, s)$-leaper, $r \le s$. Define the \emph{initial cycles} of $M$ as follows.

Provided that $r \neq 0$, the initial cycle of $M$ of type \texttt{f}, $D^M_\texttt{f}$, is the $M$-instantiation, centered at $O$, of the direction graph represented by the oriented cycle labeled \texttt{43658721}. It contains the eight squares $(r + s, 0)$, $(r, r)$, $(0, r + s)$, $(-r, r)$, $(-r - s, 0)$, $(-r, -r)$, $(0, -r - s)$, and $(r, -r)$.

Provided that $r \neq s$, the initial cycle of $M$ of type \texttt{g}, $D^M_\texttt{g}$, is the $M$-instantiation, centered at $O$, of the direction graph represented by the oriented cycle labeled \texttt{25476183}. It contains the eight squares $(s - r, 0)$, $(s, s)$, $(0, s - r)$, $(-s, s)$, $(r - s, 0)$, $(-s, -s)$, $(0, r - s)$, and $(s, -s)$.

The initial cycle of $M$ of type \texttt{h}, $D^M_\texttt{h}$, is the $M$-instantiation, centered at $O$, of the direction graph represented by the oriented cycle labeled \texttt{34567812}. It contains the eight squares $(r + s, 0)$, $(s, s)$, $(0, r + s)$, $(-s, s)$, $(-r - s, 0)$, $(-s, -s)$, $(0, -r - s)$, and $(s, -s)$. \end{definition}

\medskip

\begin{lemma} Let $M$ be an $(r, s)$-leaper, $r \le s$, $o$ one of the characters \texttt{f}, \texttt{g}, and \texttt{h} such that there exists an initial cycle of $M$ of type $o$, and $e = e_1e_2\ldots e_l$ a string composed of the characters \texttt{f}, \texttt{g}, \texttt{h}, and $\tthbar$.

Let $A_\tthbar = A_\texttt{h}$, $A_e = A_{e_1}A_{e_2}\ldots A_{e_l}$, and \[\left(\begin{array}{c} p\\ q \end{array}\right) = A_eA_o\left(\begin{array}{c} r\\ s \end{array}\right).\]

Then successively applying lifting transformations of types $e_l$, $e_{l - 1}$, \ldots, $e_1$ to $D^M_o$ yields a $(p, q)$-perfect cycle $D^M_o(e)$. \label{perfect} \end{lemma}

\medskip

\begin{proof} We proceed by induction on $e$.

When $e$ is the empty string, the lemma holds by the definition of an initial cycle.

Suppose, then, that the lemma holds for $e$. Let $e_0$ be one of the characters \texttt{f}, \texttt{g}, \texttt{h}, and $\tthbar$ and \[\left(\begin{array}{c} p'\\ q' \end{array}\right) = A_{e_0}\left(\begin{array}{c} p\\ q \end{array}\right).\]

We need to show that applying a lifting transformation of type $e_0$ to $D^M_o(e)$ yields a $(p', q')$-perfect cycle $D^M_o(e_0e)$.

\medskip

\emph{Case 1.} $s > p$ or $r + s > q - p$. By Lemma \ref{fghh}, it suffices to show that $D^M_o(e_0e)$ possesses the coherence and protocoherence properties.

By induction on $e$, if $s > p$ then $o = \texttt{f}$ and $e$ is a run of the character \texttt{f}. Analogously, if $r + s > q - p$ then $o = \texttt{g}$ and $e$ is a run of the character \texttt{g}.

We consider the former case in detail, and the latter case is analogous.

Let $o = \texttt{f}$ and $e$ be a run of the character \texttt{f} of length $n$.

Define a \emph{horizontal band} $|a; b|$ of endpoints $a$ and $b$ and length $m + 1$, with $a$ and $b$ squares and $m$ a nonnegative integer such that $b = a + m \cdot (2r, 0)$, as the set of all squares of the form $a + i \cdot (2r, 0)$, $i = 0$, 1, \ldots, $m$. Analogously, define a \emph{vertical band} $|a; b|$ of endpoints $a$ and $b$ and length $m + 1$, with $a$ and $b$ squares and $m$ a nonnegative integer such that $b = a + m \cdot (0, 2r)$, as the set of all squares of the form $a + i \cdot (0, 2r)$, $i = 0$, 1, \ldots, $m$.

(In the case of $o = \texttt{g}$ and $e$ being a run of the character \texttt{g}, \emph{forward} and \emph{backward} bands are defined analogously, with the translations $(2r, 0)$ and $(0, 2r)$ replaced by $(s - r, s - r)$ and $(r - s, s - r)$.)

By induction on $n$, the vertex set of $D^M_\texttt{f}(e_0e)$ is the union of a number of bands as follows. In each case, ``rotations'' stands for ``multiple-of-quarter-turn rotations about $O$'' and ``reflections'' stands for ``reflections in the lines $x = 0$, $y = 0$, and $y = \pm x$''. Furthermore, in each case we classify all bands in question as either \emph{inner} or \emph{outer} ones.

The vertex set of $D^M_\texttt{f}(\texttt{f}e)$ is the union of four inner bands of length $n + 3$, namely the rotations of \[|((n + 2)r, -(n + 2)r); ((n + 2)r, (n + 2)r)|,\] and four outer bands of length $n + 2$, namely the rotations of \[|((n + 2)r + s, -(n + 1)r); ((n + 2)r + s, (n + 1)r)|.\]

The vertex set of $D^M_\texttt{f}(\texttt{g}e)$ is the union of twelve inner bands of lengths $n + 1$ and $n + 2$, namely the rotations of \[|((3n + 2)r + s, -(n + 1)r); ((3n + 2)r + s, (n + 1)r)|\] and the rotations and reflections of \[|((n + 2)r + s, (n + 1)r); ((3n + 2)r + s, (n + 1)r)|,\] and twelve outer bands of length $n + 1$, namely the rotations of \[|((3n + 2)r + 2s, -nr); ((3n + 2)r + 2s, nr)|\] and the rotations and reflections of  \[|((n + 1)r + s, (n + 1)r + s); ((3n + 1)r + s, (n + 1)r + s)|.\]

The vertex set of $D^M_\texttt{f}(\texttt{h}e)$ is the union of twenty inner bands of lengths $n + 1$ and $n + 2$, namely the rotations of \[|((n + 2)r, -nr); ((n + 2)r, nr)|\] and the rotations and reflections of \[|((n + 3)r + s, (n + 1)r + s); ((3n + 3)r + s, (n + 1)r + s)|\] and \[|((3n + 3)r + s, (n + 1)r + s); ((3n + 3)r + s, (3n + 3)r + s)|,\] and twenty outer bands of lengths $n + 1$ and $n + 2$, namely the rotations of \[|((n + 2)r + s, -(n + 1)r); ((n + 2)r + s, (n + 1)r)|\] and the rotations and reflections of \[|((n + 2)r + s, (n + 1)r); ((3n + 2)r + s, (n + 1)r)|\] and \[|((3n + 3)r + 2s, (n + 2)r + s); ((3n + 3)r + 2s, (3n + 2)r + s)|.\]

\begin{figure}[ht!]\hspace*{\fill}\includegraphics[width=\textwidth]{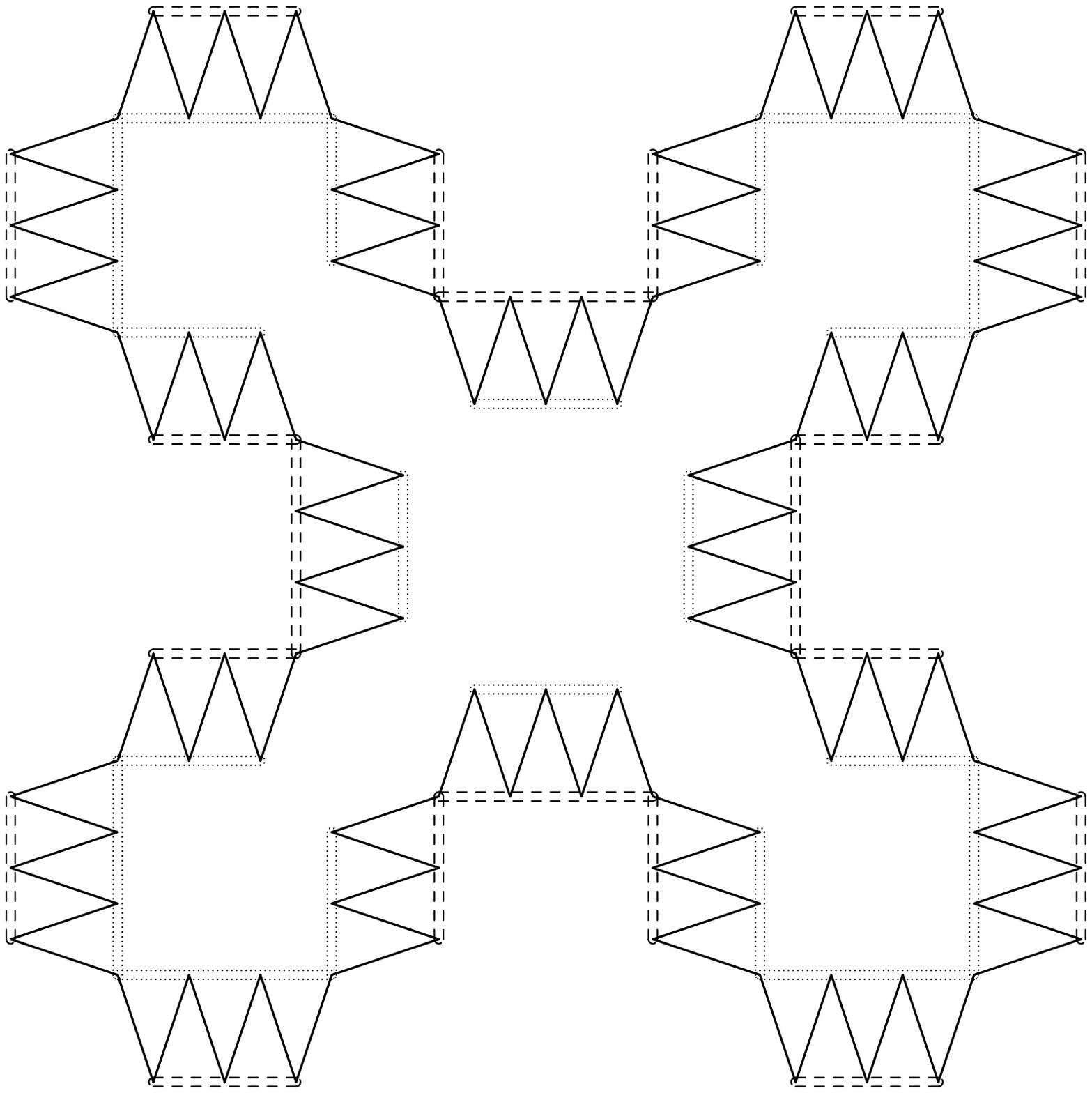}\hspace*{\fill}\caption{}\label{bands}\end{figure}

Lastly, the vertex set of $D^M_\texttt{f}(\tthbar e)$ is the union of twelve inner bands of length $n + 2$, namely the rotations of \[|((3n + 4)r + s, -(n + 1)r); ((3n + 4)r + s, (n + 1)r)|\] and the rotations and reflections of \[|((n + 2)r + s, (n + 1)r); ((3n + 4)r + s, (n + 1)r)|,\] and twelve outer bands of lengths $n + 1$ and $n + 2$, namely the rotations of \[|((3n + 4)r + 2s, -nr); ((3n + 4)r + 2s, nr)|\] and the rotations and reflections of  \[|((n + 1)r + s, (n + 1)r + s); ((3n + 3)r + s, (n + 1)r + s)|.\]

Figure \ref{bands} shows the bands making up $D^{(1, 3)}_\texttt{f}(\texttt{h}\texttt{f}\texttt{f})$, corresponding to $M$ being a $(1, 3)$-leaper, $o = \texttt{f}$, $n = 2$, $e = \texttt{f}\texttt{f}$, and $e_0 = \texttt{h}$, with all inner bands dotted and all outer bands dashed.

To establish the coherence of $D^M_o(e_0e)$, it suffices to show that it does not contain an illicit pair.

First we show that no two bands of the same type contain an illicit pair for reasons of parity.

Namely, by tracing $D^M_o(e_0e)$, $M$ travels from every square in an inner band to every other square in an inner band in an even number of moves, and similarly for outer bands. Since no odd cycle of a leaper exists, two squares in bands of the same type are never joined by a move of $M$.

Then we show that no two bands of opposite types contain an illicit pair because, roughly, they are too far away.

There are three kinds of relative positions that bands of opposite types in $D^M_o(e_0e)$ occupy, as follows.

Two bands of opposite types are \emph{tied} if they run in parallel, at a distance of $s$, and $D^M_o(e_0e)$ zigzags between them, visiting all of their squares in succession.

Two bands of opposite types are \emph{adjacent} if an edge of $D^M_o(e_0e)$ joins an endpoint $a$ of the first band to an endpoint $b$ of the second band, and no other edges of $D^M_o(e_0e)$ join any of the bands' squares. Two adjacent bands always make obtuse angles with the segment joining the centers of $a$ and $b$.

Two bands of opposite types are \emph{independent} if no edge of $D^M_o(e_0e)$ joins a square of one band to a square of the other.

All three kinds of relative positions are exemplified in Figure \ref{bands}.

No illicit pair can occur within two tied or adjacent bands.

Define the \emph{envelope} of a band $|a; b|$ as the regular chessboard of lower left corner $a + (-s, -s)$ and upper right corner $b + (s, s)$. All squares joined by a move of $M$ to a square in a band $B$ belong to the envelope of $B$.

(In the case of $o = \texttt{g}$ and $e$ being a run of the character \texttt{g}, define the envelope of a band $|a; b|$ as the set of all squares inside the convex hull of the centers of $a + (\pm(r + s), 0)$, $a + (0, \pm(r + s))$, $b + (\pm(r + s), 0)$, and $b + (0, \pm(r + s))$.)

For every band $B$ in $D^M_o(e_0e)$, the envelope of $B$ does not intersect any bands independent from $B$, except possibly at endpoints that coincide with the endpoints of bands tied or adjacent to $B$.

It follows, then, that no illicit pair can occur within two independent bands either, and that $D^M_o(e_0e)$ does not contain an illicit pair.

We are left to establish the protocoherence of $D^M_o(e_0e)$.

Consider first $D^M_o(e_0e)$ and $D^M_o(e_0e) + (2q', 0)$. Suppose that a square $a$ in the former is joined by a move of $M$ to a square $b$ in the latter.

Let $B_\texttt{E}$ be the union of all eastmost bands of $D^M_o(e_0e)$ (there is one such band when $e_0 = \texttt{f}$, \texttt{g}, or $\tthbar$, and two when $e_0 = \texttt{h}$), and define $B_\texttt{W}$ symmetrically. Furthermore, let $B'_\texttt{W} = B_\texttt{W} + (2q', 0)$ be the image of $B_\texttt{W}$ in $D^M_o(e_0e) + (2q', 0)$.

Since all squares of $D^M_o(e_0e) + (2q', 0)$ apart from the ones in $B'_\texttt{W}$ lie more than $s$ units to the east of $D^M_o(e_0e)$, $b$ belongs to $B'_\texttt{W}$.

Symmetrically, $a$ belongs to $B_\texttt{E}$.

Since both of $B_\texttt{E}$ and $B_\texttt{W}$ consist of outer bands and, by induction on $n$, the translation $(2q', 0)$ is a sum of an even number of $M$-translations, it takes $M$ an even number of moves to get from $a$ to $b$. Therefore, as above, $a$ and $b$ cannot be joined by a move of $M$. We have arrived at a contradiction.

Consider, then, $D^M_o(e_0e)$ and $D^M_o(e_0e) + (p' + q', p' + q')$.

In both cases $e_0 = \texttt{f}$ and $e_0 = \texttt{h}$, an analogous argument applies as follows.

Suppose that a square $a$ in $D^M_o(e_0e)$ is adjacent to a square $b$ in $D^M_o(e_0e) + (p' + q', p' + q')$.

Define $B_\texttt{N}$ and $B_\texttt{S}$ symmetrically to $B_\texttt{E}$ and $B_\texttt{W}$. Let $B''_\texttt{W} = B_\texttt{W} + (p' + q', p' + q')$ and $B''_\texttt{S} = B_\texttt{S} + (p' + q', p' + q')$ be the images of $B_\texttt{W}$ and $B_\texttt{S}$ in $D^M_o(e_0e) + (p' + q', p' + q')$.

Since all squares of $(D^M_o(e_0e) + (p' + q', p' + q')) \setminus B''_\texttt{W}$ lie at least $r + s$ units to the east of all squares of $D^M_o(e_0e) \setminus B_\texttt{E}$, either $a$ belongs to $B_\texttt{E}$ or $b$ belongs to $B''_\texttt{W}$.

Symmetrically, either $a$ belongs to $B_\texttt{N}$ or $b$ belongs to $B''_\texttt{S}$.

It follows that $a$ lies in the union of $B_\texttt{E}$ and $B_\texttt{N}$ and $b$ lies in the union of $B''_\texttt{W}$ and $B''_\texttt{S}$.

Since all of $B_\texttt{E}$, $B_\texttt{N}$, $B_\texttt{W}$, and $B_\texttt{S}$ consist of outer bands and, by induction on $n$, the translation $(p' + q', p' + q')$ is a sum of an even number of $M$-translations, it takes $M$ an even number of moves to get from $a$ to $b$. Therefore, as above, $a$ and $b$ cannot be joined by a move of $M$. We have arrived at a contradiction.

The cases $e_0 = \texttt{g}$ and $e_0 = \tthbar$ are handled as in the proof of part (g) of Lemma \ref{fghh} in Cases $g$ and $\hbar$. Namely, when $e_0 = \texttt{g}$ or $e_0 = \tthbar$, $D^M_o(e_0e)$ and $D^M_o(e_0e) + (p' + q', p' + q')$ are strictly separated by a diagonal strip of width $\sqrt{2}(q - p) \ge \sqrt{2}(r + s) > \frac{1}{\sqrt{2}}(r + s)$.

\medskip

\emph{Case 2.} $s \le p$ and $r + s \le q - p$. Then $D^M_o(e_0e)$ is a $(p', q')$-perfect cycle by Lemma \ref{fghh}. \end{proof}

\medskip

We have set up all the tools we need and are ready to construct dual boards by means of $f$, $g$, $h$, and $\hbar$.

Given $M$, $o$, $e$, $p$, and $q$ as in the setting of Lemma \ref{perfect}, let $L^M_o(e)$ be a $(p, q)$-leaper and $B^M_o(e)$ the board formed by the squares of $D^M_o(e)$.

\medskip

\begin{theorem} Let $M$ be an $(r, s)$-leaper, $r \le s$, $o$ one of the characters \texttt{f}, \texttt{g}, and \texttt{h} such that there exists an initial cycle of $M$ of type $o$, and $e = e_1e_2\ldots e_l$ a string composed of the characters \texttt{f}, \texttt{g}, \texttt{h}, and $\tthbar$. Then the board $B^M_o(e)$ is dual with respect to $M$ and $L^M_o(e)$. \label{fghhboard} \end{theorem}

\medskip

\begin{proof} By Lemma \ref{perfect}, $D^M_o(e)$ is a perfect cycle. By the coherence property of $D^M_o(e)$, the leaper graph of $M$ over $B^M_o(e)$ is a cycle. By Lemma \ref{longcoh} and analogously to the proofs of Lemmas \ref{f1}, \ref{g1}, and \ref{h1}, so is the leaper graph of $L^M_o(e)$ over $B^M_o(e)$. \end{proof}

\medskip

In particular, we obtain the following corollary.

\medskip

\begin{corollary} Let $L$ and $M$ be two distinct skew basic leapers such that the descent of one of them is a suffix of the descent of the other. Then there exists a board dual with respect to $L$ and $M$. \label{fghboard} \end{corollary}

\medskip

\begin{proof} Let the descent of $L$ be $e_1e_2 \ldots e_l$ and that of $M$ $e_{m + 1}e_{m + 2}\ldots e_l$, with $m \le l$. By Theorem \ref{fghhboard}, $B^M_{e_m}(e_1e_2\ldots e_{m - 1})$ is a board dual with respect to $L$ and $M$. \end{proof}

\medskip

We go on to extract dual direction graphs from the dual boards given by Theorem \ref{fghhboard}.

\medskip

\begin{theorem} Let $M$ be a skew $(r, s)$-leaper, $r < s$, $o$ one of the characters \texttt{f}, \texttt{g}, and \texttt{h}, and $e = e_1e_2\ldots e_l$ a string composed of the characters \texttt{f}, \texttt{g}, \texttt{h}, and $\tthbar$.

Then both the cycle $C^M_o(e)$ of $L^M_o(e)$ over $B^M_o(e)$ and the cycle $D^M_o(e)$ are trivial.

The direction graph extracted from $C^M_o(e)$ depends only on $e$. Extending the definition of a fundamental direction cycle, we refer to this direction graph as the fundamental direction cycle $\Phi(e)$ of descent $e$.

The direction graph extracted from $D^M_o(e)$ depends only on $o$ and $e$. Extending the definition of a second fundamental direction cycle, we refer to this direction graph as the second fundamental direction cycle $\Phi^\text{II}_o(e)$ of origin $o$ and descent $e$.

(The above continues to apply to non-skew leapers $M$ in the sense that $D^M_o(e)$ is an instantiation of $\Phi^\text{II}_o(e)$.)

Extend the definition of a flip by $\overline{\tthbar} = \tthbar$, and let the strings $e'$ and $e''$ composed of the characters \texttt{f}, \texttt{g}, \texttt{h}, and $\tthbar$ be flips of each other. Then the fundamental direction cycle of descent $e'$ and
the three second fundamental direction cycles of descent $e''$ are equivalent.

Furthermore, both of $\Phi(e)$ and $\Phi^\text{II}_o(e)$ are dual direction graphs. Each of them is a complement of the other, the duality matrix of $\Phi^\text{II}_o(e)$ being $A_eA_o$ and that of $\Phi(e)$ its inverse. \label{fghhgraph} \end{theorem}

\medskip

\begin{proof} As in the proofs of Theorems \ref{fund} and \ref{secondfund}, the moves of both $C^M_o(e)$ and $D^M_o(e)$ occur in pairs symmetric with respect to $O$ (those of $C^M_o(e)$ by induction on $e$ and those of $D^M_o(e)$ by its symmetry property) and since the sum of the associated direction matrices over each such pair is the zero matrix, both cycles are trivial.

Define the signature of $C^M_o(e)$ as in the proof of Theorem \ref{fund}. Furthermore, define an \emph{$\hbar$-rewrite} by means of the following system of rewriting rules. \begin{align*} +_\texttt{s} &\to -_\texttt{c}-_\texttt{s}-_\texttt{c} & +_\texttt{c} &\to +_\texttt{c}+_\texttt{s}+_\texttt{c}+_\texttt{s}+_\texttt{c}\\ -_\texttt{s} &\to +_\texttt{c}+_\texttt{s}+_\texttt{c} & -_\texttt{c} &\to -_\texttt{c}-_\texttt{s}-_\texttt{c}-_\texttt{s}-_\texttt{c} \end{align*}

The proof that the direction graph extracted from $C^M_o(e)$ depends only on $e$ is then analogous to the proof of the corresponding part of Theorem \ref{fund}.

Define the signatures of $D^M_\texttt{g}(e)$ and $D^M_\texttt{h}(e)$ as in the proof of Theorem \ref{secondfund}. Define the signature of $D^M_\texttt{f}(e)$ analogously, referring to the following table. \[ \begin{tabular}{c c} Label of $a$ & Directions of moves to and from $a$\\ \hline $+_\texttt{s}$ & \texttt{21}, \texttt{43}, \texttt{65}, \texttt{87}\\  $+_\texttt{c}$ & \texttt{14}, \texttt{36}, \texttt{58}, \texttt{72} \\ $-_\texttt{s}$ & \texttt{12}, \texttt{34}, \texttt{56}, \texttt{78} \\ $-_\texttt{c}$ & \texttt{27}, \texttt{41}, \texttt{63}, \texttt{85} \end{tabular} \]

Furthermore, define an \emph{$\hbar$-rearrangement} by \[(s^\text{Corner},\; s^\text{Side}) \to (\overline{s^\text{Side}}\,\overline{s^\text{Corner}}\,\overline{s^\text{Side}},\; s^\text{Side}s^\text{Corner}s^\text{Side}s^\text{Corner}s^\text{Side}).\]

The proof that the direction graph extracted from $D^M_o(e)$ depends only on $o$ and $e$ is then analogous to the proof of the corresponding part of Theorem \ref{secondfund}.

That the fundamental direction cycle of descent $e'$ and the three second fundamental direction cycles of descent $e''$ are equivalent is established as in the proof of Theorem \ref{flip}, the equivalence permutation \[\pi_\texttt{f} = \left(\begin{array}{cccccccc} \texttt{1} & \texttt{2} & \texttt{3} & \texttt{4} & \texttt{5} & \texttt{6} & \texttt{7} & \texttt{8}\\ \texttt{2} & \texttt{3} & \texttt{8} & \texttt{1} & \texttt{6} & \texttt{7} & \texttt{4} & \texttt{5} \end{array}\right),\] induced by the matrices \[P_\texttt{f} = \left(\begin{array}{cc} -\;\frac{1}{\sqrt{2}} & \frac{1}{\sqrt{2}}\\[5pt] \frac{1}{\sqrt{2}} & \frac{1}{\sqrt{2}} \end{array}\right)\] and \[Q_\texttt{f} = \left(\begin{array}{cc} \frac{1}{\sqrt{2}} & \frac{1}{\sqrt{2}}\\[5pt] -\;\frac{1}{\sqrt{2}} & \frac{1}{\sqrt{2}} \end{array}\right),\] added alongside $\pi_\texttt{g}$ and $\pi_\texttt{h}$.

Lastly, the derivation of the duality identity for $\Phi(e)$ and $\Phi^\text{II}_o(e)$ is analogous to the derivation of the duality identity in the proof of Theorem \ref{dualfund}. The concluding part of the theorem then follows from that identity and Theorem \ref{dualbasic} as in the proof of Theorem \ref{dualfund}. \end{proof}

\medskip

It can be demonstrated that the family of all dual direction graphs given by Theorem \ref{fghhgraph} is closed under equivalence. In fact, it coincides with the closure of the family of all dual direction graphs given by Theorem \ref{dualfund} under equivalence and the extended $f$, $g$, and $h$.

Theorems \ref{fghhboard} and \ref{fghhgraph} add a great variety of dual boards and dual direction graphs to the ones given by Theorems \ref{dual} and \ref{dualfund}.

Let us look more closely into what dual direction graphs given by Theorem \ref{fghhgraph} are essentially different from all dual direction graphs given by Theorem \ref{dualfund}, in the sense of not being equivalent to any of them.

Since the latter family is a subset of the former, it suffices to obtain a necessary and sufficient condition for two dual direction graphs given by Theorem \ref{fghhgraph} to be equivalent. Since every second fundamental direction cycle is equivalent to the fundamental direction cycle of flipped descent, the question reduces to obtaining a necessary and sufficient condition for two fundamental direction cycles to be equivalent.

First we show that two fundamental direction cycles never coincide except trivially.

\medskip

\begin{theorem} Two fundamental direction cycles coincide if and only if their descents do. \label{fghheq} \end{theorem}

\medskip

\begin{proof} For all strings $e$ composed of the characters \texttt{f}, \texttt{g}, \texttt{h}, and $\tthbar$, let $n_\texttt{s}(e)$ be the number of occurrences of the character $+_\texttt{s}$ in the signature of $\Phi(e)$, and $n_\texttt{c}(e)$ the number of occurrences of the character $+_\texttt{c}$.

When $e$ is the empty string, $n_\texttt{s}(e) = n_\texttt{c}(e) = 4$.

By the definition of an $f$-rewrite, \[n_\texttt{s}(\texttt{f}e) = n_\texttt{s}(e) + 2n_\texttt{c}(e)\] and \[n_\texttt{c}(\texttt{f}e) = n_\texttt{c}(e).\]

Analogously, by the definition of a $g$-rewrite \[n_\texttt{s}(\texttt{g}e) = n_\texttt{s}(e)\] and \[n_\texttt{c}(\texttt{g}e) = 2n_\texttt{s}(e) + n_\texttt{c}(e),\] by the definition of an $h$-rewrite \[n_\texttt{s}(\texttt{h}e) = 3n_\texttt{s}(e) + 2n_\texttt{c}(e)\] and \[n_\texttt{c}(\texttt{h}e) = 2n_\texttt{s}(e) + n_\texttt{c}(e),\] and by the definition of an $\hbar$-rewrite \[n_\texttt{s}(\tthbar e) = n_\texttt{s}(e) + 2n_\texttt{c}(e)\] and \[n_\texttt{c}(\tthbar e) = 2n_\texttt{s}(e) + 3n_\texttt{c}(e).\]

By induction on $e$, then, $n_\texttt{s}(e) \ge 4$ and $n_\texttt{c}(e) \ge 4$ for all $e$.

Let $e_0$ be one of the characters \texttt{f}, \texttt{g}, \texttt{h}, and $\tthbar$, $a = n_\texttt{s}(e_0e)$, and $b = n_\texttt{c}(e_0e)$.

It follows that $e_0 = \texttt{f}$ if and only if \[a > 2b,\] $e_0 = \texttt{g}$ if and only if \[2a < b,\] $e_0 = \texttt{h}$ if and only if \[2b > a > b,\] and $e_0 = \tthbar$ if and only if \[a < b < 2a.\]

Suppose that the fundamental direction cycles $\Phi(e')$ and $\Phi(e'')$ coincide, with $e' = e'_1e'_2 \ldots e'_l$ and $e'' = e''_1e''_2\ldots e''_m$. Then their signatures coincide as well. Consequently, $n_\texttt{s}(e') = n_\texttt{s}(e'')$ and $n_\texttt{c}(e') = n_\texttt{c}(e'')$.

If $n_\texttt{s}(e') = n_\texttt{s}(e'') = 4$ and $n_\texttt{c}(e') = n_\texttt{c}(e'') = 4$, then both of $e'$ and $e''$ are the empty string.

Otherwise, by the above analysis we conclude that $e'_1 = e''_1$, $n_\texttt{s}(e'_2e'_3\ldots e'_l) = n_\texttt{s}(e''_2e''_3\ldots e''_m)$, and $n_\texttt{c}(e'_2e'_3\ldots e'_l) = n_\texttt{c}(e''_2e''_3\ldots e''_m)$.

Iteration successively yields $e'_2 = e''_2$, $e'_3 = e''_3$, \ldots, $l = m$, $e'_l = e''_l$, and $e' = e''$. \end{proof}

\medskip

\begin{figure}[t!]\hspace*{\fill}\includegraphics[width=0.9\textwidth]{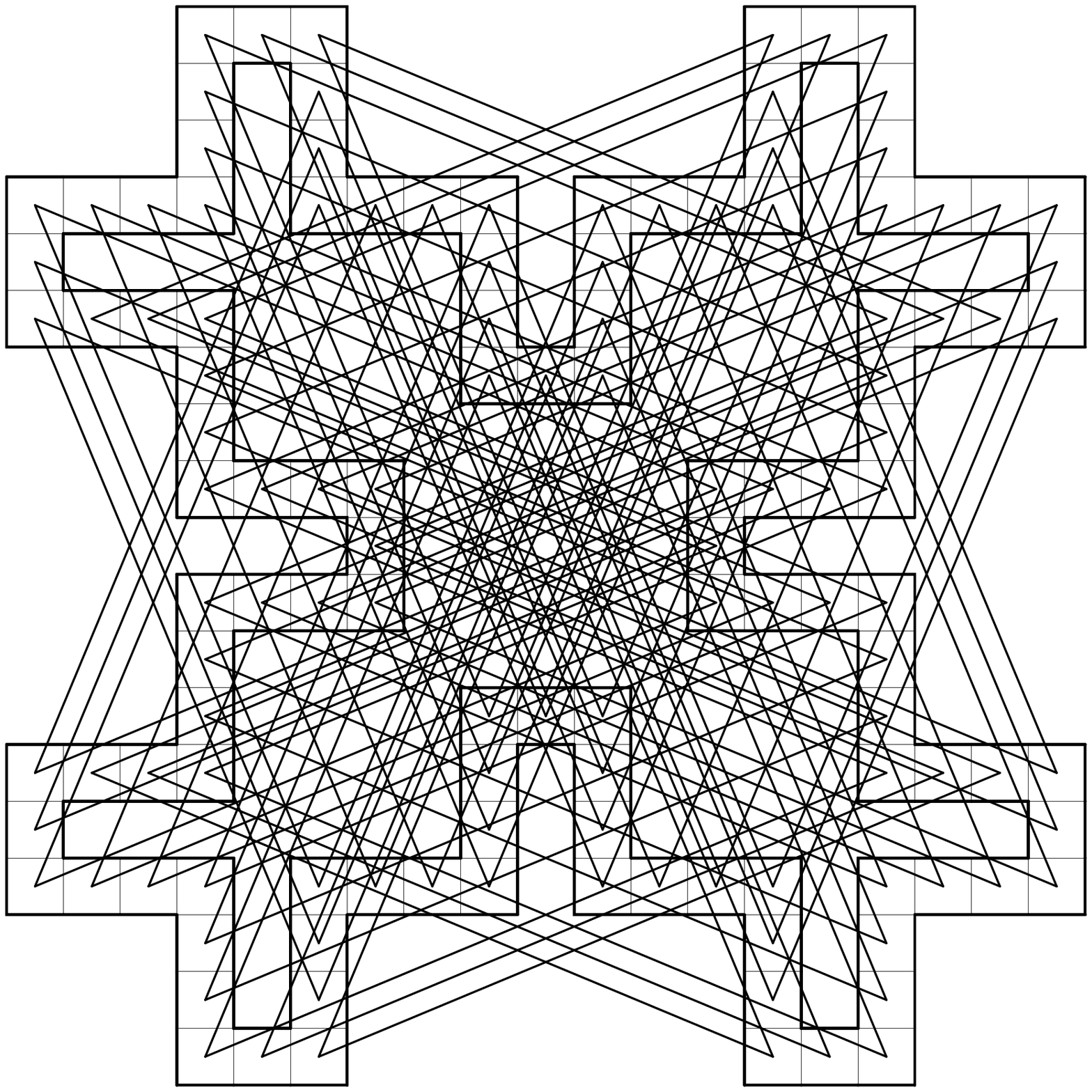}\hspace*{\fill}\caption{}\label{hhbar}\end{figure}

\begin{figure}[t!]\hspace*{\fill}\includegraphics[width=\textwidth]{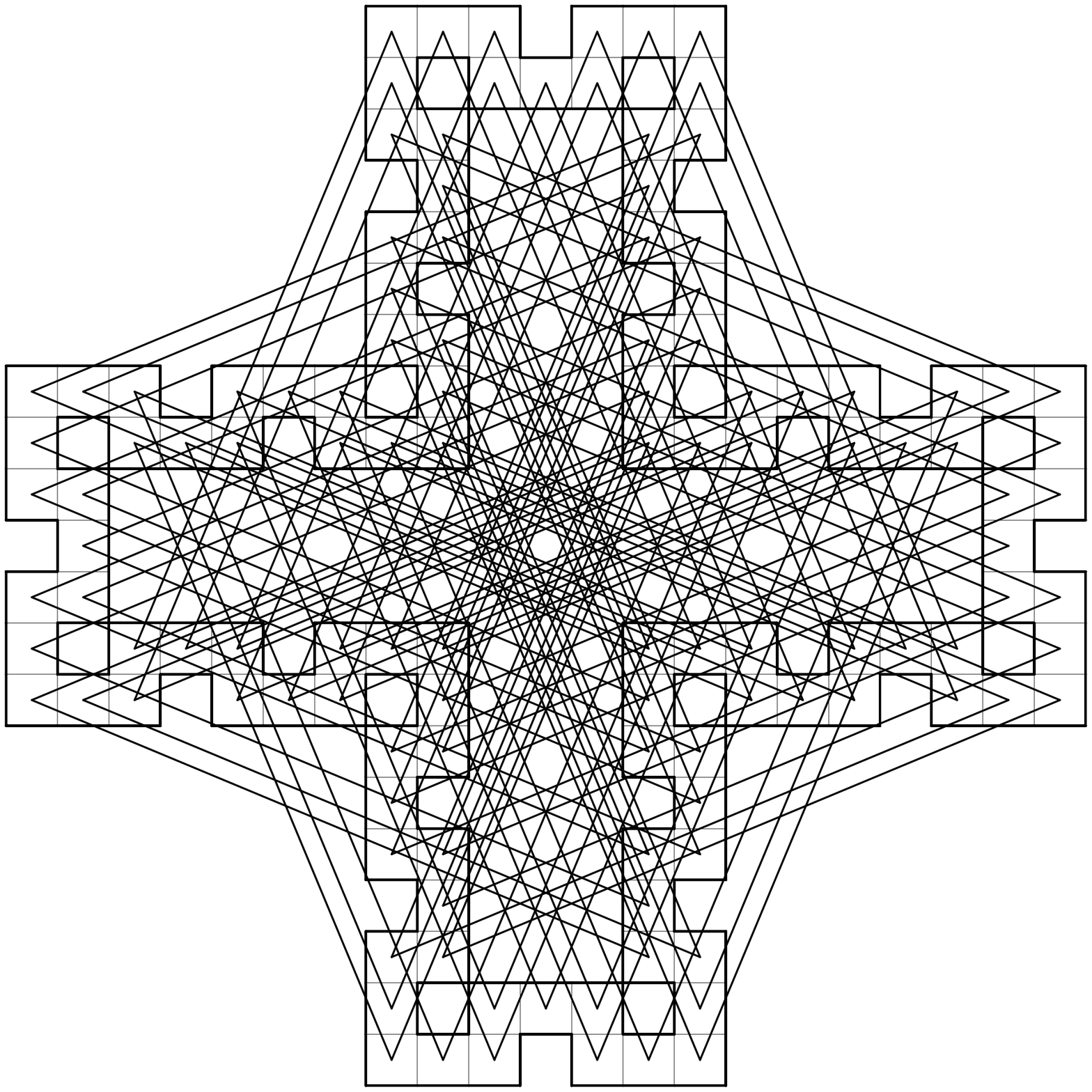}\hspace*{\fill}\caption{}\label{hbarh}\end{figure}

\begin{definition} Let $e$ be a string composed of the characters \texttt{f}, \texttt{g}, \texttt{h}, and $\tthbar$. The \emph{companion} of $e$ is the string obtained from $e$ by replacing every \texttt{f} with a \texttt{g} and vice versa, and every \texttt{h} with an $\tthbar$ and vice versa. \end{definition}

\medskip

A string $e'$ is the companion of $e''$ if and only if $e''$ is the companion of $e'$.

\medskip

\begin{theorem} Let $e'$ and $e''$ be two strings composed of the characters \texttt{f}, \texttt{g}, \texttt{h}, and $\tthbar$. Then the fundamental direction cycles of descents $e'$ and $e''$ are equivalent if and only if either $e'$ and $e''$ coincide, or $e'$ and $e''$ are companions.~\label{fghhequiv} \end{theorem}

\medskip

\begin{proof} If $e'$ and $e''$ are companions, then by induction on $e'$ the equivalence permutation $\pi_\text{Shift}$ defined by $\pi_\text{Shift}(i) = i + 1$ for all skew directions $i$ and induced by the matrices \[P_\text{Shift} = \left(\begin{array}{cc} \frac{1}{\sqrt{2}} & -\;\frac{1}{\sqrt{2}}\\[5pt] \frac{1}{\sqrt{2}} & \frac{1}{\sqrt{2}} \end{array}\right)\] and \[Q_\text{Shift} = \left(\begin{array}{cc} -\;\frac{1}{\sqrt{2}} & \frac{1}{\sqrt{2}}\\[5pt] \frac{1}{\sqrt{2}} & \frac{1}{\sqrt{2}} \end{array}\right)\] maps $\Phi(e')$ onto $\Phi(e'')$.

Let, then, $\Phi(e')$ and $\Phi(e'')$ be equivalent and $\pi$ an equivalence permutation that maps $\Phi(e')$ onto $\Phi(e'')$.

Let $C'$ be an oriented cycle representing $\Phi(e')$, and define $C''$ analogously.

The directions of the moves leading to and from every vertex of $C'$ are of the form $i$ and $i \pm 3$, and similarly for $C''$. Furthermore, by induction on $e'$, for every skew direction $i$ there exists a vertex in $C'$ such that the directions of the moves to and from it coincide with $i$ and $i + 3$. (Though, possibly, not in that order.)

Therefore, $\pi$ maps every unordered pair of skew directions of the form $i$ and $i + 3$ to an unordered pair of the same form.

There exist sixteen such permutations of the eight skew directions, namely $\pi^k_\text{Shift}$ and $\pi_\text{Reflect} \circ \pi^k_\text{Shift}$ for $k = 0$, 1, \ldots, 7, where \[\pi_\text{Reflect} = \left(\begin{array}{cccccccc} \texttt{1} & \texttt{2} & \texttt{3} & \texttt{4} & \texttt{5} & \texttt{6} & \texttt{7} & \texttt{8}\\ \texttt{8} & \texttt{7} & \texttt{6} & \texttt{5} & \texttt{4} & \texttt{3} & \texttt{2} & \texttt{1} \end{array}\right).\]

All of those are equivalence permutations. Eight of them preserve $\Phi(e')$, and eight map it onto the fundamental direction cycle of descent the companion of $e'$.

By Theorem \ref{fghheq}, either $e'$ and $e''$ coincide or $e'$ and $e''$ are companions. \end{proof}

\medskip

\begin{corollary} A dual direction graph given by Theorem \ref{fghhgraph} is not equivalent to any dual direction graph given by Theorem \ref{dualfund} if and only if its descent contains both of the characters \texttt{h} and $\tthbar$. \end{corollary}

\medskip

The first dual direction graphs given by Theorem \ref{fghhgraph} that are not equivalent to any dual direction graph given by Theorem \ref{dualfund} are thus the ones of descents $\texttt{h}\tthbar$ and $\tthbar\texttt{h}$.

Figure \ref{hhbar} shows the board $B^{(0, 1)}_\texttt{h}(\texttt{h}\tthbar)$ overlaid with the associated cycle $C^{(0, 1)}_\texttt{h}(\texttt{h}\tthbar)$, the $(5, 12)$-instantiation of $\Phi(\texttt{h}\tthbar)$.

Similarly, Figure \ref{hbarh} shows the board $B^{(0, 1)}_\texttt{h}(\tthbar\texttt{h})$ overlaid with the associated cycle $C^{(0, 1)}_\texttt{h}(\tthbar\texttt{h})$, the $(5, 12)$-instantiation of $\Phi(\tthbar\texttt{h})$.

\section[Constructions of Dual Boards and Dual Direction Graphs II]{Constructions of Dual Boards\\ and Dual Direction Graphs II}

We conclude by exhibiting one more construction of dual boards and dual direction graphs.

Given a $(p, q)$-leaper $L$, introduce a Cartesian coordinate system $Oxy$ over the infinite chessboard such that the integer points are the vertices of the squares if $p + q$ is odd, and the centers of the squares if it is even, and write $(x, y)$ for the square centered at $(x, y)$.

Given a positive integer $d$ and real numbers $e_1$ and $e_2$, write $N_d(e_1, e_2)$ for the \emph{net} formed by all squares $(x, y)$ such that $x \equiv e_1$ and $y \equiv e_2$ modulo $d$.

\medskip

\begin{definition} Let $n$ be a positive integer and $L$ a $(p, q)$-leaper with $q \neq 0$.

Let $N^\text{I}(L)$ be the net \[N_{2q}\left(\frac{1}{2}(p + q), \frac{1}{2}(p - q)\right)\] and $W^\text{I}_n(L)$ the board formed by all squares of $N^\text{I}(L)$ whose centers lie on the boundary or in the interior of the rectangle $R^\text{I}_n(L)$ bounded by \[p \le x + y \le p + 2nq\] and \[-(2n - 1)q \le y - x \le (2n - 1)q.\]

Let $N^i(L)$, $R^i_n(L)$, and $W^i_n(L)$, for $i = \text{II}$, $\text{III}$, and $\text{IV}$, be the rotations of $N^\text{I}(L)$, $R^\text{I}_n(L)$, and $W^\text{I}_n(L)$ by $90^\circ$, $180^\circ$, and $270^\circ$ counterclockwise about $O$.

Then the \emph{pinwheel board} of order $n$ for $L$, $W_n(L)$, is the disjoint union of its four \emph{wings} $W^i_n(L)$, $i = \text{I}$, $\text{II}$, $\text{III}$, and $\text{IV}$. \end{definition}

\medskip

\begin{figure}[ht!]\hspace*{\fill}\includegraphics[width=\textwidth]{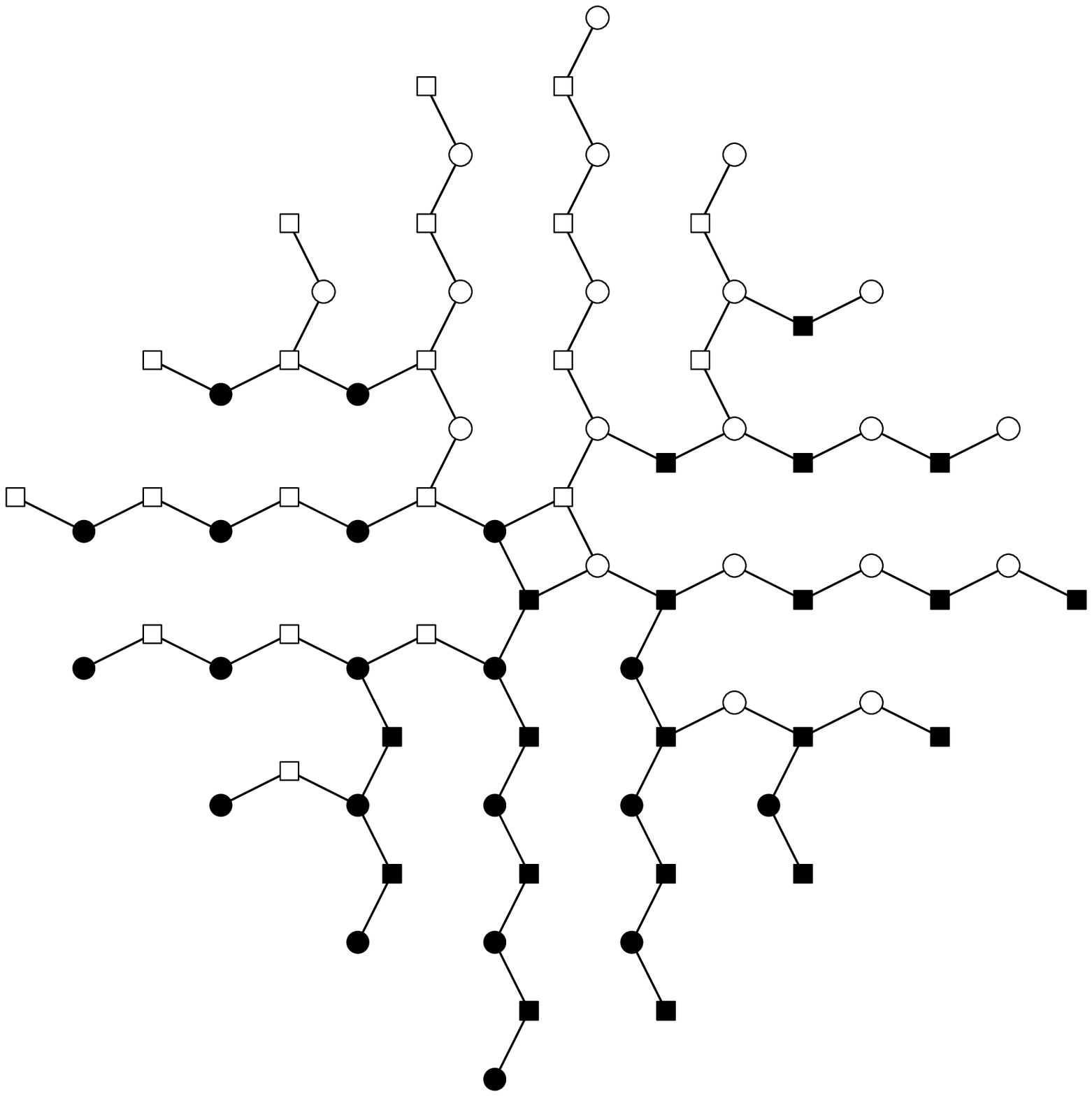}\hspace*{\fill}\caption{}\label{w412}\end{figure}

Figure \ref{w412} shows the centers of the squares of $W_4(1, 2)$, overlaid with the associated leaper graph and with the four wings labeled {\setlength{\fboxrule}{0pt}\fbox{\put(0, 4){\circle{8}}}}, $\square$, {\setlength{\fboxrule}{0pt}\fbox{\put(0, 4){\circle*{8}}}}, and $\blacksquare$.

\medskip

\begin{theorem} Let $n$ be a positive integer, \[A^\text{Pinwheel}_n = \left(\begin{array}{cc} 0 & 1\\1 & 2n \end{array}\right),\] $L$ a $(p, q)$-leaper with $p < q$, and $M$ an $(r, s)$-leaper with \[\left(\begin{array}{c} r\\ s \end{array}\right) = A^\text{Pinwheel}_n \left(\begin{array}{c} p\\ q \end{array}\right).\]

Then the pinwheel board of order $n$ for $L$ is dual with respect to $L$ and $M$.~\label{pinwheelboard} \end{theorem}

\medskip

\begin{figure}[t!]\hspace*{\fill}\includegraphics[scale=0.66]{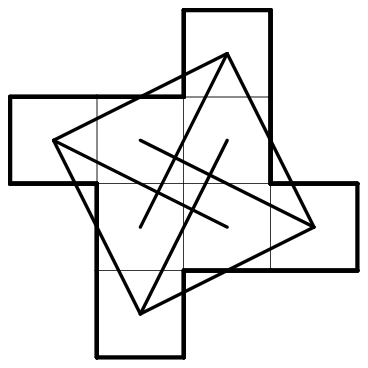}\hfill\includegraphics[scale=0.66]{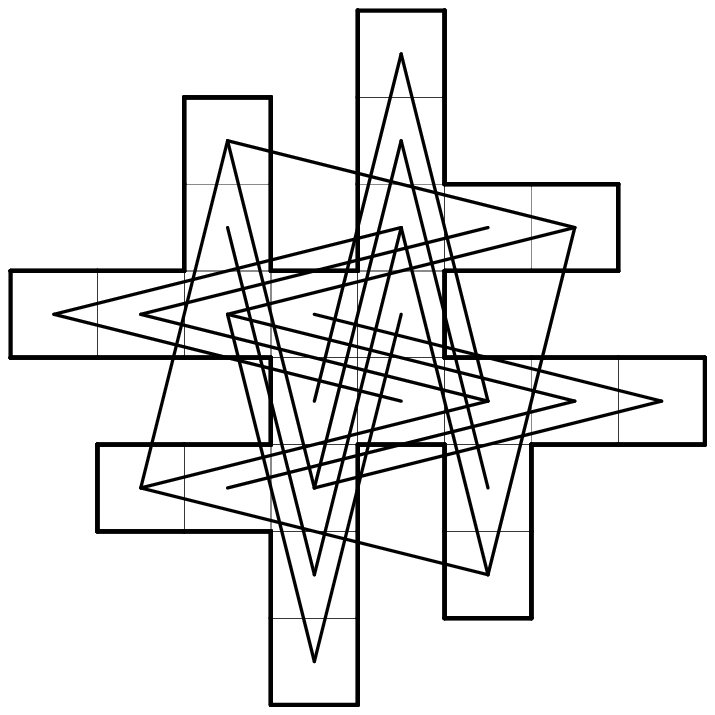}\hspace*{\fill}\newline\bigskip\newline\hspace*{\fill}\includegraphics[scale=0.66]{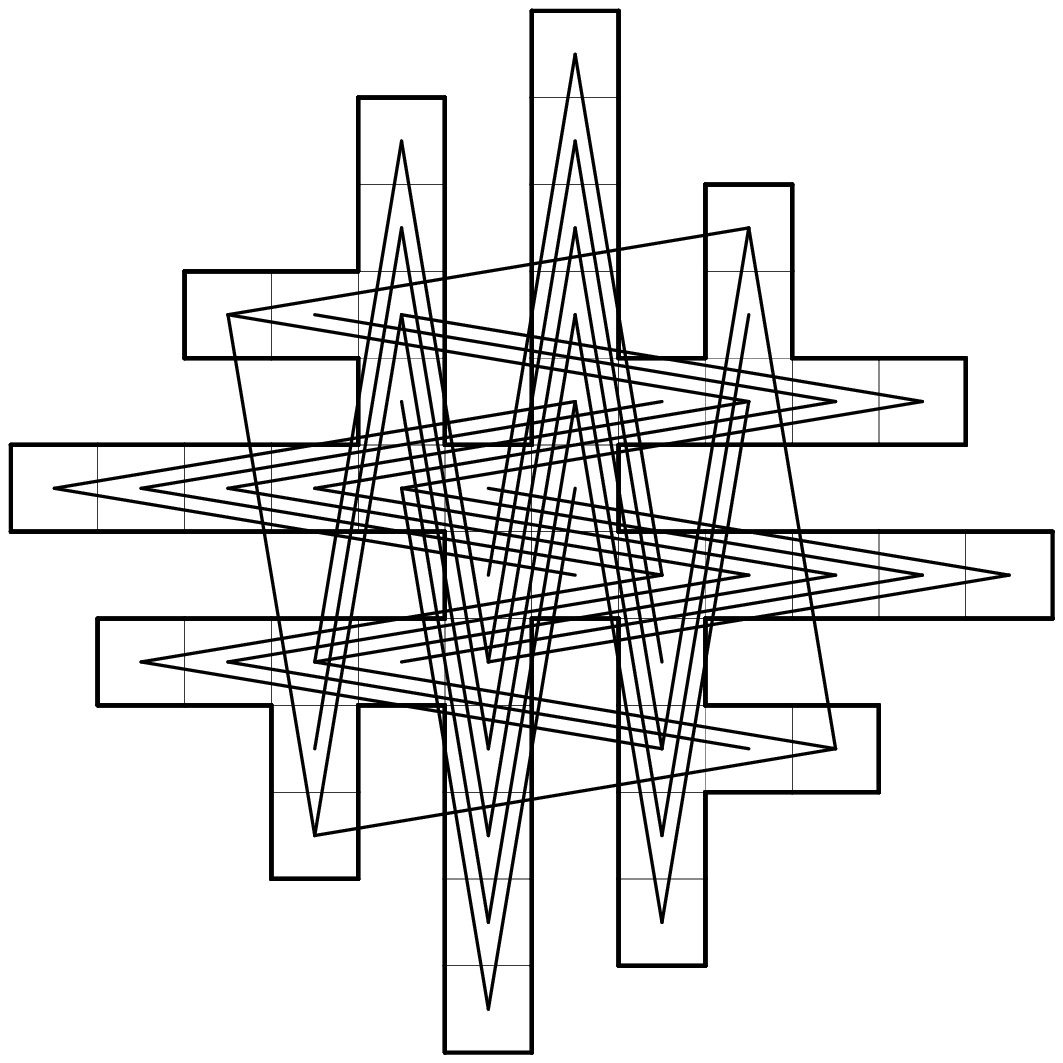}\hspace*{\fill}\caption{}\label{wn01}\end{figure}

Figure \ref{wn01} shows $W_n(0, 1)$, overlaid with the associated leaper graph of a $(1, 2n)$-leaper, for $n = 1$, 2, and 3.

\medskip

\begin{proof} First we show that the leaper graph of $L$ over $W_n(L)$ is connected.

Consider the squares $a^\text{I}$, $a^\text{II}$, $a^\text{III}$, and $a^\text{IV}$ of $W_n(L)$ of centers $(\frac{1}{2}(p + q), \frac{1}{2}(p - q))$, $(\frac{1}{2}(q - p), \frac{1}{2}(p + q))$, $(-\frac{1}{2}(p + q), \frac{1}{2}(q - p))$, and $(\frac{1}{2}(p - q), -\frac{1}{2}(p + q))$. They form a cycle $C$ of $L$.

Let $a$, centered at $(x_a, y_a)$, be any square of $W_n(L)$. It suffices to show that there exists a path of $L$ from $a$ to $C$ within $W_n(L)$.

We consider the case $0 \le x_a \le y_a$ in detail, and all other cases are analogous.

Place $L$ at $a$.

Suppose first that $a$ belongs to $W^\text{I}_n(L)$. Let $L$ move successively in directions \texttt{6}, \texttt{7}, \texttt{6}, \texttt{7}, \ldots, \texttt{6}, \texttt{7} until it cannot advance any further without leaving the board. At that point, $L$ occupies a square $b$ of the form $a^\text{I} + m \cdot (2q, 2q)$ for some nonnegative integer $m$. From $b$, let $L$ move successively in directions \texttt{5}, \texttt{4}, \texttt{6}, \texttt{7}, \texttt{5}, \texttt{4}, \texttt{6}, \texttt{7}, \ldots, \texttt{5}, \texttt{4}, \texttt{6}, \texttt{7} until it arrives at $a^\text{I}$.

Suppose, then, that $a$ belongs to $W^\text{II}_n(L)$. Then a move in direction \texttt{7} brings $L$ to a square $a'$ in $W^\text{I}_n(L)$ and the proof continues as in the previous case.

Suppose, lastly, that $a$ belongs to $W^\text{IV}_n(L)$. This is only possible if $p = 0$ and $a$ has the form $(\frac{1}{2}(4m + 3)q, \frac{1}{2}(4m + 3)q)$ for some nonnegative integer $m$. In that case, let $L$ move successively in directions \texttt{4}, \texttt{6}, \texttt{7}, \texttt{5}, \texttt{4}, \texttt{6}, \texttt{7}, \texttt{5}, \ldots, \texttt{4}, \texttt{6}, \texttt{7}, \texttt{5} until it arrives at $a^\text{IV}$.

We proceed to show that the leaper graphs of $L$ and $M$ over $W_n(L)$ are isomorphic.

Given a square $a$ of $W_n(L)$, define $\varphi(a)$ as follows. Let $i$, out of $\text{I}$, $\text{II}$, $\text{III}$, and $\text{IV}$, be such that $a$ belongs to the wing $W^i_n(L)$. Then $\varphi(a)$ is the reflection of $a$ in the center $O^i$ of $R^i_n(L)$.

\medskip

\begin{lemma} Two squares $a$ and $b$ of $W_n(L)$ are joined by a move of $L$ if and only if $\varphi(a)$ and $\varphi(b)$ are joined by a move of $M$. \label{phi} \end{lemma}

\medskip

\begin{proof} Let $b = a + u$, $u = (x_u, y_u)$, $\varphi(b) = \varphi(a) + v$, and $v = (x_v, y_v)$.

Suppose first that $u$ is an $L$-translation.

Since each of the four translations $(\pm 2q, 0)$ and $(0, \pm 2q)$ is the sum of two $L$-translations, it takes $L$ an even number of moves to travel between two squares in the same net $N^i(L)$, for each $i$ out of $\text{I}$, $\text{II}$, $\text{III}$, and $\text{IV}$. Furthermore, since it takes $L$ two moves to travel between opposite squares in $C$, it also takes $L$ an even number of moves to travel between two squares in opposite nets.

It follows that $a$ and $b$ cannot belong to the same wing or to opposite wings. Therefore, $a$ and $b$ belong to adjacent wings.

By symmetry, it suffices to consider the case when $a$ belongs to $W^\text{I}_n(L)$ and $b$ belongs to $W^\text{II}_n(L)$.

Then \[x_u \equiv \frac{1}{2}(q - p) - \frac{1}{2}(p + q) = -p \pmod{2q}\] and \[y_u \equiv \frac{1}{2}(p + q) - \frac{1}{2}(p - q) = q \pmod{2q}.\]

Consequently, $x_u = -p$ and $y_u = \pm q$.

Let $O^\text{I} + o = O^\text{II}$. Since $a$ and $\varphi(a)$ are symmetric with respect to $O^\text{I}$ and $b$ and $\varphi(b)$ are symmetric with respect to $O^\text{II}$, \[u + v = 2o.\]

Therefore, as $o = (-p - nq, 0)$, \[v = 2o - u = (-p - 2nq, \pm q) = (-s, \pm r)\] and $v$ is an $M$-translation.

Suppose, then, that $v$ is an $M$-translation.

As above, it suffices to consider the case when $a$ belongs to $W^\text{I}_n(L)$ and $b$ belongs to $W^\text{II}_n(L)$, when $x_v \equiv -p$ and $y_v \equiv q$ modulo $2q$.

If $p \neq 0$, then it follows that $x_v = -s$ and $y_v = \pm r$ and the proof continues as before.

When $p = 0$, however, we also need to rule out the possibility that $x_v = s$. This is done as follows.

Let $\varphi(a)$ be at $(x'_a, y'_a)$ and $\varphi(b)$ at $(x'_b, y'_b)$. Since $\varphi(a)$ lies on the boundary or in the interior of $R^\text{I}_n(L)$, $x'_a \ge -\frac{1}{2}(2n - 1)q$. Analogously, $x'_b \le \frac{1}{2}(2n - 1)q$. Therefore, \[x'_b - x'_a \le (2n - 1)q < s.\]

Thus necessarily $x_v = -s$ and $y_v = \pm r$ even if $p = 0$, and the proof continues as before. \end{proof}

\medskip

By Lemma \ref{phi}, the leaper graph of $L$ over $W_n(L)$ is isomorphic to the leaper graph of $M$ over $\varphi(W_n(L))$.

When $n$ is odd, $\varphi(W_n(L))$ coincides with $W_n(L)$, and when $n$ is even, $\varphi(W_n(L))$ is a reflection of $W_n(L)$ (in each of the lines $x = 0$, $y = 0$, and $y = \pm x$).

Consequently, the leaper graph of $M$ over $\varphi(W_n(L))$ is isomorphic to the leaper graph of $M$ over $W_n(L)$ and the proof is complete. \end{proof}

\medskip

Dual pinwheel boards are fundamentally different from all other dual boards we have encountered thus far.

However, they do not uncover any novel pairs of leapers $L$ and $M$ such that there exists a board dual with respect to $L$ and $M$. Since \[A^\text{Pinwheel}_n = A_\texttt{f}^{n - 1}A_\texttt{h},\] all pairs of leapers such that a board dual with respect to them exists by Theorem \ref{pinwheelboard} are already accounted for by Theorem \ref{dual}.

We go on to extract dual direction graphs from dual pinwheel boards.

\medskip

\begin{theorem} Let $n$ be a positive integer and $L$ a skew $(p, q)$-leaper with $p < q$. Then the leaper graph of $L$ over $W_n(L)$ is trivial and the direction graph extracted from it depends only on $n$. We refer to this direction graph as the \emph{pinwheel direction graph} $W_n$ of order $n$.

(The above continues to apply to non-skew leapers $L$ in the sense that the leaper graph of $L$ over $W_n(L)$ is an instantiation of $W_n$.)

The pinwheel direction graph of order $n$ is a dual direction graph of complement $\pi_\text{Pinwheel}(W_n)$ when $n$ is odd, where \[\pi_\text{Pinwheel} = \left(\begin{array}{cccccccc} \texttt{1} & \texttt{2} & \texttt{3} & \texttt{4} & \texttt{5} & \texttt{6} & \texttt{7} & \texttt{8}\\ \texttt{3} & \texttt{8} & \texttt{5} & \texttt{2} & \texttt{7} & \texttt{4} & \texttt{1} & \texttt{6} \end{array}\right)\] is the equivalence permutation induced by the unit $2 \times 2$ matrix and $\left( \begin{smallmatrix} 0 & 1\\ -1 & 0 \end{smallmatrix} \right)$, and $\pi_\text{Reflect} \circ \pi_\text{Pinwheel}(W_n)$ when $n$ is even, and duality matrix $A^\text{Pinwheel}_n$. \label{pinwheelgraph} \end{theorem}

\medskip

\begin{proof} In the setting of the proof of Theorem \ref{pinwheelboard}, let $a$ be any square of $W_n(L)$. Then there exist unique $i$ out of $\text{I}$, $\text{II}$, $\text{III}$, and $\text{IV}$ and integers $k$ and $l$ such that $a = a^i + (2kq, 2lq)$. Assign to $a$ the ordered triplet $\sigma(a) = (i, k, l)$.

Let \begin{align*} A^\text{I} &= \left(\begin{array}{cc} \frac{1}{2} & \frac{1}{2}\\[5pt] \frac{1}{2} & -\;\frac{1}{2} \end{array}\right), & A^\text{II} &= \left(\begin{array}{cc} -\;\frac{1}{2} & \frac{1}{2}\\[5pt] \frac{1}{2} & \frac{1}{2} \end{array}\right),\\ A^\text{III} &= \left(\begin{array}{cc} -\;\frac{1}{2} & -\;\frac{1}{2}\\[5pt] -\;\frac{1}{2} & \frac{1}{2} \end{array}\right), & A^\text{IV} &= \left(\begin{array}{cc} \frac{1}{2} & -\;\frac{1}{2}\\[5pt] -\;\frac{1}{2} & -\;\frac{1}{2} \end{array}\right), \end{align*} and \[A(i, k, l) = A^i + \left(\begin{array}{cc} 0 & 2k\\ 0 & 2l \end{array}\right).\]

We proceed to show that a move of $L$ of direction $j$ leads from $a'$ to $a''$ in $W_n(L)$ if and only if \[A(\sigma(a'')) - A(\sigma(a')) = A_j.\]

The ``if'' part holds as, for every square $a$ of $W_n(L)$, the transpose of the coordinates of $a$ is given by $A(\sigma(a))\left( \begin{smallmatrix} p\\ q \end{smallmatrix} \right)$.

For the ``only if'' part, let a move of $L$ lead from $a'$ to $a''$.

By the proof of Lemma \ref{phi}, $a'$ and $a''$ belong to adjacent wings. By symmetry, it suffices to consider the case when $a'$ belongs to $W^\text{I}_n(L)$ and $a''$ belongs to $W^\text{II}_n(L)$.

By the proof of Lemma \ref{phi}, the move of $L$ from $a'$ to $a''$ is of the form either $(-p, q)$ or $(-p, -q)$. In the former case, it is of direction \texttt{3} and $A(\sigma(a'')) - A(\sigma(a')) = A_\texttt{3}$, and in the latter case it is of direction \texttt{6} and $A(\sigma(a'')) - A(\sigma(a')) = A_\texttt{6}$, as needed.

It follows that the leaper graph of $L$ over $W_n(L)$ is trivial and the direction of a move of $L$ from a square $a'$ to a square $a''$ in $W_n(L)$ depends only on $\sigma(a')$ and $\sigma(a'')$, and not on $p$ and $q$.

(In fact, the leaper graph of $L$ over $W_n(L)$ contains a single simple cycle, namely $C$. However, the argument above continues to apply in generalizations of pinwheel boards replacing $R^i_n(L)$, $i = \text{I}$, $\text{II}$, $\text{III}$, and $\text{IV}$, with other regions, as in Theorem \ref{pinwheelexp}.)

Consider, then, the direction graph $W_n$ of vertices $\sigma(W_n(L))$ such that an arc labeled $i$ points from $\sigma(a')$ to $\sigma(a'')$ if and only if an $L$-move of direction $i$ leads from $a'$ to $a''$. It is extracted from the leaper graph of $L$ over $W_n(L)$ and does not depend on $L$, settling the first part of the theorem.

For the second part of the theorem, we consider the case when $n$ is odd, and the opposite case is analogous.

Define the mapping $\eta$ over the vertices of $W_n$ by \[\eta(i, k, l) = (i,\; \tfrac{1}{2}(n - 1) - k,\; \tfrac{1}{2}(n + 1) - l).\]

By the proof of Lemma \ref{phi}, $W_n$ satisfies the definition of a dual direction graph of complement $\pi_\text{Pinwheel}(W_n)$ and duality matrix $A^\text{Pinwheel}_n$ with the one-to-one mapping $\eta$ between the vertex sets of $W_n$ and $\pi_\text{Pinwheel}(W_n)$. \end{proof}

\medskip

Theorem \ref{pinwheelboard} does not list all values of $n$, $p$, and $q$ such that $W_n(L)$ is dual with respect to $L$ and $M$. For instance, $n = 2$, $p = 2$, and $q = 1$ yield a pinwheel board dual with respect to a $(1, 2)$-leaper and a $(1, 6)$-leaper.

The construction of a dual pinwheel board admits a variety of modifications.

For instance, in the setting of Theorem \ref{pinwheelboard} and its proof, adjoining $a^\text{I} + (2nq, 0)$ and its multiple-of-quarter-turn rotations about $O$ to $W_n(L)$ yields an \emph{augmented} pinwheel board dual with respect to $L$ and $M$.

\begin{figure}[ht!]\hspace*{\fill}\includegraphics[scale=0.66]{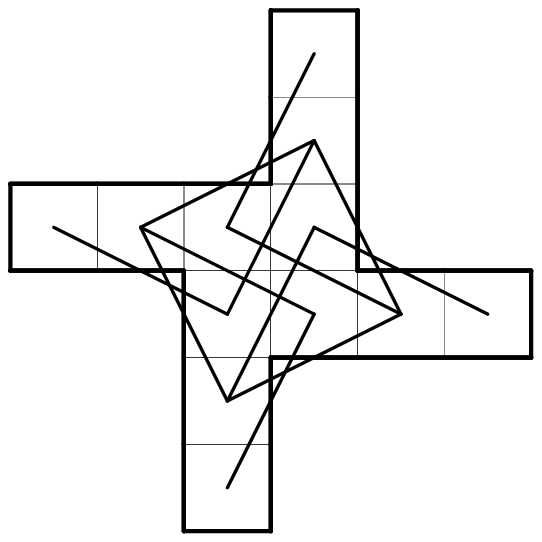}\hspace*{\fill}\caption{}\label{w10112aug}\end{figure}

The first dual board discovered by the author, in October 2005, was the augmented $W_1(0, 1)$ in Figure \ref{w10112aug}.

A different species of modification proceeds by expanding the regions $R^i_n(L)$, $i = \text{I}$, $\text{II}$, $\text{III}$, and $\text{IV}$, in the definition of a pinwheel board.

For instance, in the setting of Theorem \ref{pinwheelboard}, require additionally that $p \neq 0$. Then replacing each of $R^i_n(L)$, $i = \text{I}$, $\text{II}$, $\text{III}$, and $\text{IV}$, with the complete plane yields an infinite board, the disjoint union of $N^i(L)$ over $i = \text{I}$, $\text{II}$, $\text{III}$, and $\text{IV}$, dual with respect to $L$ and $M$.

We give a construction of this type of arbitrarily large finite dual pinwheel-like boards.

\begin{figure}[ht!]\hspace*{\fill}\includegraphics[width=\textwidth]{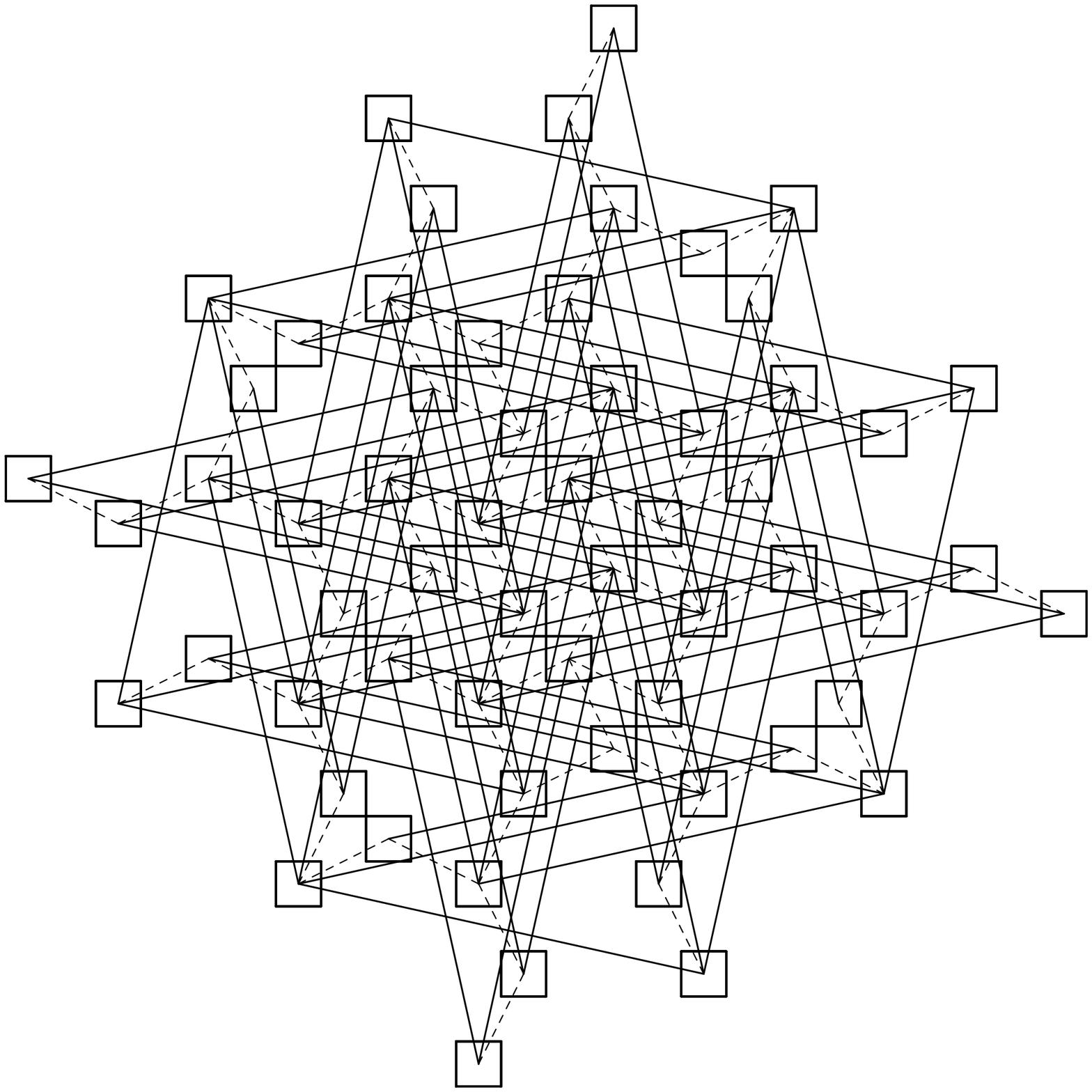}\hspace*{\fill}\caption{}\label{w211229exp}\end{figure}

\medskip

\begin{definition} Let $n$ be a positive integer, $d$ a nonnegative integer, and $L$ a $(p, q)$-leaper with $q \neq 0$.

Let $W^\text{I}_{n, d}(L)$ be the board formed by all squares of $N^\text{I}(L)$ whose centers lie on the boundary or in the interior of the rectangle $R^\text{I}_{n, d}(L)$ bounded by \[p - 2dq \le x + y \le p + 2(n + d)q\] and \[-(2n + 2d - 1)q \le y - x \le (2n + 2d - 1)q.\]

Let $R^i_{n, d}(L)$ and $W^i_{n, d}(L)$, for $i = \text{II}$, $\text{III}$, and $\text{IV}$, be the rotations of $R^\text{I}_{n, d}(L)$ and $W^\text{I}_{n, d}(L)$ by $90^\circ$, $180^\circ$, and $270^\circ$ counterclockwise about $O$.

Then the \emph{expanded} pinwheel board of order $n$ and margin $d$ for $L$, $W_{n, d}(L)$, is the (disjoint if $p < q$) union of its four wings $W^i_{n, d}(L)$, $i = \text{I}$, $\text{II}$, $\text{III}$, and $\text{IV}$. \end{definition}

\medskip

The expanded pinwheel board of order $n$ and margin 0 for $L$ coincides with the pinwheel board of order $n$ for $L$.

\medskip

\begin{theorem} Let $n$ and $d$ be positive integers, $L$ a $(p, q)$-leaper with $p \neq 0$ and $p < q$, and $M$ an $(r, s)$-leaper with \[\left(\begin{array}{c} r\\ s \end{array}\right) = A^\text{Pinwheel}_n \left(\begin{array}{c} p\\ q \end{array}\right).\]

Then the expanded pinwheel board of order $n$ and margin $d$ for $L$ is dual with respect to $L$ and $M$.

Furthermore, the leaper graph of $L$ over $W_{n, d}(L)$ is trivial and the direction graph extracted from it depends only on $n$ and $d$. We refer to this direction graph as the \emph{expanded} pinwheel direction graph $W_{n, d}$ of order $n$ and margin $d$. It is a dual direction graph of complement $\pi_\text{Pinwheel}(W_{n, d})$ when $n$ is odd and $\pi_\text{Reflect} \circ \pi_\text{Pinwheel}(W_{n, d})$ when $n$ is even, and duality matrix $A^\text{Pinwheel}_n$. \label{pinwheelexp} \end{theorem}

\medskip

Figure \ref{w211229exp} shows $W_{2, 1}(1, 2)$ overlaid with the associated leaper graphs of a $(1, 2)$-leaper and a $(2, 9)$-leaper.

The proof of Theorem \ref{pinwheelexp} is analogous to the proofs of Theorems \ref{pinwheelboard} and \ref{pinwheelgraph}.

In particular, we obtain the following corollary.

\medskip

\begin{corollary} Let $n$, $p$, and $q$ be positive integers with $p < q$. Then there exist arbitrarily large boards dual with respect to a $(p, q)$-leaper and a $(q, p + 2nq)$-leaper and arbitrarily large dual direction graphs of duality matrix $A^\text{Pinwheel}_n$.~\label{pinwheelunbound} \end{corollary}

\end{document}